\documentclass[11pt, oneside]{article}   	
\usepackage{geometry}                		
\usepackage{graphicx}
\usepackage{float}
		
\usepackage{amssymb}
\usepackage{amsthm}
\usepackage{amsmath}
\usepackage{mathtools}
\usepackage{enumerate}
\usepackage{mathdots}
\usepackage{extarrows}
\usepackage[dvipsnames]{xcolor}

\usepackage{tikz}
\usetikzlibrary{arrows,decorations.pathmorphing,decorations.markings,backgrounds,positioning,fit,petri}
\usetikzlibrary{calc,intersections,through,backgrounds}
\usepackage{caption}
\usepackage{subcaption}
\usepackage{float}
\usetikzlibrary{calc} 
\usepackage{hyperref}

\usepackage{tikz-cd}
\usepackage{standalone}
\usepackage{titlesec}

\newtheorem{theorem}{Theorem}[section]
\newtheorem{lemma}[theorem]{Lemma}
\newtheorem{proposition}[theorem]{Proposition}
\newtheorem{corollary}[theorem]{Corollary}

\theoremstyle{definition}
\newtheorem{definition}[theorem]{Definition}

\theoremstyle{remark}
\newtheorem{remark}[theorem]{Remark}
\newtheorem*{notation}{Notation}

\newcommand{\A}{\mathbb{A}}

\newcommand{\C}{\mathbb{C}}
\newcommand{\D}{\mathbb{D}}
\newcommand{\Z}{\mathbb{Z}}
\newcommand{\N}{\mathbb{N}}
\newcommand{\CR}{\mathbb{R}}
\newcommand{\Sym}{\text{Sym}}
\newcommand{\f}{\tilde{f}}
\newcommand{\g}{\tilde{g}}

\newcommand{\phii}{\tilde{\varphi}}
\newcommand{\tg}{\tilde{\Gamma}}
\newcommand{\F}{\mathcal{F}}
\newcommand{\W}{\mathcal{W}}

\newcommand{\RR}{\text{R}}
\newcommand{\Xn}{\ensuremath{X_{n}}}
\newcommand{\Xnn}{\ensuremath{X_{n+1}}}
\newcommand{\Yn}{\ensuremath{Y_{n}}}
\newcommand{\Ynn}{\ensuremath{Y_{n+1}}}

\newcommand{\End}{\mathcal{E}}
\newcommand{\G}{\mathcal{G}}
\newcommand{\Gg}{\mathcal{G}'}

\newcommand{\LAMBDA}{\mathbf{\Lambda}}
\newcommand{\Hom}{\text{Hom}}
\newcommand{\Ext}{\text{Ext}}
\newcommand{\fSym}{f} 
\newcommand{\fC}{g} 
\newcommand{\pp}{\textbf{p}} 
\newcommand{\Q}{Q} 
\newcommand{\K}{K} 
\newcommand{\Pp}{P} 
\newcommand{\Cone}{\textnormal{Cone}} 
\newcommand{\tA}{\tilde{A}}

\newcommand{\Iyama}{\Lambda} 
\newcommand{\DIy}{\Delta}
\newcommand{\DAC}{{D}}
\newcommand{\LL}{\Lambda}
\newcommand{\LAC}{L}
\newcommand{\SIy}{\Sigma}
\newcommand{\SSn}{\Sigma_n}

\newcommand{\n}{^{-}}
\newcommand{\p}{^{+}}
\newcommand{\s}{^{\boldsymbol{\cdot}}}

\newcommand{\y}{^{\circ}}
\newcommand{\z}{^{\scriptscriptstyle{\triangle}}}

\newcommand{\perf}{\textnormal{\textbf{perf}}}
\newcommand{\kk}{\boldsymbol{k}}

\DeclarePairedDelimiter\floor{\lfloor}{\rfloor}

\newcommand{\comment}[1]{}
\title{Realising perfect derived categories of Auslander algebras of type A as Fukaya-Seidel categories}
\date{}		
\author{Ilaria Di Dedda}	

\begin{document}
\maketitle

ABSTRACT. We prove that the Fukaya-Seidel categories of a certain family of Lefschetz fibrations on $\C^2$ are equivalent to the perfect derived categories of Auslander algebras of Dynkin type $\A$. We give an explicit equivalence between these categories and the partially wrapped Fukaya categories considered in \cite{DJL}. We provide a complete description of the Milnor fibre of such fibrations.

\tableofcontents

\section{Introduction}\label{section1}
\subsection{Main results}\label{section:mainresults}
Let $f: \C^m \to \C$ be a polynomial with an isolated singularity at the origin. The category we naturally associate to it is the Fukaya-Seidel category $\F(f)$, as constructed in \cite[Chapter~3]{SeidelBk}. The objects of study of this paper are the Fukaya-Seidel categories of a family of polynomial singularities $f_n: \C^2 \to \C$, indexed by natural numbers, which we will explicitly define in Section \ref{section:2}. Due to the natural identification $\Sym^2(\C)\cong \C^2$, $f_n$ have a very simple presentation in terms of polynomials defined on the symmetric product, and are given by the following collection of maps:
\begin{equation*}
    \Sym^2(\C) \to \C \quad \{(x,y)\} \mapsto x^n+y^n.
\end{equation*}

We use constructions on curve singularities developed, independently, by A'Campo (\cite{AC75}, \cite{AC99}) and Gusein-Zade (\cite{GZ}) to compute a favourite collection of generators of such categories, and we prove the following.

\begin{figure}[H]
    \centering
    \begin{minipage}{.45 \textwidth}
        \centering
        \begin{tikzpicture}
            [scale=.8,align=center, v/.style={draw,shape=circle, fill=black, minimum size=1.2mm, inner sep=0pt, outer sep=0pt},
            every path/.style={shorten >= 1mm, shorten <= 1mm},
            font=\small, label distance=1pt,
            every loop/.style={distance=1cm, label=right:}
            ]
            \node[v] (00) at (0,0) {};
            \node[v] (01-) at (0,1) {};
            \node[v] (11) at (1,1) {};
            \node[v] (02) at (0,2) {};
            \node[v] (12+) at (1,2) {};
            \node[v] (22) at (2,2) {};
            \node[v] (03-) at (0,3) {};
            \node[v] (13) at (1,3) {};
            \node[v] (23-) at (2,3) {};
            \node[v] (33) at (3,3) {};

            \node[v] (aa) at (0,5) {};
            \node[v] (bb) at (1,5) {};
            \node[v] (cc) at (2,5) {};
            \node[v] (dd) at (3,5) {};
           
            \node[v] (ee) at (5,5) {};
    
            \path
            (01-) edge[color=black,->] (00)
            (01-) edge[color=black,->] (11)
            (01-) edge[color=black,->] (02)
    
            (02) edge[color=black,->] (12+)
            (11) edge[color=black,->] (12+)
            (22) edge[color=black,->] (12+)
            (13) edge[color=black,->] (12+)
    
            (03-) edge[color=black,->] (02)
            (03-) edge[color=black,->] (13)

            (23-) edge[color=black,->] (13)
            (23-) edge[color=black,->] (22)
            (23-) edge[color=black,->] (33)
            ;

            \path
            (0,3) edge[color=black,->] (0,4)
            (1,3) edge[color=black,->] (1,4)
            (2,3) edge[color=black,->] (2,4)
            (3,3) edge[color=black,->] (3,4)
            ; 

            \node at(0,4.5){$\vdots$};
            \node at(1,4.5){$\vdots$};
            \node at(2,4.5){$\vdots$};
            \node at(3,4.5){$\vdots$};

            \node at(3.5,5){$\dots$};
            \node at(2,4){$\vdots$};
            \node at(4,4){$\iddots$};

            \path
            (aa) edge[color=black,->] (bb)
            (cc) edge[color=black,->] (bb)
            (cc) edge[color=black,->] (dd)
            (4,5) edge[color=black,->] (ee)
            ;
                        
        \end{tikzpicture}
    \end{minipage}
    \begin{minipage}{.45 \textwidth}
        \centering
    \begin{tikzpicture}
        [scale=.8,align=center, v/.style={draw,shape=circle, fill=black, minimum size=1.2mm, inner sep=0pt, outer sep=0pt},
        every path/.style={shorten >= 1mm, shorten <= 1mm},
        font=\small,
        ]
        \node[v] (00) at (0,0) {};
        \node[v] (01-) at (0,1) {};
        \node[v] (11) at (1,1) {};
        \node[v] (02) at (0,2) {};
        \node[v] (12+) at (1,2) {};
        \node[v] (22) at (2,2) {};
        \node[v] (03-) at (0,3) {};
        \node[v] (13) at (1,3) {};
        \node[v] (23-) at (2,3) {};
        \node[v] (33) at (3,3) {};

        \node[v] (aa) at (0,5) {};
        \node[v] (bb) at (1,5) {};
        \node[v] (cc) at (2,5) {};
        \node[v] (dd) at (3,5) {};
       
        \node[v] (ee) at (5,5) {};

        \node (10) at (1,0) {0};
        \node (21) at (2,1) {0};
        \node (32) at (3,2) {0};

        \path
        (00) edge[color=black,->] (01-)
        (01-) edge[color=black,->] (11)
        (01-) edge[color=black,->] (02)

        (02) edge[color=black,->] (12+)
        (11) edge[color=black,->] (12+)
        (12+) edge[color=black,->] (22)
        (12+) edge[color=black,->] (13)

        (02) edge[color=black,->] (03-)
        (03-) edge[color=black,->] (13)

        (13) edge[color=black,->] (23-)
        (22) edge[color=black,->] (23-)
        (23-) edge[color=black,->] (33)

        (00) edge[color=black,->] (10)
        (10) edge[color=black,->] (11)
        (11) edge[color=black,->] (21)
        (21) edge[color=black,->] (22)
        (22) edge[color=black,->] (32)
        (32) edge[color=black,->] (33)
        ;             

        \path
        (0,3) edge[color=black,->] (0,4)
        (1,3) edge[color=black,->] (1,4)
        (2,3) edge[color=black,->] (2,4)
        (3,3) edge[color=black,->] (3,4)
        ; 

        \node at(0,4.5){$\vdots$};
        \node at(1,4.5){$\vdots$};
        \node at(2,4.5){$\vdots$};
        \node at(3,4.5){$\vdots$};

        \node at(3.5,5){$\dots$};
        \node at(2,4){$\vdots$};
        \node at(4,4){$\iddots$};

        \path
        (aa) edge[color=black,->] (bb)
        (bb) edge[color=black,->] (cc)
        (cc) edge[color=black,->] (dd)
        (4,5) edge[color=black,->] (ee)
        ;
    \end{tikzpicture}
    \end{minipage}
    \caption{The two quivers that are central to our results. (left) The quiver whose path algebra is $\tg_n$, where relations are given by commutativity of the squares. (right) The quiver whose path algebra is $\Gamma_n$, with all possible commutativity relations. Both quivers have $n-2$ rows and columns.}
    \label{algebrasgamman}
\end{figure}
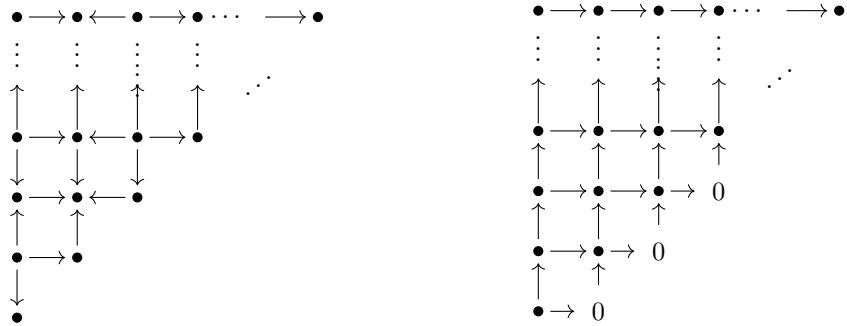

\begin{theorem}[Theorem \ref{thm1intext}]\label{thm1}
    $\F(f_n)$, as a triangulated $A_{\infty}$-category, is quasi-equivalent to the perfect derived category $\perf(\tg_n)$ of the path algebra $\tg_n$ associated to the quiver in Figure \ref{algebrasgamman} (left).
\end{theorem} 

Following \cite{Auroux2}, we equip the Liouville domain $\Sym^2(\D)$ (as in \cite{GPS2}) with a collection of stops $\Lambda^{(2)}_n := \bigcup_{p\in \Lambda_n} \{p\} \times \D$, where $\Lambda_n \subset \partial \D$ is a set of $n$ marked points; we consider the \emph{partially wrapped Fukaya category} $\W(\Sym^2(\D),\Lambda^{(2)}_n)$, as first constructed by Auroux in \cite{Auroux1} for the case at hand and later generalised in \cite{GPS1}. Our second result is the following.

\begin{theorem}[Theorem~\ref{qefn}]\label{thm2}
        $\perf(\tg_n)$ is quasi-equivalent to $\W(\Sym^2(\D),\Lambda^{(2)}_n)$.
    \end{theorem}

\begin{figure}
\centering
    \begin{minipage}{.45 \textwidth}
    \centering
    \begin{tikzpicture}[scale=.6,align=center, v/.style={draw,shape=circle, fill=black, minimum size=1.2mm, inner sep=0pt, outer sep=0pt}, font=\small,
        ]
    \draw (0,0) circle (3cm);

    \node[v] at(90:3cm){};
    \node[v] at(30:3cm){};
    \node[v] at(-30:3cm){};
    \node[v] at(-90:3cm){};
    \node[v] at(150:3cm){};
    \node[v] at(-150:3cm){};

    \draw[color=blue] (70:3cm)[inner sep=0pt] to[bend left=35] (110:3cm);
    \draw[color=blue] (115:3cm)[inner sep=0pt] to[bend right=20] (10:3cm);
    \draw[color=blue] (180:3cm)[inner sep=0pt] to[bend right=0] (0:3cm);
    \draw[color=blue] (185:3cm)[inner sep=0pt] to[bend left=20] (-50:3cm);
    \draw[color=blue] (-110:3cm)[inner sep=0pt] to[bend left=35] (-70:3cm);
    \end{tikzpicture}
    \end{minipage}
    \qquad
    \begin{minipage}{.45 \textwidth}
        \centering
        \begin{tikzpicture}[scale=.6,align=center, v/.style={draw,shape=circle, fill=black, minimum size=1.2mm, inner sep=0pt, outer sep=0pt}, font=\small,
            ]
    \draw (0,0) circle (3cm);

    \node[v] at(90:3cm){};
    \node[v] at(30:3cm){};
    \node[v] at(-30:3cm){};
    \node[v] at(-90:3cm){};
    \node[v] at(150:3cm){};
    \node[v] at(-150:3cm){};

     \draw[color=blue] (70:3cm)[inner sep=0pt] to[bend left=35] (110:3cm);
     \draw[color=blue] (65:3cm)[inner sep=0pt] to[bend left=20] (180:3cm);
     \draw[color=blue] (60:3cm)[inner sep=0pt] to[bend right=0] (240:3cm);
     \draw[color=blue] (55:3cm)[inner sep=0pt] to[bend right=20] (-70:3cm);
     \draw[color=blue] (50:3cm)[inner sep=0pt] to[bend right=35] (10:3cm);
    \end{tikzpicture}
    \end{minipage}
\caption{Two collections of arcs on $\D$ giving rise to two different collections of generators of the partially wrapped Fukaya category $\W(\Sym^2(\D),\Lambda^{(2)}_n)$, here for $n=6$.}
\label{fig:introgensdisk}
\end{figure}
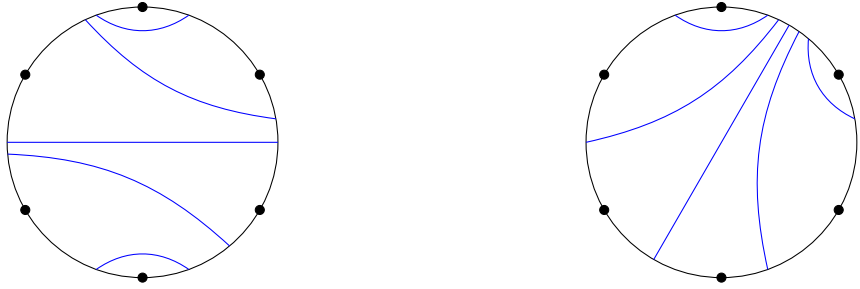

Theorem \ref{thm2} amounts to an explicit computation of the derived endomorphism algebra $\mathcal{B}_n$ of a favourite set of generators of $\W(\Sym^2(\D),\Lambda^{(2)}_n)$, given by pairs of arcs on $\D$ as in Figure \ref{fig:introgensdisk} (left). Our computations provide an isomorphism between $\mathcal{B}_n$ and the algebra $\tg_n$ given by Theorem \ref{thm1}; this will be the object of discussion of Section \ref{section:W2n}.

The main result in \cite{DJL} provides a quasi-equivalence of triangulated $A_{\infty}$-categories between $\W(\Sym^2(\D),\Lambda^{(2)}_n)$ and the derived $A_{\infty}$-category of the algebra $\Gamma_{n-2}$; we construct the latter as follows. Consider the path algebra of the linearly oriented $A_n$-quiver, which is known to have ${n+1} \choose 2$ isomorphism classes of indecomposable modules (see, for example, \cite[Theorem~3.18]{Kirillov}). Let $\{M_i\}$ for $i=1, \dots, {{n+1}\choose{2}}$ be a complete set of such representatives, and define $\Gamma_n$ to be the endomorphism algebra of a such a collection:
\begin{equation}\label{eq:endomalgebra}
    \Gamma_n := \bigoplus_{i,j} \text{Hom}(M_i,M_j)
\end{equation}
where the sum ranges over all such elements. Equivalently (\cite[Section~VII.1]{ARS}), the $\kk$-algebra ($\kk$ any field) $\Gamma_n$ is the path algebra of the quiver depicted in Figure \ref{algebrasgamman} (right), with all possible commutativity relations. We call the collection of algebras $\Gamma_n$ \emph{Auslander algebras of Dynkin type $\A$}, and we refer to Section \ref{section:background} for the proper contextualisation of these structures. We consider the perfect derived category $\perf(\Gamma_n)$ of $\Gamma_n$, whose objects are bounded complexes $\dots~\to~M_0~\to~M_1~\to~\dots$ of projective $\Gamma_n$-modules. Of this category, we have a natural exceptional collection of generators, which we will call \emph{Iyama generators} (due to Iyama's work on \emph{higher} Auslander algebras, see Sections \ref{section:background} and \ref{section:future}), given by the collection of indecomposable $A_n$-modules.

Together with \cite[Theorem~1]{DJL}, which provides the derived equivalence of $\perf(\Gamma_n)$ and $\W(\Sym^2(\D),\Lambda^{(2)}_n)$, Theorem \ref{thm2} allows us to complete the diagram of quasi-equivalences:

\begin{equation*}
  \begin{tikzcd}
    \perf(\tg_n) \arrow{r}{\simeq} \arrow[swap,dashed]{dr}{\simeq} & \W(\Sym^2(\D),\Lambda^{(2)}_n) \arrow{d}{\simeq} \\
     & \perf(\Gamma_n)
  \end{tikzcd}
\end{equation*}

The derived equivalence of the algebras $\tg_n$ and $\Gamma_n$ can be interpreted as the higher dimensional version of the following well-known fact. Consider the $A_n$ quiver with alternating and linear orientations of the arrows, denoted as $\tA_n$ and $A_n$ respectively (Figure \ref{fig:orientationAn}). $\tA_n$ arises more naturally from the perspective of singularity theory (\cite[Section~2, Exemple~1]{AC75}), and is the most natural presentation of the quiver associated to singularity of the real polynomial in one variable $x^{n+1}$ (the so-called $A_n$ singularity). On the other hand, the linearly oriented quiver $A_n$ arises more naturally in relation to representation theory (this is known as the \emph{standard} presentation of the quiver). The path algebras of these quivers are well-known to be derived equivalent (see, for example, \cite[Theorem~3.2]{KY11} and \cite[Theorem~8.6]{FZ03}). Furthermore, they are also known (\cite{Auroux1} and \cite{HKK17}) to be derived equivalent to the partially wrapped Fukaya category $\W(\D, \Lambda_{n+1})$; the (single) arcs depicted in Figure \ref{fig:introgensdisk} represent two generating collections, whose endomorphism algebras are isomorphic to the path algebras of $\tA_n$ and $A_n$ respectively. We provide an explicit computation of the derived equivalence of $\tg_n$ and $\Gamma_n$, which amounts to relating the two collections of generators depicted in Figure \ref{fig:introgensdisk} to each other (Proposition \ref{prop:equivperf}). 

\begin{figure}
    \centering
    \begin{minipage}{.45 \textwidth}
        \centering
        \begin{tikzpicture}
            [scale=.8,align=center, v/.style={draw,shape=circle, fill=black, minimum size=1.2mm, inner sep=0pt, outer sep=0pt},
            every path/.style={shorten >= 1mm, shorten <= 1mm},
            font=\small,
            ]

            \node at (-1,0) {$\tA_n:$};
       
            \node[v] (0) at (0,0) {};
            \node[v] (1) at (1,0) {};
            \node[v] (2) at (2,0) {};
            \node[v] (3) at (3,0) {};
           
            \node[v] (5) at (5,0) {};
            \node[v] (6) at (6,0) {};
 
            \path
            (0) edge[color=black,->] (1)
            (2) edge[color=black,->] (1)
            (2) edge[color=black,->] (3)

            (5) edge[color=black,->] (4,0)
            (5) edge[color=black,->] (6)
            ;

            \node at (3.65,0) {$\dots$};
                        
        \end{tikzpicture}
    \end{minipage}
    \quad
    \begin{minipage}{.45 \textwidth}
        \centering
        \begin{tikzpicture}
            [scale=.8,align=center, v/.style={draw,shape=circle, fill=black, minimum size=1.2mm, inner sep=0pt, outer sep=0pt},
            every path/.style={shorten >= 1mm, shorten <= 1mm},
            font=\small,
            ]

            \node at (0,0) {$A_n:$};
       
            \node[v] (1) at (1,0) {};
            \node[v] (2) at (2,0) {};
            \node[v] (3) at (3,0) {};
           
            \node[v] (5) at (5,0) {};
            \node[v] (6) at (6,0) {};
 
            \path
            (1) edge[color=black,->] (2)
            (2) edge[color=black,->] (3)
            (5) edge[color=black,->] (6)
            (4,0) edge[color=black,->] (5)
            ;

            \node at (3.65,0) {$\dots$}; 
        \end{tikzpicture}
    \end{minipage}
\caption{The $A_n$ quiver, with alternating and linear orientations of the arrows.}
\label{fig:orientationAn}
\end{figure}

Immediate consequences of Theorem \ref{thm1}, Theorem \ref{thm2} and \cite[Theorem~1]{DJL} are the following quasi-equivalences of $A_{\infty}$-categories.
\begin{corollary}\label{cor:maineq}
    $\F(f_n) \simeq \perf(\Gamma_n) \simeq \W(\Sym^2(\D),\Lambda^{(2)}_n)$.
\end{corollary}

Corollary \ref{cor:maineq} allows us to construct a restriction functor
\begin{equation*}
    \perf(\Gamma_{n}) \to \F(\Sigma_n)
\end{equation*}
to the compact Fukaya category of the regular fibre of $f_n$, which is a punctured surface of genus $\frac{(n-2)^2}{4}$ (resp. $\frac{(n-1)(n-3)}{4}$) for $n$ even (resp. odd) and $\floor{\frac{n+1}{2}}$ punctures, equipped with a canonical grading coming from $\C^2$. We construct such functor in Section \ref{section:Fcatregfibre}, where we also prove that it is essentially surjective. 

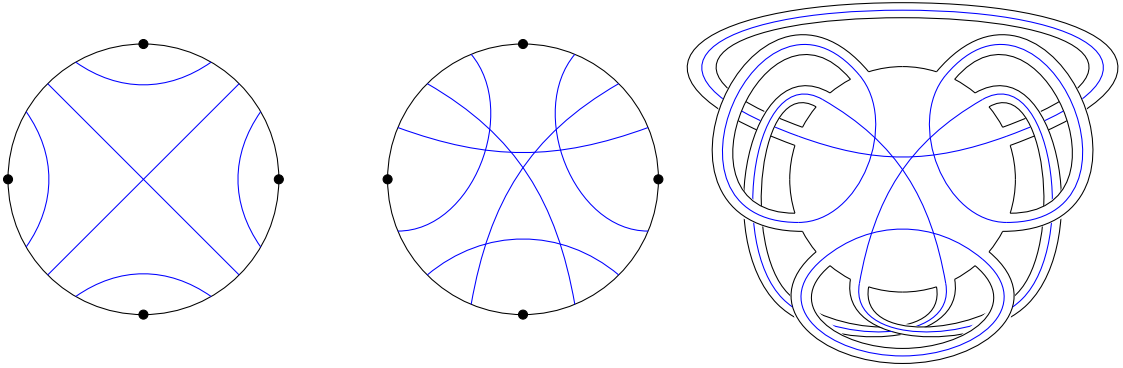
\begin{figure}
    \centering
    \begin{minipage}{.3 \textwidth}
        \centering
        \begin{tikzpicture}[scale=.6,align=center, v/.style={draw,shape=circle, fill=black, minimum size=1.2mm, inner sep=0pt, outer sep=0pt}, font=\small,
            ]
        \draw (0,0) circle (3cm);
    
        \node[v] at(90:3cm){};
        \node[v] at(0:3cm){};
        \node[v] at(180:3cm){};
        \node[v] at(270:3cm){};
    
        \draw[color=blue] (60:3cm)[inner sep=0pt] to[bend left=35] (120:3cm);
        \draw[color=blue] (-30:3cm)[inner sep=0pt] to[bend left=35] (30:3cm);
        \draw[color=blue] (-120:3cm)[inner sep=0pt] to[bend left=35] (-60:3cm);
        \draw[color=blue] (150:3cm)[inner sep=0pt] to[bend left=35] (210:3cm);

        \draw[color=blue] (45:3cm)[inner sep=0pt] to (225:3cm);
        \draw[color=blue] (-45:3cm)[inner sep=0pt] to (-225:3cm);
        \end{tikzpicture}
        \end{minipage}
        \quad
        \begin{minipage}{.3 \textwidth}
            \centering
            \begin{tikzpicture}[scale=.6,align=center, v/.style={draw,shape=circle, fill=black, minimum size=1.2mm, inner sep=0pt, outer sep=0pt}, font=\small,
                ]
        \draw (0,0) circle (3cm);
    
        \node[v] at(90:3cm){};
        \node[v] at(0:3cm){};
        \node[v] at(180:3cm){};
        \node[v] at(270:3cm){};
    
         \draw[color=blue] (67.5:3cm) to[out=230,in=180] (-22.5:3cm);
         \draw[color=blue] (112.5:3cm) to[out=-50,in=0] (202.5:3cm);
        \draw[color=blue] (22.5:3cm) to[out=200,in=-20] (157.5:3cm);
        \draw[color=blue] (45:3cm) to[out=210,in=80] (247.5:3cm);
        \draw[color=blue] (135:3cm) to[out=-30,in=100] (-67.5:3cm);
        \draw[color=blue] (-45:3cm) to[out=140,in=40] (225:3cm);

        \end{tikzpicture}
        \end{minipage}
        \quad
        \begin{minipage}{.3 \textwidth}
        \centering
        \begin{tikzpicture}[remember picture,overlay,scale=.5,align=center, v/.style={draw,shape=circle, fill=black, minimum size=1.2mm, inner sep=0pt, outer sep=0pt}, font=\small,]
        \tikzstyle{reverseclip}=[insert path={(current page.north east) --
        (current page.south east) -- (current page.south west) -- (current page.north west) -- (current page.north east)}]
            \draw (0:3cm) arc (0:17.5:3cm);
            \draw (27.5:3cm) arc (27.5:40:3cm);
            \draw (50:3cm) arc (50:62.5:3cm);
            \draw (72.5:3cm) arc (72.5:90:3cm);
            \draw (90:3cm) arc (90:107.5:3cm);
            \draw (117.5:3cm) arc (117.5:130:3cm);
            \draw (140:3cm) arc (140:152.5:3cm);
            \draw (162.5:3cm) arc (162.5:180:3cm);
            \draw (180:3cm) arc (180:197.5:3cm);
            \draw (207.5:3cm) arc (207.5:220:3cm);
            \draw (230:3cm) arc (230:242.5:3cm);
            \draw (252.5:3cm) arc (252.5:270:3cm);
            \draw (270:3cm) arc (270:287.5:3cm);
            \draw (297.5:3cm) arc (297.5:310:3cm);
            \draw (320:3cm) arc (320:332.5:3cm);
            \draw (342.5:3cm) arc (342.5:360:3cm);

    \draw[color=blue] (67.5:3cm) to[out=230,in=180] (-22.5:3cm);
    \draw[color=blue] (67.5:3cm) to[out=50,in=360, looseness=2.5] (-22.5:3cm);
    \draw (72.5:3cm) to[out=50,in=360, looseness=2.7] (-27.5:3cm);
    \draw (62.5:3cm) to[out=50,in=360, looseness=2.3] (-17.5:3cm);

    \draw[color=blue] (112.5:3cm) to[out=-50,in=0] (202.5:3cm);
    \draw[color=blue]  (112.5:3cm) to[out=130,in=180, looseness=2.5] (202.5:3cm);
    \draw (107.5:3cm) to[out=130,in=180, looseness=2.7] (207.5:3cm);
    \draw (117.5:3cm) to[out=130,in=180, looseness=2.3] (197.5:3cm);

    \draw[color=blue] (-45:3cm) to[out=140,in=40] (225:3cm);
    \draw[color=blue]  (-45:3cm) to[out=320,in=0,looseness=1.5] (-90:4.7cm)  to[out=180,in=220, looseness=1.5] (225:3cm);
    \draw(-40:3cm) to[out=320,in=0,looseness=1.5] (-90:4.9cm) to[out=180,in=220, looseness=1.5]  (220:3cm);
    \draw(-50:3cm)  to[out=320,in=0,looseness=1.5] (-90:4.5cm)  to[out=180,in=220, looseness=1.5] (230:3cm);

    \draw[color=blue] (22.5:3cm) to[out=200,in=-20] (157.5:3cm);
    \draw[color=blue] (45:3cm) to[out=210,in=80] (247.5:3cm);
    \draw[color=blue] (135:3cm) to[out=-30,in=100] (-67.5:3cm);

    \begin{pgfinterruptboundingbox}
    \path [clip]  (68:3.8cm) --(71:3.35cm) -- (60:3.55cm)  -- (58.5:3.8cm) -- (49:4cm) -- (35:4.5cm)-- (27:4.65cm) --  (20:4.65cm)-- (-1:4.3cm) --(-5:4.2cm)--  (-11:3.8cm)-- (-19:3.9cm)-- (-14:4.35cm)--  (-10:4.7cm) -- (0:4.8cm)--  (21:5.25cm) -- (28:5.25cm)-- (42:4.4cm) -- (50:4.4cm) --   (60:4.1cm) -- cycle[reverseclip];
    \end{pgfinterruptboundingbox}

    \begin{pgfinterruptboundingbox}
    \path [clip]  (180-60:4.1cm)--(180-50:4.4cm) --(180-42:4.4cm) -- (180-28:5.25cm)-- (180-21:5.25cm) -- (180-0:4.8cm)--  (180--10:4.7cm)-- (180--14:4.35cm)--(180--19:3.9cm) --  (180--11:3.8cm) --(180--5:4.2cm) -- (180--1:4.3cm)-- (180-20:4.65cm)--(180-27:4.65cm) -- (180-35:4.5cm)-- (180-49:4cm) -- (180-58.5:3.8cm) -- (180-60:3.55cm)-- (180-71:3.35cm) --  (180-68:3.8cm)  -- cycle[reverseclip];
        \end{pgfinterruptboundingbox}

    \begin{pgfinterruptboundingbox}
    \path [clip] (305:4.15cm) -- (290:4.45cm) -- (300:4.7cm) --(307:4.75cm)-- (315:4.25cm) -- cycle[reverseclip];
    \end{pgfinterruptboundingbox}

    \begin{pgfinterruptboundingbox}
    \path [clip](180-315:4.25cm) --(180-307:4.75cm)-- (180-300:4.7cm)  --  (180-290:4.45cm)-- (180-305:4.15cm)-- cycle[reverseclip];
    \end{pgfinterruptboundingbox}
                
    \draw[color=blue] (45:3cm) to[out=30,in=45] (-45:4.5cm) to[out=225,in=260] (247.5:3cm);
    \draw(50:3cm) to[out=30,in=45,looseness=1.1] (-45:4.7cm) to[out=225,in=260] (242.5:3cm);
    \draw(40:3cm) to[out=30,in=45,looseness=.9] (-45:4.3cm) to[out=225,in=260] (252.5:3cm);

    \begin{pgfinterruptboundingbox}
    \path [clip] (17.5:3.7cm) -- (18:4.1cm) -- (27:4.1cm) -- (27:3.65cm)  -- cycle [reverseclip];
    \end{pgfinterruptboundingbox}

    \begin{pgfinterruptboundingbox}
    \path [clip] (270:4.15cm)--  (275:4.15cm) -- (279.5:4cm) -- (265:3.8cm) --  (255:3.9cm)-- cycle[reverseclip];
    \end{pgfinterruptboundingbox}

    \draw[color=blue] (135:3cm) to[out=150,in=135] (225:4.5cm) to[out=-45,in=280] (-67.5:3cm);
    \draw(130:3cm) to[out=150,in=135,looseness=1.1] (225:4.7cm) to[out=-45,in=280] (-62.5:3cm);
    \draw(140:3cm) to[out=150,in=135,looseness=0.9] (225:4.3cm) to[out=-45,in=280] (-72.5:3cm);

    \begin{pgfinterruptboundingbox}
    \path [clip] (180-27:3.65cm) --  (180-27:4.1cm) --  (180-18:4.1cm) --  (180-17.5:3.7cm) -- cycle [reverseclip];
    \end{pgfinterruptboundingbox}

    \draw[color=blue] (22.5:3cm) to[out=20,in=0,looseness=3] (90:4.5cm) to[out=180,in=160,looseness=3] (157.5:3cm);

    \draw(17.5:3cm)  to[out=20,in=0,looseness=3] (90:4.7cm) to[out=180,in=160,looseness=3] (162.5:3cm);
    \draw(27.5:3cm)  to[out=20,in=0,looseness=3] (90:4.3cm)  to[out=180,in=160,looseness=3] (152.5:3cm);
        
    \end{tikzpicture}
    \end{minipage}
    \vskip .5cm 
    \caption{\emph{A posteriori} construction of the Milnor fibre $\Sigma_n$, here for $n=4$.}
    \label{fig:allarcs}
\end{figure}

The description of the Milnor fibre and vanishing cycles given in Section \ref{section:vcyclesalgorithm} allows us to provide a geometric motivation for the algebraic definition (\ref{eq:endomalgebra}). Following from our computations, we \emph{a posteriori} obtain the following description of $\Sigma_n$. We consider $(\D,\Lambda_{n})$, the disk equipped with $n$ stops, and all possible (isotopy classes of) arcs with endpoints on $\partial \D \setminus \Lambda_n$ (not necessarily pairwise disjoint, see Figure \ref{fig:allarcs}, left), here denoted as $\lambda_i$. \cite[Theorem~4.3]{HKK17}, together with \cite[Theorem~1]{DJL}, states that such arcs represent all irreducible $A_{n-1}$-modules. The disk model captures morphisms (computed in the derived category of finitely generated $A_{n-1}$-modules) in all degrees, as given by either intersection points between arcs or Reeb chords obtained by ``flowing'' counter-clockwise along the boundary of $\D$. As we are interested in recovering only the endomorphism algebras of irreducible modules (i.e.\ morphisms in degree zero), we perform small positive perturbations (in the sense of \cite[Definition~7]{Auroux1}) on the arcs $\lambda_i$, for each of the ordered pairs $\lambda_i<\lambda_j$ (order which is detailed in Section \ref{section:vcyclesalgorithm}). After such perturbations, $\{\lambda_i\}$ are as in Figure \ref{fig:allarcs} (centre), and morphisms are generated by intersection points \cite[Definition~8]{Auroux1}. $\hat{\Sigma}_n$ is then constructed by attaching 1-handles along $\partial\D$, so that arcs become closed circles, intersecting transversely only away from the handles (Figure \ref{fig:allarcs}, right). Handle attachment can be done symplectically (\cite{Weinstein}) and so that $\hat{\Sigma}_n$ is orientable, and the surface with boundary we obtain can be completed to ${\Sigma}_n=\hat{\Sigma}_n \cup \left(\partial \hat{\Sigma}_n \times [1,\infty)\right)$, which is a punctured surface.

\subsection{Background}\label{section:background}
\subsubsection{Fukaya-Seidel categories}
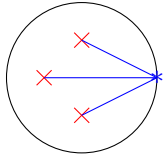
\begin{figure}
   \centering
   \begin{tikzpicture}
        \draw (0,0) circle (1);

        \draw[color=red] (.1,.6) -- (-.1,.4);
        \draw[color=red] (-.1,.6) -- (.1,.4);

        \draw[color=red] (-.4,.1) -- (-.6,-.1);
        \draw[color=red] (-.6,.1) -- (-.4,-.1);

        \draw[color=red] (.1,-.6) -- (-.1,-.4);
        \draw[color=red] (-.1,-.6) -- (.1,-.4);

        \node[color=blue] (star) at (1,0) {$*$};

        \draw[color=blue] (1,0) -- (0,.5);
        \draw[color=blue] (1,0) -- (-.5,0);
        \draw[color=blue] (1,0) -- (0,-.5);
    \end{tikzpicture}
    \caption{Base of a Lefschetz fibration, with regular value taken to be far away from the critical ones.}
    \label{fig:diagramvpaths}
\end{figure}

Picard-Lefschetz theory can be formulated in various contexts and degrees of generality. As per \cite{SeidelBk}, \cite{SeidelArt1}, \cite{SeidelArt2}, we equip the open, exact symplectic manifold $\C^2$ with standard symplectic form $\omega$ and complex structure $J$. We consider a \emph{symplectic Lefschetz fibration} $f: \C^2 \to \C$, i.e.\ a Lefschetz fibration compatible with $\omega$. Fixing a regular value $*$ and corresponding regular fibre $F_*:=f^{-1}(*)$, we consider a distinguished collection of embedded paths $\gamma_i: [0,1] \to \C$ indexed by the critical values (called \emph{vanishing paths}) on the base, such that $F_*$ is the fibre above $\gamma_i(0)$ for all $i$, $\{\gamma_i(1)\}_i$ is the set of critical values of $f$, and the paths are pairwise disjoint away from $\gamma_i(0)$ (Figure~\ref{fig:diagramvpaths}). For $t$ close to critical values of $f$, we construct \emph{vanishing cycles} as embedded curves (half dimensional Lagrangian submanifolds) on $F_t=f^{-1}(t)$, such that these collapse to a point as we approach each critical fibre $F_{\gamma_i(1)}=f^{-1}(\gamma_i(1))$. Using symplectic parallel transport in $\C^2$, we can transport all vanishing cycles along their respective vanishing paths so that they all lie in the fibre $F_*$. To each vanishing path we can associate a \emph{Lefschetz thimble} $\Delta_{\gamma_i}$ as the union of all corresponding vanishing cycles above that path: these are embedded Lagrangian disks in $\C^2$, whose boundaries $\partial \Delta_{\gamma_i} = \Delta_{\gamma_i} \cap F_*$ are the vanishing cycles $V_i \subset F_*$. Vanishing cycles and thimbles are naturally ordered by the (clockwise) ordering of the vanishing paths, as given by the clockwise ordering of the angles at the common intersection point $\gamma_i(0)$ (see Figure~\ref{fig:diagramvpaths}). We fix an indexing of the cycles such that $V_i < V_j$ whenever $i < j$.

As detailed in \cite[Chapter~3]{SeidelBk}, one can define the Fukaya-Seidel category $\F(f)$ associated to a symplectic Lefschetz fibration $f$. Objects of this category are \emph{Lagrangian branes} $\Delta^{\#}_{\gamma_i}$, consisting of Lefschetz thimbles $\Delta_{\gamma_i}$ equipped with additional spin structures and gradings (these two are what constitutes a so-called \emph{brane structure} on a Lagrangian submanifold). Morphisms between Lagrangian branes are given by Floer complexes. Moreover, $\F(f)$ is generated (as a triangulated category) by a distinguished collection of Lagrangian branes, and is independent (up to quasi-equivalence) of the choices of vanishing paths after taking twisted complexes (\cite[Section~18j]{SeidelBk}). We can also equip each vanishing cycle $V_{i}$ with a brane structure $V^{\#}_i$ induced by that of its Lefschetz thimble; this turns each vanishing cycle into an object of the (compact) Fukaya category $\F(F_*)$ of the regular fibre of $f$. Essentially by construction (due to Seidel \cite[Section~18e]{SeidelBk}), we have an isomorphism of Floer complexes:
\begin{equation*}
    CF^*_{\F(f)}(\Delta^{\#}_{\gamma_i},\Delta^{\#}_{\gamma_j}) \cong  CF^*_{\F(F_*)}(V^{\#}_i,V^{\#}_j)
\end{equation*}
whenever $i<j$, while each morphism space $\Hom_{\F(f)}(\Delta^{\#}_{\gamma_i},\Delta^{\#}_{\gamma_j})$ vanishes whenever $i>j$ and is one-dimensional for $i=j$. Concretely, this allows us to carry out our Floer cohomology computations in the directed $A_{\infty}$-subcategory of $\F(F_*)$ associated to an ordered collection of vanishing cycles.

The (derived) Fukaya-Seidel category and directed category of cycles are invariants of the Lefschetz fibration: the independence, up to quasi-equivalence, of these categories on the choice of vanishing paths (\cite[Sections~16, 18]{SeidelBk}) is guaranteed by the fact that two different choices of (clockwise ordered) collections of paths can be related to each other through a series of \emph{mutations} (see, for example, \cite[Lemma~2.23]{Keating}). The existence of such mutations relies on the existence of a simply-transitive action of the braid group on the set of all distinguished collections of vanishing paths. In the remainder of this section, we will review what the geometric effect of this action on a given collection is, for which we refer to Seidel's book; this will not come into use until Section \ref{section:vcyclesalgorithm}.

\begin{figure}
    \centering
    \begin{minipage}{.3 \textwidth}
        \centering
        \begin{tikzpicture}[v/.style={draw,shape=circle, fill=black, minimum size=1.2mm, inner sep=0pt, outer sep=0pt}]
            \node[v] at (0,0) {};
            \node[v] at (-.5,1) {};
            \node[v] at (.5,1) {};

            \draw (0,0) to (-.5,1);
            \draw (0,0) to (.5,1);

            \node[label=180:$\gamma_{i-1}$] at (-.25,.5) {};
            \node[label=0:$\gamma_{i}$] at (.25,.5) {};
        \end{tikzpicture}
    \end{minipage}
    \begin{minipage}{.3 \textwidth}
        \centering
        \begin{tikzpicture}[v/.style={draw,shape=circle, fill=black, minimum size=1.2mm, inner sep=0pt, outer sep=0pt}]
            \node[v] at (0,0) {};
            \node[v] at (-.5,1) {};
            \node[v] at (.5,1) {};

            \draw (0,0) to (-.5,1);
            \draw (0,0) to[out=150,in=225] (-1,1.5) to[out=45] (.5,1);

            \node[label=0:$\gamma_{i-1}$] at (-.5,.5) {};
            \node[label=45:$\tau_{\gamma_{i-1}}(\gamma_{i})$] at (.25,1) {};
        \end{tikzpicture}
    \end{minipage}
    \begin{minipage}{.3 \textwidth}
        \centering
        \begin{tikzpicture}[v/.style={draw,shape=circle, fill=black, minimum size=1.2mm, inner sep=0pt, outer sep=0pt}]
            \node[v] at (0,0) {};
            \node[v] at (-.5,1) {};
            \node[v] at (.5,1) {};

            \draw (0,0) to[out=30,in=-45] (1,1.5) to[out=135,in=45] (-.5,1);
            \draw (0,0) to (.5,1);

            \node[label=135:$\tau^{-1}_{\gamma_{i}}(\gamma_{i-1})$] at (-.25,1) {};
            \node[label=180:$\gamma_{i}$] at (.25,.5) {};
        \end{tikzpicture}
    \end{minipage}
    \caption{Hurwitz moves on the vanishing paths.}
    \label{fig:Hurwitz}
\end{figure}
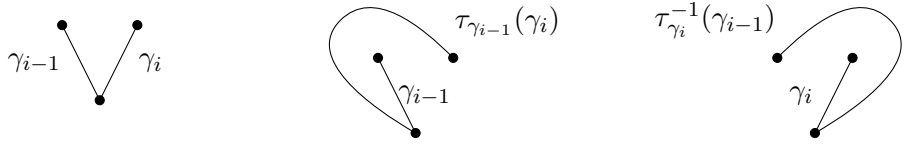

Let $V_0$ and $V_1$ be two vanishing cycles, equipped with a brane structure $\{V_i^{\#}\}_{i=0,1}$, associated to a choice of vanishing paths $\gamma_0$ and $\gamma_1$. It is a non-trivial result by Seidel (\cite[Corollary~17.17]{SeidelBk}) that, in the derived Fukaya category of the regular fibre,
\begin{equation*}
    T_{V_0^{\#}}(V_1^{\#})\cong \tau_{V_0^{\#}}(V_1^{\#}),
\end{equation*}
where $T$ denotes the twist functor around a spherical object of a triangulated $A_{\infty}$-category \cite[Section~5h]{SeidelBk}, and $\tau$ denotes the symplectic Dehn twist \cite[Section~16c]{SeidelBk}. In other words, $\tau_{V_0^{\#}}(V_1^{\#})$ and $\tau^{-1}_{V_1^{\#}}(V_0^{\#})$ fit into exact triangles:
\begin{equation*}
    V_0^{\#} \to V_1^{\#} \to \tau_{V_0^{\#}} V_1^{\#} \to V_0^{\#}[1] \qquad V^{\#}_0 \to V^{\#}_1 \to \tau^{-1}_{V^{\#}_1} V^{\#}_0\to V^{\#}_0[1]
\end{equation*}
where $[1]$ denotes a shift in grading by one and the morphisms are the Floer complexes. Fix an ordered collection of vanishing paths $\gamma_1, \dots, \gamma_m$ and corresponding vanishing cycles $V_1, \dots, V_m$. The braid group $Br_m$ acts freely on the set of all distinguished collections, and the action of the standard ${(i-1)}^{th}$ generator of $Br_m$ gives rise to the Hurwitz move:
\begin{equation}\label{hurwitzmove1}
    (\gamma_1, \dots, \gamma_{i-2}, \gamma_{i-1}, \gamma_{i}, \gamma_{i+1},\dots \gamma_m)\mapsto (\gamma_1, \dots, \gamma_{i-2}, \tau_{\gamma_{i-1}}(\gamma_i), \gamma_{i-1}, \gamma_{i+1},\dots \gamma_m)
\end{equation}
where $\tau_{\gamma_{i-1}}(\gamma_i)$ is the vanishing path obtained by precomposing $\gamma_i$ with a clockwise loop around $\gamma_{i-1}$ (Figure \ref{fig:Hurwitz}). The ordered collection of vanishing paths on the right of (\ref{hurwitzmove1}) is a new distinguished collection. Similarly, the inverse of the standard ${(i-1)}^{th}$ generator of $Br_m$ gives rise to the Hurwitz move
\begin{equation}\label{hurwitzmove2}
    (\gamma_1, \dots, \gamma_{i-2}, \gamma_{i-1}, \gamma_{i}, \gamma_{i+1},\dots \gamma_m)\mapsto(\gamma_1, \dots, \gamma_{i-2},\gamma_{i}, \tau^{-1}_{\gamma_{i}}(\gamma_{i-1}), \gamma_{i+1},\dots \gamma_m)
\end{equation}
where $\tau^{-1}_{\gamma_i}(\gamma_{i-1})$ is obtained by precomposing $\gamma_{i-1}$ with a counter-clockwise loop around $\gamma_i$ (Figure \ref{fig:Hurwitz}); the right-hand side of (\ref{hurwitzmove2}) is a new distinguished collection. Such actions lift to Hurwitz-type moves on the vanishing cycles, which relate any two bases of such objects. Furthermore, if $V'_{i}$ and $V''_{i-1}$ are the vanishing cycles associated to $\tau_{\gamma_{i-1}}(\gamma_i)$ and $\tau^{-1}_{\gamma_i}(\gamma_{i-1})$ respectively, then $V'_{i}=\tau_{V_{i-1}}(V_i)$ and $V''_{i-1}=\tau^{-1}_{V_i}(V_{i-1})$ (\cite[Section~16c]{SeidelBk}).

\begin{remark}
    A Hurwitz move on two consecutive and disjoint vanishing cycles leaves the geometric objects underlying the Lagrangian branes intact, but it inverts the order of the two.
\end{remark}

\subsubsection{Auslander algebras}
The Auslander algebras of Dynkin type $\A$ we have mentioned in Section \ref{section:mainresults} belong to a wider class of algebras. In the more general context of Auslander-Reiten theory, \emph{Auslander algebras} $\Gamma$ are finite-dimensional Artin algebras characterised by (a) $\text{gl.}\dim \Gamma \leq 2$ and (b) if $0 \to \Gamma \to I_0 \to I_1 \to I_2$ is a minimal injective resolution of $\Gamma$, then $I_0$ and $I_1$ are projective $\Gamma$-modules.

Auslander's correspondence (\cite[Section~VI.5]{ARS}) establishes a bijection between the Morita equivalence classes of finite-dimensional $\kk$-algebras ($\kk$ any field) of finite representation type and the Morita equivalence classes of Auslander algebras. Under this correspondence, the Auslander algebra $\Gamma_{\Lambda}$ associated to a $\kk$-algebra $\Lambda$ of finite representation type is the endomorphism algebra of any additive generator of $\textbf{mod}\Lambda$; since every indecomposable $\Lambda$-module $M$ is isomorphic to a summand of any additive generator of $\textbf{mod}\Lambda$, we have that $\Gamma_{\Lambda}=\text{End}_{\Lambda}(\oplus_{[M]}M)$, where $M$ are indecomposable $\Lambda$-modules and the sum runs over all isomorphism classes of them. For us, $\Lambda=A_{n,1}$ will be the path algebra of the linearly oriented $A_n$-quiver, whose corresponding Auslander algebra is $\Gamma_{n+1}:=\Gamma_{A_{n,1}}$.

\subsection{Further directions}\label{section:future}
Auslander algebras can be generalised to \emph{higher Auslander algebras} as algebras $\Gamma$ characterised by (a) $\text{gl.}\dim \Gamma \leq d$ and (b) if $0 \to \Gamma \to I_0 \to I_1 \to \cdots \to I_d$ is a minimal injective resolution of $\Gamma$, then $I_0, \dots I_{d-1}$ are projective $\Gamma$-modules. These are due to Iyama's work for which we refer to \cite{Iyama1} and \cite{Iyama3}. We refer to \cite{DJL} for the definition of $A_{n,d}$, the $d$-dimensional Auslander algebra of Dynkin type $\A$, which is proven to satisfy above inequalities in \cite{Iyama2}. We know from \cite{DJL} that we have an equivalence of categories $\W(\Sym^d(\D),\Lambda_n^{(d)})\simeq \perf(A_{n,d})$, where the collection of stops is defined as 
    \begin{equation*}
        \Lambda^{(d)}_n := \bigcup_{p\in \Lambda_n} \{p\} \times \Sym^{d-1}(\D).
    \end{equation*}
Furthermore, the two-fold symmetric product $\Sym^2(\C)$ can be generalised to the $d$-fold product $\Sym^d(\C)$, which is naturally isomorphic to $\C^d$. We expect results in analogy to Theorem \ref{thm1} and \ref{thm2} in this direction.

\subsection{Structure of the paper}
In Section \ref{section:2}, we provide the necessary background for the computation of the Fukaya-Seidel category of $f_n$, and we use tools developed by A'Campo and Gusein-Zade to prove Theorem \ref{thm1}. Section \ref{section:3} is dedicated to the proof of Theorem \ref{thm2} and Corollary \ref{cor:maineq}; in particular, in Section \ref{section:algorithm} we provide an explicit equivalence of categories $\perf({\tg_n}) \simeq \perf({\Gamma_n})$.

\subsection{Acknowledgements}
The author would like to thank her PhD supervisor Yank{\i} Lekili for suggesting the project, for carefully reading its previous versions, and for the invaluable feedback. She would also like to thank Matthew Habermann for the helpful discussions in the early stages of this project, and Ailsa Keating for her interest in the project and constructive conversations. The author is indebted to the anonymous referee for the generous and extensive comments. The author is supported by the Faculty of Natural, Mathematical \& Engineering Sciences [NMESFS], King's College London.

\section{The Fukaya-Seidel category $\F(f_n)$}\label{section:2}
From now on, we will refer to $f_n$ as the family of maps defined on the symmetric product:
\begin{equation*}
    f_n:\Sym^2(\C) \to \C, \quad (x,y) \mapsto x^n+y^n
\end{equation*}
and to $\fC_n$ as the corresponding polynomial singularities defined on $\C^2$. Under the natural identification:
\begin{equation}\label{SymIsom}
    \varphi: \Sym^2(\C) \xrightarrow{\cong} \C^2 \quad (x,y) \mapsto (x+y, xy),
\end{equation}
 $f_n$ and $\fC_n$ satisfy $\fC_n=\fSym_n\circ \varphi^{-1}$. We define the Fukaya-Seidel category $\F(f_n):=\F(g_n)$, as defined in \cite[Chapter~3]{SeidelBk}.

\subsection{A Morsification of $f_n$}\label{section:morsification}
We make some preliminary considerations on the local behaviour of a holomorphic function  $\fC: \C^2 \to \C$ with an isolated singularity at the origin, and we refer to \cite{Ebeling}, \cite{Dimca} and \cite{Keating} for the results used in this section.

We define an \emph{unfolding} of such $\fC$ to be a holomorphic function germ $F:\C^2 \times \C \to \C$ satisfying $F(z,0)=\fC(z)$, $z\in \C^2$. Any two such unfoldings $F,G$ of $\fC$ are \emph{equivalent} if there exists a transition function between them, i.e.\ a holomorphic map germ $\psi: \C^2 \times \C \to \C^2$, satisfying $\psi(z,0)=z$ and $G(z,u)=F(\psi(z,u),u)$, for any $z\in\C^2$, $u \in \C$. A \emph{Morsification} of $\fC$ is a representative $\C^2 \times \C \to \C$, $(z,\epsilon) \mapsto \fC_{\epsilon}(z)$ of such an unfolding, such that $\fC_{\epsilon}:\C^2 \to \C$ is a Morse function for almost all $\epsilon$ in a neighbourhood of zero. We often call $\fC_{\epsilon}$ a Morsification of $\fC$.  The  \emph{Milnor fibre} of $\fC$ is the fibre above a (sufficiently small) regular value of $\fC_{\epsilon}$.

In general, $\fC_{\epsilon}$ is a Lefschetz fibration when restricted to a suitable open subset of $\C^2$, and one needs to check the behaviour at infinity to extend it to a Lefschetz fibration on the whole $\C^2$; in our case, the holomorphic functions under consideration are Lefschetz fibrations defined on $\C^2$ because they are \emph{tame} (in the sense of \cite{Broughton}, property that directly follows from \cite[Proposition~3.1]{Broughton}). When this holds, the fibre $\fC_{\epsilon}^{-1}(*)$ above a regular value is symplectomorphic (with respect to the standard symplectic structure) to the Milnor fibre of $\fC$ in the usual sense, as shown in \cite[Lemma~2.18]{Keating} (\cite[Lemma~3.3]{Dimca} for the smooth version of the statement).

As a map defined on the symmetric product, $f_n(x,y):=x^n+y^n$ has an isolated singularity at the origin. Under (\ref{SymIsom}), $\Sym^2(\C)$ inherits the standard symplectic structure of $\C^2$. Holomorphicity of the perturbations we consider ensures that Morsifications of $\fC_n:=\fSym_n\circ \varphi^{-1}$ are symplectic Lefschetz fibrations. Under these conditions, we can define parallel transport and Lefschetz thimbles as in Section \ref{section:background}. 

\begin{remark}
    The space of deformations of $\fC_n$ which are Morse is path connected (\cite[Section~2.3]{Keating}), hence any two such unfoldings are equivalent and define the same Morsification, so the $A_{\infty}$-category $\F(f_n)$ is independent of a choice of such unfolding. In this section we will consider two (equivalent) Morsifications of $\fC_n$.
\end{remark}

\subsubsection{Description of critical values}
We first consider a linear Morsification of $\fC_n$, defined as the representative 
\begin{equation*}
    \C^2 \times \C \to \C, \quad (u,v,\epsilon)\mapsto \fC_n(u,v)-\epsilon u=:\fC_{n,\epsilon}(u,v),
\end{equation*}
which we use to study the topology of the Milnor fibre of $\fSym_n$. Explicit computations show that $\fC_{n,\epsilon}$ is Morse for any $\epsilon \neq 0$. By fixing, once and for all, such a generic $\epsilon$, we define $\fC_{n,M} := \fC_{n,\epsilon}$ and $\fSym_{n,M}:=\fC_{n,M}\circ \varphi$, where $\varphi$ is the isomorphism (\ref{SymIsom}).

\begin{lemma}\label{lemma:Milnornumber}
    The Milnor number of $\fSym_n$ is ${n-1} \choose 2$.
\end{lemma}

\begin{proof}
We give a full description of the critical points of $\fSym_{n,M}$ in terms of critical points of $\fC_{n,M}$. To do so, take the lifts to $\C^2$ of $\fSym_{n,M}$, $\fC_{n,M}$ and $\varphi$ and denote these respectively by $\f_n, \g_n$ and $\phii$. If $\pi: \C^2 \to \Sym^2(\C)$ is the branched covering map, then:
\begin{equation*}
    \f_n=f_{n,M}\circ \pi, \quad \phii=\varphi \circ \pi, \quad \g_n=g_{n,M}
\end{equation*}
so that $\f_n=\g_n\circ \phii$ holds. By the chain rule and by surjectivity of $\phii$, it follows that the critical points of $\g_n$ are contained in the image under $\phii$ of the critical points of $\f_n$. On the other hand, the Jacobian of $\phii$
\begin{equation*}
    J_{\phii}=\begin{pmatrix}
        1 & y \\
        1 & x 
        \end{pmatrix}
\end{equation*}
is invertible whenever $x\neq y$. This, together with the chain rule, implies that the images under $\phii$ of the critical points of $\f_n$ away from the diagonal are critical points of $\g_n$. As a map defined on $\C^2$, $\f_n(x,y)=x^n+y^n-\epsilon(x+y)$ has $(n-1)^2$ critical points, given by $\{ (\xi_i,\xi_j)\}:=\{ (\xi^i,\xi^j)\}$ with $\xi$ the scaled $(n-1)^{th}$ root of unity; of these, exactly $n-1$ are of the form $(x,x)$. In particular,
\begin{alignat*}{3}
    \{ (\xi_i+\xi_j,\xi_{i}\xi_j) | i\neq j\} &\subset &\text{Crit}(\g_n) &\subset \{ (\xi_i+\xi_j,\xi_i\xi_j)\} \\
   {n-1}\choose 2 &\leq &|\text{Crit}(\g_n)| &\leq {{n}\choose 2}.
\end{alignat*}

Finally, we claim that all pairs $\{(2\xi_i,\xi_i^2)\}$ are not critical points of $\g_n$. It then follows that $\fC_{n,M}$ (and $f_{n,M}$) has exactly ${n-1 \choose 2}$ critical points, which, for this choice of Morsification, are $\{(\xi^i+\xi^j,\xi^{i+j}) \mid i\neq j\}$.

To prove the claim we observe that, using the binomial expansion of $(x+y)^n$, $\fC_n=\fSym_n \circ \varphi^{-1}$ can be written recursively as follows:
\begin{equation*}
    g_n(u,v)=
    \begin{dcases}
        1 \quad &\text{ if $n=0$}\\
        u \quad &\text{ if $n=1$}\\
        u^n -\sum_{k=1}^{\frac{n}{2}} {n \choose k} g_{n-2k}(u,v)v^k  \quad &\text{ if $n>1$ even}\\[2pt]
        u^n - \sum_{k=1}^{\frac{n-1}{2}} {n \choose k} g_{n-2k}(u,v)v^k \quad &\text{ if $n>1$ odd} \\
    \end{dcases}
\end{equation*}

Differentiating the polynomial $\g_n$ with respect to the first coordinate, evaluating it at $(\xi_i+\xi_j,\xi_i\xi_j)$ and fixing a value $\xi_i$ yields a non-zero polynomial in the one variable $\xi_j$ of degree $n-2$:
\begin{equation}\label{eq:evaluatedpoly}
    n(\xi_i+\xi_j)^{n-1} - \sum {n \choose k} \frac{\partial g_{n-2k}}{\partial u}|_{(\xi_i+\xi_j,\xi_i\xi_j)}(\xi_i\xi_j)^k - \epsilon
\end{equation}
with $\xi_i^{n-1}=\xi_j^{n-1}=-\epsilon/n$. Hence for each fixed $\xi_i$, (\ref{eq:evaluatedpoly}) has at most $n-2$ solutions. Of these, we know $\xi_j$ to be a solution, for all $j\neq i$. Hence $\xi_j=\xi_i$ cannot be a solution for (\ref{eq:evaluatedpoly}), and $(2\xi_i,\xi_i^2)$ cannot be a critical point of $\g_n$, for any index $i$. The claim follows.
\end{proof}

\begin{lemma}\label{factorisation}
        $\fSym_n$ factorises over $\Sym^2(\C)$ in $\floor{ \frac{n+1}{2}}$ terms.
    \end{lemma} 
    
\begin{proof}
    The lift  $\f_n$ of $\fSym_n$ to $\C^2$ admits a linear factorisation:
    \begin{equation*}
        \f_n(x,y)=x^n+y^n=\prod_{k=1}^n(x-\zeta_k y)
    \end{equation*}
    where $\zeta_k$ denotes the $k^{th}$ root of $-1$. Excluding the root $\zeta=-1$ for $n$ odd, these come in pairs of distinct roots $\zeta$ and $\bar{\zeta}$, with $\zeta \bar{\zeta} =1$. 
    The above is not a factorisation of $\fSym_n$ over $\Sym^2(\C)$, but the following factor is:
    \begin{equation*}
        (x-\zeta y)(x-\bar{\zeta} y)=x^2-2\mathrm{Re}(\zeta)xy+y^2 \in \Sym^2(\C)[x,y].
    \end{equation*}
    It follows that
    \begin{equation*}
    \fSym_n(x,y)=
    \begin{dcases}
        \prod_{k=1}^{\frac{n}{2}}\left(x^2-2\mathrm{Re}(\zeta_k)xy+y^2\right) &\quad \text{for $n$ even}\\
        (x+y)\prod_{k=1}^{\frac{n-1}{2}}\left(x^2-2\mathrm{Re}(\zeta_k)xy+y^2\right) &\quad \text{for $n$ odd}\\
    \end{dcases}
    \end{equation*}
    is a factorisation of $f_n$ as a symmetric polynomial. This gives a factorisation of $g_n$:
\begin{equation}\label{eq:factorisation}
    \fC_n(u,v)=
    \begin{dcases}
        \prod_{k=1}^{\frac{n}{2}}\left(u^2-2(1+\mathrm{Re}(\zeta_k))v\right) &\quad \text{for $n$ even}\\
        u\prod_{k=1}^{\frac{n-1}{2}}\left(u^2-2(1+\mathrm{Re}(\zeta_k))v\right) &\quad \text{for $n$ odd}\\
    \end{dcases}
\end{equation}
\end{proof}  

\subsubsection{The Milnor fibre}\label{section:regularfibre}
Fix $\Sigma_n$ to be the regular fibre above zero of $\fC_{n,M}$. In order to give a topological description of $\Sigma_n$, we introduce a second fibration $\rho$ on $\Sigma_n$ given by projection to second coordinate. Since
\begin{equation*}
    \chi(\Sigma_n)+\#\{\text{critical points}\} =\chi(\C^2),
\end{equation*}
$\chi(\Sigma_n)=\frac{n(3-n)}{2}$. We compute the number of punctures by looking at $\rho(u,v)$ approaching infinity. Assuming $n$ even, and considering a circle $\{|v| = N \gg 0\}$, $\rho$ defines a branched cover of degree $n$ of the circle (the degree of the defining polynomial of $\Sigma_n$). The equation defining $\rho^{-1}(N)$ is:
\begin{equation}\label{eq:branchindex}
    0=\fC_{n,M}(u,N) \approx u^n-{n\choose {\frac{n}{2}}}N^{\frac{n}{2}}
\end{equation}
which has at most $\frac{n}{2}=\gcd\left(n,\frac{n}{2}\right)$ solutions for $u$, for fixed $N$. In fact, it has exactly as many: one can check from (\ref{eq:branchindex}) and from the factorisation (\ref{eq:factorisation}) that the preimage of the circle is a collection of $\frac{n}{2}$ circles, each of ramification index 2. Similarly, for $n$ odd, the defining equation for $\rho^{-1}(N)$ is:
\begin{equation*}
    0=\fC_{n,M}(u,N) \approx u^n-{n\choose {\frac{n-1}{2}}}uN^{\frac{n-1}{2}}
\end{equation*}
which has $\frac{n-1}{2}=\gcd\left(n-1,\frac{n-1}{2}\right)$ solutions for $u\neq0$, with an additional solution given by $u=0$. We conclude that the number of punctures is $p=\floor{\frac{n+1}{2}}$. From the Euler characteristic of $\Sigma_n$, if follows that the genus of the regular fibre is

\begin{equation*}
    \text{genus}(\Sigma_n)=
         \begin{dcases}
             \frac{(n-2)^2}{4} \quad &\text{for $n$ even}\\
             \frac{(n-1)(n-3)}{4} \quad &\text{for $n$ odd}.\\
          \end{dcases}
 \end{equation*}

 \subsubsection{The Fukaya category of the Milnor fibre}\label{section:Fcatregfibre}

Given a Lefschetz fibration $f$, this comes with a restriction functor
\begin{equation*}
    \F(f) \to \F(F_*)
\end{equation*}
to the Fukaya category of its regular fibre, given by restricting Lefschetz thimbles to their boundaries on $F_*$ (see Section \ref{section:background}). This allows us to reduce Floer cohomology computations carried out in the total space to ones carried out in the fibre which, for our purposes, is a Riemann surface. Unlike $\F(\fSym_n)$, for which we have generation results, the Fukaya category of the regular fibre is not always generated by a distinguished collection of vanishing cycles; see \cite[Theorem~6.2]{Keating} for a counter-example. As it turns out, in the case at hand, the (compact) Fukaya category of the regular fibre of $f_n$ is indeed generated, as an $A_{\infty}$-category, by the images under the restriction functor of a distinguished collection of thimbles. Let $r={{n-1}\choose 2}$ be the Milnor number of $\fSym_n$, and $V_1<\dots < V_{r}$ the ordered collection of vanishing cycles associated to a Morsification of $f_n$. The following holds.

\begin{proposition}
    A distinguished collection of vanishing cycles on the regular fibre $\Sigma_n$ of $f_n$ generates $\F(\Sigma_n)$ for all $n > 3$. Moreover, there is a quasi-isomorphism
    \begin{equation*}
        [(T_{V_1}\dots T_{V_{r}})^n] \cong [2(n-3)]
     \end{equation*}
     of functors $\F(\Sigma_n) \to \F(\Sigma_n)$, where $T$ denotes a spherical twist as in \cite[Section~5h]{SeidelBk}, and the composition of such twists is the symplectic monodromy.
\end{proposition}

\begin{proof}
    The first part of the statement is a direct consequence of \cite[Theorem~3.3]{Keating}, which follows from Seidel's \cite[Section~4c]{SeidelArt3} and \cite[Proposition~18.17]{SeidelBk}; a combination of these results can be stated as follows. Given a weighted homogeneous polynomial $\pp$ in two (complex) variables of weights $(w_1, w_{1})$ and total weight $w$, with an isolated singularity at the origin and such that the sum of the weights is not equal to $w$, the Fukaya category of the Milnor fibre of $\pp$ is generated by a distinguished collection of vanishing cycles. Indeed, the polynomial expression of $\fSym_n$ as a function on the symmetric product satisfies a quasi-homogeneity condition:
    \begin{equation*}
        \fSym_n(tx,ty)=t^n \fSym_n(x,y).
    \end{equation*}
    It then follows from the isomorphism (\ref{SymIsom}) that $\fC_n=\fSym_n \circ \varphi^{-1}$ is weighted homogeneous, of weights $(1,2)$ and of total weight $n$. This is also clear from the factorisation (\ref{eq:factorisation}) of $\fC_n$. The remaining part of the statement follows directly from \cite[Theorem~4.17]{SeidelArt3} and \cite[Lemma~4.15]{SeidelArt3}.
    \end{proof}

\begin{remark}
    This generation result does not hold in the case of $n=3$. As we will see (Section \ref{ACconf}), the Milnor fibre of $f_3$ is a cylinder $T^*S^1$ with the (single) vanishing cycle that is the zero section equipped with a (fixed) $U(1)$-local system. This does not generate (or even split-generate) the Fukaya category, whose objects include the zero section equipped with \underline{any} $U(1)$-local system.
\end{remark}

\subsection{A Morsification following A'Campo}\label{section:ACmorsification}

The linear Morsification we chose in the Section \ref{section:morsification} allowed us to compute the number of vanishing cycles, as well as the topology of the regular fibre. We will now follow methods developed by A'Campo \cite{AC99} to describe the Milnor fibre of our singularity, together with a favourite collection of vanishing cycles. Such results rely on the notions of \emph{r-divides} and \emph{real deformations} developed independently by A'Campo \cite{AC75} and by Gusein-Zade \cite{GZ}.

\subsubsection{$r$-divides}\label{section:rdivides}
We refer to \cite{AC75}, \cite{GZ} and \cite{AC99} for the complete list of definitions and proofs mentioned in this section, and to \cite[Section~2.4]{Keating} for the symplectic version of A'Campo's work. 

    \begin{definition}[\cite{AC75}, Section~1]
        Let $J$ be the disjoint union of $r$ copies of the interval $[0,1]$ and $\D_{\epsilon}\subset \CR^2$ the closed disk of radius $\epsilon$. An \emph{$r$-divide} of $\D_{\epsilon}$ is an immersion $\alpha: J   \to \D_{\epsilon}$ such that:
        \begin{itemize}
            \item $\alpha(\partial J) \subset \partial \D_{\epsilon}$, $\alpha(\mathring{J}) \subset \mathring{\D}_{\epsilon}$ and $\alpha(J)$ is connected;
            \item $\alpha$ is generic, in the sense that it only has ordinary double points, none of which lie on $\partial \D_{\epsilon}$;
            \item The closures of two distinct \emph{regions} (defined as the connected components of $\D_{\epsilon}\setminus \alpha(J)$ disjoint from $\partial \D_{\epsilon}$) are either disjoint or such that their intersection is either a point or the image of a segment $\alpha(I)$ of $J$.
            \end{itemize}

            A \emph{signed $r$-divide} is an $r$-divide equipped with a sign for each region, such that two regions sharing an edge have different signs.
    \end{definition}

A'Campo considers $r$-divides associated to real deformations of (polynomial) isolated singularities in two variables $\pp(x,y)$. The existence of such a divide follows from results in \cite[Theorem~1]{AC75} and \cite[Section~5]{GZ}, where both authors prove the existence of a real polynomial deformation $\pp_t(x,y):=\pp(x,y;t)$, $t\in \CR$ of $\pp(x,y)$ satisfying the following two conditions for all sufficiently small $t\neq 0$:
\begin{itemize}
    \item The zero locus of $\pp(x,y;t)$ is an $r$-divide;
    \item The number of regions and of double points of the $r$-divide add up to the Milnor number of $\pp(x,y)$ at the origin.
\end{itemize}
As per \cite{Keating}, we call a real deformation satisfying the above conditions a \emph{good real deformation} or a \emph{good real Morsification} of $\pp$. Given such a deformation $\pp_t$, \cite[Theorem~1]{AC99} constructs the real curve $C_t=\{\pp_t(x,y)=0 \mid x,y \in \CR \} \cap \D_{\varepsilon}$, which gives a signed $r$-divide. This associates to every critical point of a good real Morsification of $\pp$ either a region or an intersection point of the $r$-divide. More specifically, it associates to each double point of the divide a critical point of $\pp_t$ whose critical value is zero, and to each positive (resp. negative) region a critical point whose value is positive (resp. negative). We call the former ``saddles'' and the latter ``positive'' (resp. ``negative'') critical points.

\begin{remark}\label{rk:furtherperturbation}
    It should be clear that $\pp_t$, as a polynomial in complex variables, is also a Morsification of $\pp$ as a complex function. We fix, once and for all, a sufficiently small $t \neq 0$, and denote such Morsification as $\pp^{AC}$. $\pp^{AC}$ is not required to have distinct critical values, and in fact this is not the case in our computations. When necessary, one can further perform a small perturbation of $\pp^{AC}$ to separate the critical values \cite[Remark~2.32]{Keating}.
\end{remark}

A'Campo further associates to such $r$-divides an oriented (topological) Riemann surface $\Sigma$, the Milnor fibre of $\pp$, and a collection of vanishing cycles associated to a collection of vanishing paths (see \cite[Example~1]{AC99}, which also describes how to obtain the A'Campo-Gusein-Zade diagram of the singularity). $\Sigma$ is constructed by taking, for each double point of the given $r$-divide, a cylinder embedded into $\CR^3$ with four half twist; additionally, for each segment of the divide, we glue a ribbon-like strip with one half twist to the cylinders that correspond to the boundary points of each segment. The vanishing cycles are given as follows. For each twisted cylinder, draw the curve circling its waist. For each negative (resp. positive) critical point, draw a curve going along the ribbon-like strips and circling the corresponding negative (resp. positive) region. These vanishing cycles correspond to the vanishing paths that are straight lines from the critical values (after possibly considering a further small perturbation, as per Remark \ref{rk:furtherperturbation}) to the regular value $ -i \eta$ for some small $\eta \in \CR_{>0}$  (\cite[Section~2.4.2]{Keating} and again \cite[Example~1]{AC99}). The $\Sigma$ thus constructed is then $(\pp^{AC})^{-1}(-i \eta)$, and is understood as the smoothing of the critical fibre above zero. 
    
There are three families of vanishing cycles arising from A'Campo's construction:
\begin{itemize}
    \item ``Saddle'' vanishing cycles, as waist curves of each twisted cylinder;
    \item ``Negative'' vanishing cycles, as the curves associated to a negative region;
    \item ``Positive'' vanishing cycles, as the curves associated to a positive region.
\end{itemize}
We further denote by saddle (resp. negative, positive) vanishing paths and saddle (resp. negative, positive) thimbles the paths and thimbles associated to saddle (resp. negative, positive) vanishing cycles.

The vanishing paths described above give a total order of the vanishing cycles (and corresponding thimbles) as:
\begin{equation*}
    \{\text{negatives}\} < \{\text{saddles}\} < \{\text{positives}\}.
\end{equation*}
The order of the negative, saddle and positive Lagrangians of the same type does not matter, as all of the vanishing cycles of same type can be Hamiltonian isotoped so that they are pairwise disjoint (\cite{Keating}, Prop. 2.35).

\begin{remark}
    The ``twisting'' we refer to above is meant to provide an embedding of $\Sigma$ into $\CR^3$; changing the direction of the twists changes the embedding, but not $\Sigma$ itself. See \cite[Figure~7]{AC99}, \cite[Figure~6]{Keating} and Figure \ref{ACsurface5} in the next Section for a pictorial description of $\Sigma$.
\end{remark}

This description provides preferred orientations of the vanishing cycles, given by the counter-clockwise orientation on the plane onto which we project the surface (\cite[Section~2.4.3]{Keating}). With the Milnor fibre naturally carrying an orientation, we can define an intersection number between vanishing cycles. This is provided by A'Campo as follows:
\begin{itemize}
    \item $V_i \cdot V_j=-1$ if $V_i$ is a negative vanishing cycle and $V_j$ is a saddle cycle whose corresponding double point in the divide is in the boundary of the region corresponding to $V_i$;
    \item $V_i \cdot V_j=-1$ if $V_j$ is a positive cycle and $V_i$ is a saddle cycle whose corresponding double point is in the boundary of the region corresponding to $V_j$;
    \item $V_i \cdot V_j=-1$ if $V_i$ is a negative cycle,  $V_j$ is a positive one, and the two corresponding regions share an edge.
    \item $V_i \cdot V_j=+1$ otherwise.
\end{itemize}

\subsubsection{A'Campo regular fibre, vanishing cycles and quiver}\label{ACconf}
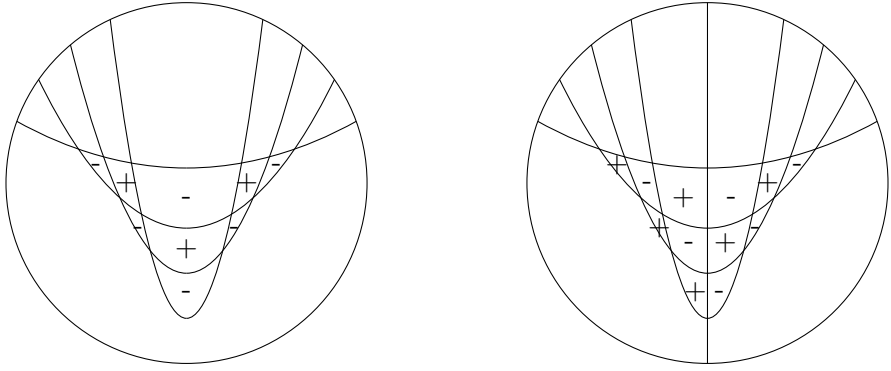
\begin{figure}
    \centering
    \begin{minipage}{.45 \textwidth}
    \centering
\begin{tikzpicture}[scale=.8]
    \draw (0,0) circle (3cm);

    \draw (0,.25) parabola (20:3cm);
    \draw (0,.25) parabola (160:3cm);

    \draw (0,-.75) parabola (35:3cm);
    \draw (0,-.75) parabola (145:3cm);

    \draw (0,-1.5) parabola (50:3cm);
    \draw (0,-1.5) parabola (130:3cm);

    \draw (0,-2.25) parabola (65:3cm);
    \draw (0,-2.25) parabola (115:3cm);

    \node at (0,-1.8){-};

    \node at (0,-1.1){+};
    \node at (.8,-.75){-};
    \node at (-.8,-.75){-};

    \node at (0,-.25){-};
    \node at (1,0){+};
    \node at (-1,0){+};
    \node at (1.5,.3){-};
    \node at (-1.5,.3){-};
    
    \end{tikzpicture}
    \end{minipage}
    \centering
    \begin{minipage}{.45 \textwidth}
        \centering
        \begin{tikzpicture}[scale=.8]
            \draw (0,0) circle (3cm);
    
    \draw (0,.25) parabola (20:3cm);
    \draw (0,.25) parabola (160:3cm);
    
    \draw (0,-.75) parabola (35:3cm);
    \draw (0,-.75) parabola (145:3cm);
    
    \draw (0,-1.5) parabola (50:3cm);
    \draw (0,-1.5) parabola (130:3cm);
    
    \draw (0,-2.25) parabola (65:3cm);
    \draw (0,-2.25) parabola (115:3cm);
    
    \draw (90:3cm) to (-90:3cm);
    
    
    \node at (.2,-1.8){-};
    \node at (-.2,-1.8){+};
    
    \node at (.3,-1){+};
    \node at (.8,-.75){-};
    \node at (-.3,-1){-};
    \node at (-.8,-.75){+};
    
    \node at (.4,-.25){-};
    \node at (1,0){+};
    \node at (-.4,-.25){+};
    \node at (1.5,.3){-};
    \node at (-1,0){-};
    \node at (-1.5,.3){+};
    
\end{tikzpicture}
\end{minipage}
\caption{(left) $\frac{n}{2}$-divide for $n$ even ($n=8$), (right) $\frac{n+1}{2}$-divide for $n$ odd ($n=9$).}
\label{rdivide}
\end{figure}

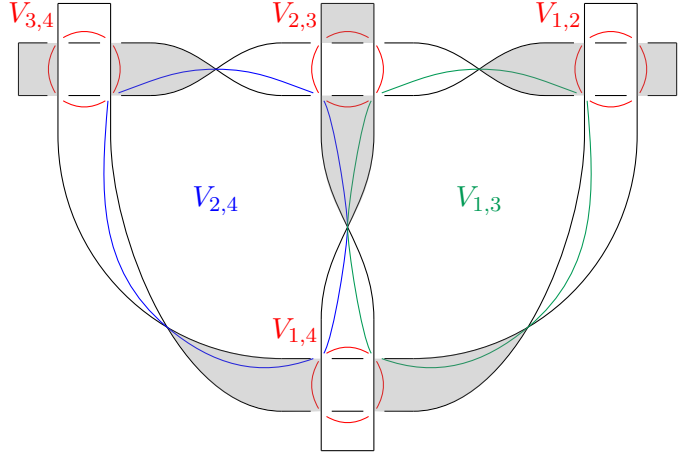
\begin{figure}
    \centering
    \begin{minipage}{1\textwidth}
    \centering
    \begin{tikzpicture}[scale=.7,
        cylinder/.pic={
        \draw (-.35,.875) to (-.35,-.875); 
        \draw (.35,.875) to (.35,-.875); 
        \draw (-.875,.35) to (-.49,.35); 
        \draw (-.21,.35) to (.21,.35); 
        \draw (.875,.35) to (.49,.35); 
        \draw (-.875,-.35) to (-.49,-.35); 
        \draw (-.21,-.35) to (.21,-.35); 
        \draw (.875,-.35) to (.49,-.35);
        \draw[color=red](-.385,.315) to[bend right=30] (-.385,-.315);
        \draw[color=red](-.28,-.42) to[bend right=30] (.28,-.42);
        \draw[color=red](.385,.315) to[bend left=30] (.385,-.315);
        \draw[color=red](-.28,.42) to[bend left=30] (.28,.42);
            }]
            
    \draw (-5,4) pic {cylinder} (0,4) pic {cylinder} (5,4) pic {cylinder} (0,-2) pic {cylinder};

    \draw (-6.25,4.5) to (-6.25,3.5);
    \draw (-5.5,5.25) to (-4.5,5.25);
    \draw (-.5,5.25) to (.5,5.25);
    \draw (5.5,5.25) to (4.5,5.25);
    \draw (6.25,4.5) to (6.25,3.5);
    \draw (-.5,-3.25) to (.5,-3.25);

    \draw[name path=line11] (-3.75,4.5) to[in=180,out=0] (-1.25,3.5);
    \draw[name path=line12] (-3.75,3.5) to[in=180,out=0] (-1.25,4.5);
    \path [name intersections={of=line11 and line12,by=int1}];
            
    \draw[name path=line21] (3.75,4.5) to[in=0,out=180] (1.25,3.5);
    \draw[name path=line22] (3.75,3.5) to[in=0,out=180] (1.25,4.5);
    \path [name intersections={of=line21 and line22,by=int2}];

    \draw[name path=line41] (-.5,2.75) to[in=90,out=270] (.5,-.75);
    \draw[name path=line42] (.5,2.75) to[in=90,out=270] (-.5,-.75);
    \path [name intersections={of=line41 and line42,by=int4}];

    \draw[name path=line31] (-5.5,2.75) to[in=180,out=270,looseness=1] (-1.25,-1.5);
    \draw[name path=line32] (-4.5,2.75) to[in=180,out=270,looseness=.8] (-1.25,-2.5);
    \path [name intersections={of=line31 and line32,by=int3}];

    \draw[name path=line51] (5.5,2.75) to[in=0,out=270,looseness=1] (1.25,-1.5);
    \draw[name path=line52] (4.5,2.75) to[in=0,out=270,looseness=.8] (1.25,-2.5);
    \path [name intersections={of=line51 and line52,by=int5}];
            
    \draw[color=blue] (-4.55,3.4) to[out=265,in=140] (int3) to[out=-40,in=200] (-.65,-1.5);
    \draw[color=blue] (-.45,-1.4) to[out=60,in=270,looseness=.4] (int4) to[out=90,in=-60,looseness=.4] (-.45,3.4);
    \draw[color=blue] (-.65,3.55) to[out=160,in=0,looseness=1] (int1) to[out=180,in=20,looseness=1] (-4.35,3.55);
            
    \draw[color=ForestGreen] (4.55,3.4) to[out=-85,in=40] (int5) to[out=220,in=-20] (.65,-1.5);
    \draw[color=ForestGreen] (.45,-1.4) to[out=120,in=-90,looseness=.4] (int4) to[out=90,in=240,looseness=.4] (.45,3.4);
    \draw[color=ForestGreen] (.65,3.55) to[out=20,in=180,looseness=1] (int2) to[out=0,in=160,looseness=1] (4.35,3.55);

    \draw [fill=gray,draw=none, fill opacity=.3] (-6.25,4.5) -- (-6.25,3.5) -- (-5.5,3.5) -- (-5.5,4.5);

    \draw[fill=gray,draw=none, fill opacity=.3] (-4.5,4.5)  to[in=145,out=0,looseness=1.2] (int1) to[in=0,out=-145,looseness=1.2] (-4.5,3.5);

    \draw[fill=gray,draw=none, fill opacity=.3] (4.5,4.5)  to[in=35,out=180,looseness=1.2] (int2) to[in=180,out=-35,looseness=1.2] (4.5,3.5);

    \draw [fill=gray,draw=none, fill opacity=.3] (6.25,4.5) -- (6.25,3.5) -- (5.5,3.5) -- (5.5,4.5);
    
    \draw[fill=gray,draw=none, fill opacity=.3] (-.5,3.5)  to[in=120,out=270,looseness=1.1] (int4) to[in=270,out=60,looseness=1.1] (.5,3.5);

    \draw [fill=gray,draw=none, fill opacity=.3] (-.5,4.5) -- (-.5,5.25) -- (.5,5.25) -- (.5,4.5);

    \draw[fill=gray,draw=none, fill opacity=.3] (-.5,-1.5)  to[in=-30,out=180,looseness=1] (int3) to[in=183,out=-56,looseness=1.2] (-.5,-2.5);

    \draw[fill=gray,draw=none, fill opacity=.3] (.5,-1.5)  to[in=210,out=0,looseness=1] (int5) to[in=-3,out=236,looseness=1.2] (.5,-2.5);

    \node[color=red] at (-6,5) {$V_{3,4}$};
    \node[color=red] at (-1,5) {$V_{2,3}$};
    \node[color=red] at (4,5) {$V_{1,2}$};
    \node[color=red] at (-1,-1) {$V_{1,4}$};

    \node[color=blue] at (-2.5,1.5) {$V_{2,4}$};
    \node[color=ForestGreen] at (2.5,1.5) {$V_{1,3}$};
\end{tikzpicture}
\end{minipage}
\caption{The Milnor fibre of $g_5$. Shaded is the area of the surface that has same orientation of the ``paper'' onto which we have projected it. The labelling of the vanishing cycles is detailed in Proposition \ref{Dynkinquiverprop}.}
\label{ACsurface5}
\end{figure}

Following the previous section, we consider small, real deformations of the singularities $\fC_n$ whose zero loci, intersected with a small disk in $\CR^2$, give rise to $r$-divides. In particular, given the factorisation of $\fC_n$ described in Lemma \ref{factorisation}, we can choose generic values $h_k \in \CR$ such that, near the origin and intersecting with a small disk in $\CR^2$, the real zero locus in the $uv$-plane of the real deformation
\begin{equation*}
    \fC_n^{AC}(u,v)=
    \begin{dcases}
        \prod_{k=1}^{\frac{n}{2}}(u^2-2(1+\mathrm{Re}(\zeta_k))v+h_k) &\quad \text{for $n$ even}\\
        u\prod_{k=1}^{\frac{n-1}{2}}(u^2-2(1+\mathrm{Re}(\zeta_k))v+h_k) &\quad \text{for $n$ odd}\\
    \end{dcases}
\end{equation*}
is as in Figure \ref{rdivide}. The regular surface and vanishing cycles are then given as in Figure \ref{ACsurface5}, with vanishing paths as straight lines from the regular value $-i\eta$ to the negative, saddle and positive critical values (which have been further perturbed so that they are all distinct). 
    
\begin{proposition}
    $\fC_n^{AC}$ is a good real Morsification of $\fC_n$.
\end{proposition}
\begin{proof}
    It suffices to count the number of vanishing cycles and check that these add up to ${n-1}\choose 2$ (Lemma \ref{lemma:Milnornumber}). One can, for example, do this inductively (for either $n$ even or odd, with an induction step of size 2) on the number of copies of the interval $[0,1]$ embedded into the corresponding divide, which is equal to the number of factors in (\ref{eq:factorisation}). If $n$ is even, one can check that the divide has exactly $\frac{n(n-2)}{4}$ double points, $\frac{n(n-2)}{8}$ negative regions and $\frac{(n-4)(n-2)}{8}$ positive ones; if $n$ is odd, the divide has $\frac{(n-1)^2}{4}$ double points, $\frac{(n-1)(n-3)}{8}$ negative regions and as many positive ones. In both cases, these add up to the Milnor number of $\fC_n$.
\end{proof}

Equipped with arbitrary brane structures, the collection of thimbles associated to the vanishing cycles (and paths) constructed above constitutes a collection of objects generating the appropriate Fukaya-Seidel category, with morphism spaces given by Floer complexes. Denote by $\End_n$ the endomorphism algebra of such collection.

\begin{figure}
        \centering
        \begin{minipage}{.47 \textwidth}
        \centering
        \begin{tikzpicture}
            [scale=.83,align=center, v/.style={draw,shape=circle, fill=black, minimum size=1.1mm, inner sep=0pt, outer sep=0pt},
            every path/.style={shorten >= 1mm, shorten <= 1mm},
            font=\small, label distance=1pt,
            ]
            \node[v,label=180:$\LAC_{1,2}$] (00) at (0,0) {};
            \node[draw, shape=circle, scale=0.9,label=180:$\LAC_{1,3}$] (01-) at (0,1) {-};
            \node[v,label=0:$\LAC_{2,3}$] (11) at (1,1) {};
            \node[v,label=180:$\LAC_{1,4}$] (02) at (0,2) {};
            \node[draw, shape=circle, scale=0.7] (12+) at (1,2) {+};
            \node[v] (22) at (2,2) {};
            \node[draw, shape=circle, scale=0.9] (03-) at (0,3) {-};
            \node[v] (13) at (1,3) {};
            \node[draw, shape=circle, scale=0.9] (23-) at (2,3) {-};
            \node[v] (33) at (3,3) {};

            \node[draw, shape=circle, scale=0.9] (0h1-) at (0,5) {-};
            \node[v] (1h1) at (1,5) {};
            \node[draw, shape=circle, scale=0.9] (h0h1-) at (4,5) {-};
            \node[v,label=-20:$\LAC_{n-3,n-2}$] (h1h1) at (5,5) {};
            \node[v,label=90:$\LAC_{1,n-1}$] (0h2) at (0,6) {};
            \node[draw, shape=circle, scale=0.7] (1h2+) at (1,6) {+};
            \node[v] (h0h2) at (4,6) {};
            \node[draw, shape=circle, scale=0.7] (h1h2+) at (5,6) {+};
            \node[v,label=90:$\LAC_{n-2,n-1}$] (h2h2) at (6,6) {};
                        
            \node at(0,4){$\vdots$};
            \node at(2.5,6){$\dots$};
            \node at(2,4){$\vdots$};
            \node at(4,4){$\iddots$};

            \path
            (01-) edge[->] (00)
            (01-) edge[->] (11)
            (01-) edge[->] (02)

            (02) edge[->] (12+)
            (11) edge[->] (12+)
            (22) edge[->] (12+)
            (13) edge[->] (12+)

            (03-) edge[->] (02)
            (03-) edge[->] (13)

            (23-) edge[->] (13)
            (23-) edge[->] (22)
            (23-) edge[->] (33)

            (0h1-) edge[->] (0h2)
            (0h1-) edge[->] (1h1)

            (0h2) edge[->] (1h2+)
            (1h1) edge[->] (1h2+)

            (h0h1-) edge[->] (h0h2)
            (h0h1-) edge[->] (h1h1)

            (h0h2) edge[->] (h1h2+)
            (h1h1) edge[->] (h1h2+)
            (h2h2) edge[->] (h1h2+)
            ;             
        \end{tikzpicture}
        \end{minipage}
        \qquad
        \begin{minipage}{.47 \textwidth}
        \centering
        \begin{tikzpicture}
            [scale=.83,align=center, v/.style={draw,shape=circle, fill=black, minimum size=1.1mm, inner sep=0pt, outer sep=0pt},
            every path/.style={shorten >= 1mm, shorten <= 1mm},
            font=\small, label distance=1pt,
            every loop/.style={distance=1cm, label=right:}
            ]
            \node[v,label=180:$\LAC_{1,2}$] (00) at (0,0) {};
            \node[draw, shape=circle, scale=0.9,label=180:$\LAC_{1,3}$] (01-) at (0,1) {-};
            \node[v,label=0:$\LAC_{2,3}$] (11) at (1,1) {};
            \node[v,label=180:$\LAC_{1,4}$] (02) at (0,2) {};
            \node[draw, shape=circle, scale=0.7] (12+) at (1,2) {+};
            \node[v] (22) at (2,2) {};
            \node[draw, shape=circle, scale=0.9] (03-) at (0,3) {-};
            \node[v] (13) at (1,3) {};
            \node[draw, shape=circle, scale=0.9] (23-) at (2,3) {-};
            \node[v] (33) at (3,3) {};

            \node[v] (0h1) at (0,5) {};
            \node[draw, shape=circle, scale=0.7] (1h1+) at (1,5) {+};
            \node[draw, shape=circle, scale=0.7] (h0h1+) at (4,5) {+};
            \node[v,label=-20:$\LAC_{n-3,n-2}$] (h1h1) at (5,5) {};
            \node[draw, shape=circle, scale=0.9,label=90:$\LAC_{1,n-1}$] (0h2-) at (0,6) {-};
            \node[v] (1h2) at (1,6) {};
            \node[v] (h0h2) at (4,6) {};
            \node[draw, shape=circle, scale=0.9] (h1h2-) at (5,6) {-};
            \node[v,label=90:$\LAC_{n-2,n-1}$] (h2h2) at (6,6) {};

            \node at(0,4){$\vdots$};
            \node at(2.5,6){$\dots$};
            \node at(2,4){$\vdots$};
            \node at(4,4){$\iddots$};

            \path
            (01-) edge[->] (00)
            (01-) edge[->] (11)
            (01-) edge[->] (02)

            (02) edge[->] (12+)
            (11) edge[->] (12+)
            (22) edge[->] (12+)
            (13) edge[->] (12+)
    
            (03-) edge[->] (02)
            (03-) edge[->] (13)

            (23-) edge[->] (13)
            (23-) edge[->] (22)
            (23-) edge[->] (33)
            
            (0h1) edge[->] (1h1+)
            (1h2) edge[->] (1h1+)

            (h0h2) edge[->] (h0h1+)
            (h1h1) edge[->] (h0h1+)

            (0h2-) edge[->] (0h1)
            (0h2-) edge[->] (1h2)

            (h1h2-) edge[->] (h0h2)
            (h1h2-) edge[->] (h1h1)
            (h1h2-) edge[->] (h2h2)
            ;             
        \end{tikzpicture}
        \end{minipage}
    \caption{Quiver of the $f_n$ singularity for (left) $n$ odd and (right) $n$ even. The vertices denoted by $+$ and $-$ symbols represent, respectively, sinks and sources of the quiver. For generic $n$, the quiver associated to the $f_n$ singularity has $n-2$ rows.}
    \label{ACquiver}
    \end{figure}
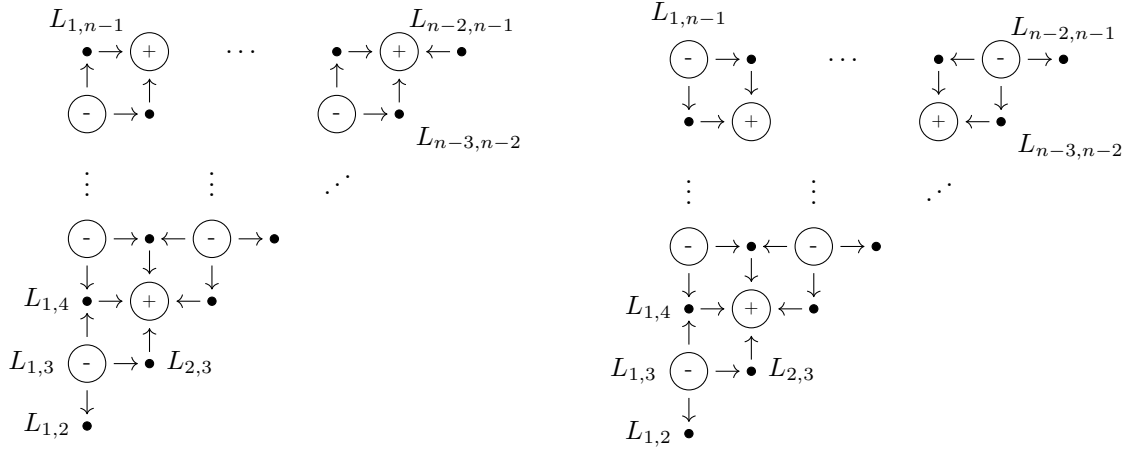

\begin{proposition}\label{Dynkinquiverprop}
    $\End_n$, as an ungraded algebra, is isomorphic to $\tg_n$, the path algebra of the quiver in Figure \ref{ACquiver} modulo commutativity of the squares.
\end{proposition}

\begin{remark}\label{AGdiagram}
    Our quiver in Figure \ref{ACquiver} differs from the A'Campo-Gusein-Zade diagram of the singularity by the following:
    \begin{itemize}
        \item It is directed, with the direction of the arrows indicating a negative intersection number;
        \item Arrows from sources to sinks are suppressed, instead viewed as paths of length two giving rise to a generator in the path algebra.
    \end{itemize}
    We will refer to this as the \emph{quiver associated to the singularity $f_n$}.
\end{remark}

\begin{proof}
The undirected graph underlying the quiver  in Figure \ref{ACquiver} (after taking into consideration Remark \ref{AGdiagram}) is the one prescribed by A'Campo in \cite[Example~1]{AC99}. From this, we denote the bottom-left vertex as the double-indexed $\LAC_{1,2}$, and we increase the first (resp. second) index as we move right (resp. up). These vertices correspond to the vanishing cycles given by A'Campo's construction, which we now denote by the double-indexed $V_{I,J}$:
    \begin{equation*}
        \{V_{I,J} \mid 1\leq I< J\leq n-1\}
    \end{equation*}
    and to their corresponding thimbles, now denoted by the double-indexed $\DAC_{I,J}$:
    \begin{equation*}
        \{\DAC_{I,J} \mid 1\leq I< J\leq n-1\}.
    \end{equation*}
    We claim that the bijection:
    \begin{equation*}
        L_{I,J} \leftrightarrow \DAC_{I,J}
    \end{equation*}
    gives rise to an isomorphism $\End_n \cong \tg_n$ of ungraded algebras. By construction, we have the following:
    \begin{itemize}
        \item For $I,J \equiv 1 \mod 2$, $\DAC_{I,J}$ are negative thimbles, corresponding to a source $\LAC_{I,J}$ of the quiver;
        \item For $I,J \equiv 0 \mod 2$, $\DAC_{I,J}$ are positive thimbles corresponding to a sink  $\LAC_{I,J}$;
        \item In all other cases $\DAC_{I,J}$ (and corresponding vertices $\LAC_{I,J}$) are saddles.
    \end{itemize}

    The path algebra $\tg_n$ of the given quiver with relations is generated by the following paths:
    \begin{itemize}
        \item For each vertex $\LAC_{I,J}$, a path of length 0;
        \item For each source $\LAC_{I,J}$, a path of length 1 from $\LAC_{I,J}$ to each of the saddles $\LAC_{I\pm 1,J}, \LAC_{I,J\pm 1}$ (where these exist);
        \item For each sink $\LAC_{I,J}$, a path of length 1 from each of the saddles $\LAC_{I\pm 1,J}, \LAC_{I,J\pm 1}$ (where these exist) to $\LAC_{I,J}$;
        \item For each source $\LAC_{I,J}$, a path of length 2 from $\LAC_{I,J}$ to each of the four sinks $\LAC_{I\pm 1,J\pm 1}$ (where these exist).
    \end{itemize}

    By construction, each arrow in the quiver corresponds to an intersection point between two vanishing cycles, and therefore a generator of the Floer complex between the corresponding thimbles. Moreover, each path of length 2 corresponds exactly to a morphism from a negative thimble $\DAC_{I,J}$ to each of the four positive thimbles (when these exist) $\DAC_{I\pm 1,J\pm 1}$. Fix a negative vanishing cycle $V_{I,J}$ and assume $V_{I+1,J+1}$ (resp. $V_{I+1,J-1}$, $V_{I-1,J+1}$, $V_{I-1,J-1}$) exists. The intersection between the vanishing cycles arises as in Figure \ref{fig:trianglesAC}, where we can observe two immersed triangles bounded by, respectively, the counter-clock ordered unions $V_{\cup}:= V_{0} \cup V_{1} \cup V_{2}$ and  $V'_{\cup}:= V_{0} \cup V'_{1} \cup V_{2}$, where $V_0:=V_{I,J}$, $V_{1}:=V_{I,J+1}$ (resp. $V_{I,J-1}$, $V_{I,J+1}$, $V_{I,J-1}$), $V'_{1}:=V_{I+1,J}$ (resp. $V_{I+1,J}$, $V_{I-1,J}$, $V_{I-1,J}$) and $V_{2}:=V_{I+1,J+1}$ (resp. $V_{I+1,J-1}$, $V_{I-1,J+1}$, $V_{I-1,J-1}$). It is known (\cite[Section~13b]{SeidelBk}) that a signed count of such immersed triangles (up to the boundary, as discussed in Seidel's book) contributes to the products in the appropriate Fukaya category:
    \begin{align}
        CF^*(V_1,V_2) \otimes CF^*(V_0,V_1) \to& CF^*(V_0,V_2) \nonumber \\
        \label{eq:products} CF^*(V'_1,V_2) \otimes CF^*(V_0,V'_1) \to& CF^*(V_0,V_2).
    \end{align}

    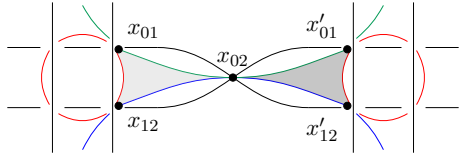
\begin{figure}
    \centering
        \begin{minipage}{1\textwidth}
        \centering
        \scalebox{.8}{
         \begin{tikzpicture}[v/.style={draw,shape=circle, fill=black, minimum size=1.2mm, inner sep=0pt, outer sep=0pt},
                        cylinder/.pic={
        \draw (-.5,1.25) to (-.5,-1.25); 
        \draw (.5,1.25) to (.5,-1.25); 
        \draw (-1.25,.5) to (-.7,.5); 
        \draw (-.3,.5) to (.3,.5); 
        \draw (1.25,.5) to (.7,.5); 
        \draw (-1.25,-.5) to (-.7,-.5); 
        \draw (-.3,-.5) to (.3,-.5); 
        \draw (1.25,-.5) to (.7,-.5);
        \draw[color=red](-.55,.45) to[bend right=30] (-.55,-.45);
        \draw[color=red](-.4,-.6) to[bend right=30] (.4,-.6);
        \draw[color=red](.55,.45) to[bend left=30] (.55,-.45);
        \draw[color=red](-.4,.6) to[bend left=30] (.4,.6);
        }]
                
        \draw (-5,4) pic {cylinder} (0,4) pic {cylinder};
                
        \draw[name path=line11] (-3.75,4.5) to[in=180,out=0] (-1.25,3.5);
        \draw[name path=line12] (-3.75,3.5) to[in=180,out=0] (-1.25,4.5);
        \path [name intersections={of=line11 and line12,by=int1}];
                
        \draw[color=blue] (-.65,3.55) to[out=160,in=0,looseness=1] (int1) to[out=180,in=20,looseness=1] (-4.35,3.55); 
        \draw[color=blue] (-.45,3.4) to[bend left=10] (0,2.8);
        \draw[color=blue] (-4.55,3.4) to[bend right=10] (-5,2.8);  
        \draw[color=ForestGreen] (-.65,4.45) to[out=-160,in=0,looseness=1] (int1) to[out=180,in=-20,looseness=1] (-4.35,4.45);  
        \draw[color=ForestGreen] (-.45,4.65) to[bend right=10] (0,5.2);
        \draw[color=ForestGreen] (-4.55,4.65) to[bend left=10] (-5,5.2); 
        \draw [fill=lightgray,draw=none, fill opacity=.3]  (-4.35,3.55) to[out=20,in=180,looseness=1] (int1) to[out=180,in=-20,looseness=1] (-4.35,4.45) -- (-4.45,4.45) to[bend left=30] (-4.45,3.55) ;
    
        \draw [fill=darkgray,draw=none, fill opacity=.3]  (-.65,3.55)to[out=-200,in=0,looseness=1] (int1) to[out=0,in=-160,looseness=1] (-.65,4.45) -- (-.55,4.45) to[bend right=30] (-.55,3.55) ;

        \node[v,label=90:$x_{02}$] at (int1) {};
        \node[v,label=80:$x_{01}$] at (-4.4,4.47) {};
        \node[v,label=280:$x_{12}$] at (-4.4,3.53) {};
        \node[v,label=100:$x'_{01}$] at (-.6,4.47) {};
        \node[v,label=260:$x'_{12}$] at (-.6,3.53) {};
                
        \end{tikzpicture}
            }
        \end{minipage}
        \caption{In green and blue respectively, a negative and positive vanishing cycle $V_{I,J}$ and $V_{I+1,J+1}$; in red, the two vanishing cycles $L_{I+1,J}$ and $L_{I,J+1}$. Shaded, the triangles contributing to the compositions.}
        \label{fig:trianglesAC}
    \end{figure}

    As we can observe from the explicit description of the fibre, there is no other triangle contributing to either of (\ref{eq:products}), so the products are given by $(x_{12},x_{01}) \mapsto \pm x_{02}$ and $(x'_{12},x'_{01}) \mapsto \pm x_{02}$ respectively, where $x_{01},x'_{01},x_{12}, x'_{12}$ and $x_{02}$ are the generators of the appropriate Floer complexes. The $\pm$ signs here depend on the orientation of the moduli spaces of such holomorphic triangles, and denote two possible choices of generators of the corresponding morphism spaces. We claim that these can all arranged to be positive, so that the products (\ref{eq:products}) reflect the commutativity of the squares
    \begin{equation*}
        \centering
        \begin{tikzcd}
            \LAC_{I,J+1} \arrow{r}  & \LAC_{I+1,J+1} \\
            \LAC_{I,J}  \arrow{u} \arrow{r}& \LAC_{I+1,J} \arrow{u}
        \end{tikzcd}
        \begin{tikzcd}
            \LAC_{I,J} \arrow{r} \arrow{d} & \LAC_{I+1,J}\arrow{d} \\
            \LAC_{I,J-1} \arrow{r}& \LAC_{I+1,J-1} 
        \end{tikzcd}
        \begin{tikzcd}
            \LAC_{I-1,J+1} & \LAC_{I,J+1} \arrow{l} \\
            \LAC_{I-1,J}  \arrow{u}& \LAC_{I,J} \arrow{u} \arrow{l}
        \end{tikzcd}
        \begin{tikzcd}
            \LAC_{I-1,J}\arrow{d} & \LAC_{I,J}\arrow{d} \arrow{l} \\
            \LAC_{I-1,J-1} & \LAC_{I,J-1} \arrow{l}
        \end{tikzcd}
    \end{equation*}
    in the quiver associated to the singularity. Specifically, this shows that the isomorphism between $\End_n$ and $\tg_n$ respects compositions. The proof of this claim is inductive on the number of rows of the quiver in Figure \ref{ACquiver}, and follows from the fact that the quiver itself is planar. Fix arbitrary generators of the Floer complexes. Starting from $V_{1,2}$ and increasing (in order) the second and first indexes, the first compositions we encounter appear in the square:
    \begin{equation*}
        \begin{tikzcd}
            V_{1,4} \arrow{r}  & V_{2,4} \\
            V_{1,3}  \arrow{u} \arrow{r}& V_{2,3} \arrow{u}
            \end{tikzcd}
    \end{equation*}
    If the square commutes keep the signs of the corresponding Floer complexes unchanged. Otherwise, reverse the sign of the composition involving the top arrow. By iteratively increasing the indexes of $V_{I,J}$ (proceeding right and up on the quiver), whenever we encounter a square we either leave the signs unchanged (if it commutes) or reverse the sign of the top horizontal morphism (so that it does).
\end{proof}

From now on, we denote by $\DAC_{I,J}$ the Lefschetz thimbles whose corresponding vanishing cycles are $V_{I,J}$, and we follow the notation of the indexes detailed in the Proof of Proposition \ref{Dynkinquiverprop}.

\begin{proposition}\label{prop:degree0}
    $\End_n=\textnormal{End}(\bigoplus_{I,J}\DAC_{I,J}) $ is an $A_{\infty}$-algebra, concentrated in degree 0 and with vanishing differential and higher $A_{\infty}$-products.
\end{proposition}

\begin{proof}
    The only polygons we can observe on the Milnor fibre are the triangles described in the proof of Proposition \ref{Dynkinquiverprop}, so all $A_{\infty}$-products apart from composition vanish. The thimbles, as object of the Fukaya-Seidel category, are gradable. In fact,  since $2c_1(\C^2)$ and $H^1(\C^2)$ both vanish, $\C^2$ carries a canonical grading (\cite[Section~12]{SeidelBk} and \cite{SeidelArt3}), which induces a $\Z$-grading on the Floer complexes $CF^*(V_{I,J},V_{I',J'})$ and which can be lifted to a grading on the Floer complexes between thimbles. It follows that $\End_n$ admits a unique (up to a global shift) $\Z$-grading. Having fixed an arbitrary grading on the objects, we iteratively perform a series of shifts and prove that all morphisms are concentrated in degree 0.

    Starting from $\DAC_{1,2}$, leave its grading unchanged. There is only one morphism space involving $\DAC_{1,2}$: $CF^*(\DAC_{1,3},\DAC_{1,2})$; we can shift the grading of $\DAC_{1,3}$ so that this lies in degree 0. Similarly, we can (independently) shift the gradings of $\DAC_{1,4}$, $\DAC_{2,3}$ and $\DAC_{2,4}$ so that $CF^*(\DAC_{1,3},\DAC_{1,4})$, $CF^*(\DAC_{1,3},\DAC_{2,3})$ and $CF^*(\DAC_{2,3},\DAC_{2,4})$ all lie in degree 0. As the products (\ref{eq:products}) have the same output, $CF^*(\DAC_{1,4},\DAC_{2,4})$ also has to lie in degree 0. This reflects the commutativity of the square:
    \begin{equation*}
        \begin{tikzcd}
            \DAC_{1,4} \arrow[r]  &\DAC_{2,4} \\
            \DAC_{1,3}  \arrow[u,"0"] \arrow[r,"0"]& \DAC_{2,3} \arrow[u,"0"]
        \end{tikzcd}
    \end{equation*}
   By repeating the process of iteratively increasing the indexes of $\DAC_{I,J}$, we can shift all degrees on the objects so that all morphisms are concentrated in degree 0.
\end{proof}

\begin{theorem}\label{thm1intext}
    For $n\geq 3$, there is a quasi-equivalence of triangulated $A_{\infty}$-categories
    \begin{equation*}
        \F(f_n) \xlongrightarrow{\simeq} \perf(\tg_n).
    \end{equation*}
\end{theorem}

\begin{proof}
    This directly follows from Propositions \ref{Dynkinquiverprop} and \ref{prop:degree0} and generation results of the two triangulated categories by the constructed collection of thimbles and by projective $\tg_n$-modules respectively.
\end{proof}

\section{The derived equivalence $\F(f_n) \simeq \perf(\Gamma_n)$}\label{section:3}

\subsection{The derived equivalence $\F(f_n) \simeq \W^2_n$}\label{section:W2n}

This section is dedicated to proving Theorem \ref{thm2} and Corollary \ref{cor:maineq}.

For $n \in \N_{\geq 3}$, we consider the disk $\D$ (with standard orientation of the boundary) and a set of $n$ points $\Lambda_n$ on $\partial \D$. We equip the symmetric product $\Sym^2(\D)$ with a natural symplectic structure $\omega$ coming from a choice of positive area form $\alpha$ on $\D$, as prescribed by Perutz in \cite[Corollary~7.2]{Perutz}. Away from the diagonal, $\omega$ is the smooth pushforward $\pi_*(\alpha^{\times 2})$, where $\pi: \D^2 \to \Sym^2(\D)$ is the branched covering map which ramifies along the diagonal of $\Sym^2(\D)$.

\begin{remark}
	$\Sym^2(\D)$ can also be equipped with the standard symplectic structure coming from the identification $\Sym^2(\C) \cong \C^2$ (the former being the Liouville completion of $\Sym^2(\D)$). From (\cite[Proposition~1.1]{Perutz}), it follows that the two symplectic structures on $\Sym^2(\D)$ are equivalent, as they both tame the complex structure $\Sym^2(J)$ induced by the standard complex structure $J$ on $\D$ (\cite[Remark~1.1.1]{DJL}). Motivated by \cite{Auroux2} and \cite{DJL}, we prefer to use the symplectic structure obtained from the symmetric product construction, as this allows us to consider Lagrangian submanifolds of $\Sym^2(\D)$ arising from products of pairwise disjoint Lagrangian submanifolds of $\D$.
\end{remark}

Consider the symplectic hypersurfaces $\Lambda^{(2)}_n$ defined as:
\begin{equation*}
    \Lambda^{(2)}_n := \bigcup_{p\in \Lambda_n} \{p\} \times \D;
\end{equation*}
we call these \emph{stops} and, following \cite{Auroux2} (and the more general theory developed in \cite{GPS1} and \cite{GPS2}), we construct the \emph{partially wrapped Fukaya category} $\W^2_n:=\W(\Sym^2(\D),\Lambda^{(2)}_n)$. Products of arcs in $\D \setminus \Lambda_n$ are contractible, hence admit a unique choice of $\Z$-grading, up to a global shift; a choice of grading on such Lagrangian subspaces determines an object of the partially wrapped Fukaya category. From now on, we equip $\W^2_n$ with this canonical $\Z$-grading. We have the following result by Auroux on generation of $\W^2_n$ in terms of products of arcs.

\begin{theorem}[\cite{Auroux1}, Theorem 1]\label{Auroux}
    Let $\Lambda$ be a finite set of points on $\partial \D$ as above, $L_1,...,L_m$ a collection of disjoint properly embedded arcs in $\D$ with endpoints in $\partial \D\setminus \Lambda$. Assume that $\D \setminus (L_1 \cup ...\cup L_m)$ is a disjoint union of disks, each of which contains at most one point of $\Lambda$. Then, the partially wrapped Fukaya category $\W^2_n$ is generated by the $m \choose 2$ Lagrangian submanifolds $L_{ij} := L_i \times L_j$, products of distinct pairs of arcs.
\end{theorem}

\begin{remark}
    The above Theorem holds (and is, in fact, stated by Auroux), for any compact Riemann surface with non-empty boundary $\Sigma$, and for the partially wrapped Fukaya category $\W(\Sym^d(\Sigma),\Lambda^{(d)})$, where $\Lambda^{(d)}$ is the set of symplectic hypersurfaces $\Lambda^{(d)} := \bigcup_{p\in \Lambda} \{p\} \times \Sym^{d-1}(\Sigma)$ defined for any $d\geq 1$ (see Section \ref{section:future}). In this case, a set of generators is given by the $m \choose d$ Lagrangian submanifolds $L_I:=\prod_{i\in I}L_i$, where $I$ ranges over the $d$-element subsets of $\{1,\dots m\}$.
\end{remark}
    
\begin{notation}
    We introduce here some notation on the objects of $\W^2_n$ of the form $L=L_i \times L_j$. Fix $\Lambda_n$ as above and label each boundary component of $\D$ clockwise, from $0$ to $n-1$. We denote $L$ by $i_0i_1 \times j_0j_1$, where $L_i$ has endpoints on the boundary arcs labelled $i_0$ and $i_1$, and $L_j$ has  endpoints on $j_0$ and $j_1$. We fix, once and for all, such labelling such that $i_0< i_1$ and $j_0<j_1$.
\end{notation}

Morphisms between products of Lagrangian arcs satisfying the hypotheses of Theorem \ref{Auroux} are generated by all products of Reeb chords induced by the Reeb flows along $\partial \D$, which are the rotational flows in the counter-clockwise direction of the boundary (\cite[Sections~4.1, 4.2]{Auroux2}). More in general, if the Lagrangian arcs in $(\D,\partial \D)$ are not disjoint, we can observe intersection points in the interior of $\D$, which can contribute to the morphism spaces between corresponding Lagrangians. This will never arise in our computations, as we will only work with generators of $\W^2_n$ satisfying the hypotheses of Theorem \ref{Auroux}.

With respect to suitable grading structures, Auroux provides (\cite[Lemma~5.2]{Auroux2}) the following quasi-isomorphism in $\W^2_n$, for any pairwise distinct $i,j,k$:
\begin{equation}\label{isoAuroux}
    ij \times jk \simeq ij \times ik
\end{equation}
as well as the exact triangles in $\W^2_n$, for any $p,q$, $0\leq i < j < k\leq n-1$:
\begin{equation}\label{trianglesAuroux}
    pq \times ij \xrightarrow{id \otimes x} pq \times ik  \xrightarrow{id \otimes y} pq \times jk  \xrightarrow{id \otimes z} (pq \times ij)[1]
\end{equation}
where $x,y,z$ denote the morphism spaces generated by the appropriate Reeb chords, and $id$ denotes the identity morphism.

\begin{figure}
    \begin{minipage}{.45 \textwidth}
        \centering
        \begin{tikzpicture}[scale=.7,align=center, v/.style={draw,shape=circle, fill=black, minimum size=1.2mm, inner sep=0pt, outer sep=0pt}, font=\small, label distance=1pt,
            ]
        \draw (0,0) circle (3cm);
    
        \node[v] at(165:3cm){};
        \node[v] at(140:3cm){};
        \node[v] at(115:3cm){};
        \node[v] at(90:3cm){};
        \node[v] at(65:3cm){};
        \node[v] at(40:3cm){};
        \node[v] at(15:3cm){};
    
        \node[v] at(207.5:3cm){};
        \node[v] at(232.5:3cm){};
        \node[v] at(257.5:3cm){};
        \node[v] at(282.5:3cm){};
        \node[v] at(307.5:3cm){};
        \node[v] at(332.5:3cm){};
    
        \node at (183:3.2cm){$\vdots$};
        \node at (-3:3.2cm){$\vdots$};
    
        \node at (77.5:3.4cm){0};
        \node at (52.5:3.4cm){1};
        \node at (27.5:3.4cm){2};
    
        \node at (320:3.4cm){$\frac{n-5}{2}$};
        \node at (295:3.4cm){$\frac{n-3}{2}$};
        \node at (270:3.4cm){$\frac{n-1}{2}$};
        \node at (245:3.4cm){$\frac{n+1}{2}$};
        \node at (220:3.5cm){$\frac{n+3}{2}$};
    
        \node at (152.5:3.6cm){$n-3$};
        \node at (130:3.6cm){$n-2$};
        \node at (102.5:3.4cm){$n-1$};
            
        \draw[color=blue] (80:3cm)[inner sep=0pt] to[bend left=45] (100:3cm);
        \draw[color=blue] (78:3cm)[inner sep=0pt] to[bend left=30] (128:3cm);
        \draw[color=blue] (130:3cm)[inner sep=0pt] to[bend right=25] (55:3cm);
        \draw[color=blue] (53:3cm)[inner sep=0pt] to[bend left=20] (150:3cm);
        \draw[color=blue] (152:3cm)[inner sep=0pt] to[bend right=20] (30:3cm);

        \draw[color=blue] (220:3cm)[inner sep=0pt] to[bend left=25] (322:3cm);
        \draw[color=blue] (248:3cm)[inner sep=0pt] to[bend left=30] (320:3cm);
        \draw[color=blue] (250:3cm)[inner sep=0pt] to[bend left=35] (297:3cm);
        \draw[color=blue] (275:3cm)[inner sep=0pt] to[bend left=45] (295:3cm);    
        \end{tikzpicture}
        \end{minipage}
        \qquad
        \begin{minipage}{.45 \textwidth}
            \centering
            \begin{tikzpicture}[scale=.7,align=center, v/.style={draw,shape=circle, fill=black, minimum size=1.2mm, inner sep=0pt, outer sep=0pt}, font=\small, label distance=1pt, every loop/.style={distance=1cm, label=right:}
                ]
                \draw (0,0) circle (3cm);
    
        \node[v] at(165:3cm){};
        \node[v] at(140:3cm){};
        \node[v] at(115:3cm){};
        \node[v] at(90:3cm){};
        \node[v] at(65:3cm){};
        \node[v] at(40:3cm){};
        \node[v] at(15:3cm){};
    
        \node[v] at(220:3cm){};
        \node[v] at(245:3cm){};
        \node[v] at(270:3cm){};
        \node[v] at(295:3cm){};
        \node[v] at(320:3cm){};
    
        \node at (190:3.2cm){$\vdots$};
        \node at (-10:3.2cm){$\vdots$};
    
        \node at (77.5:3.4cm){0};
        \node at (52.5:3.4cm){1};
        \node at (27.5:3.4cm){2};
    
        \node at (307.5:3.6cm){$\frac{n}{2}-2$};
        \node at (282.5:3.4cm){$\frac{n}{2}-1$};
        \node at (257.5:3.4cm){$\frac{n}{2}$};
        \node at (232.5:3.6cm){$\frac{n}{2}+1$};
    
        \node at (152.5:3.6cm){$n-3$};
        \node at (130:3.6cm){$n-2$};
        \node at (102.5:3.4cm){$n-1$};
    
        \draw[color=blue] (80:3cm)[inner sep=0pt] to[bend left=45] (100:3cm);
        \draw[color=blue] (103:3cm)[inner sep=0pt] to[bend right=30] (57:3cm);
        \draw[color=blue] (130:3cm)[inner sep=0pt] to[bend right=25] (55:3cm);
        \draw[color=blue] (33:3cm)[inner sep=0pt] to[bend left=20] (132:3cm);
        \draw[color=blue] (152:3cm)[inner sep=0pt] to[bend right=20] (30:3cm);
    
        \draw[color=blue] (210:3cm)[inner sep=0pt] to[bend left=25] (313:3cm);
        \draw[color=blue] (233:3cm)[inner sep=0pt] to[bend left=30] (310:3cm);
        \draw[color=blue] (235:3cm)[inner sep=0pt] to[bend left=35] (285:3cm);
        \draw[color=blue] (260:3cm)[inner sep=0pt] to[bend left=45] (280:3cm);
    
    \end{tikzpicture}
    \end{minipage}
\caption{A collection of arcs satisfying the hypotheses of Theorem \ref{Auroux} for (left) $n$ odd and (right) $n$ even.}
\label{swingingarcs}
\end{figure}
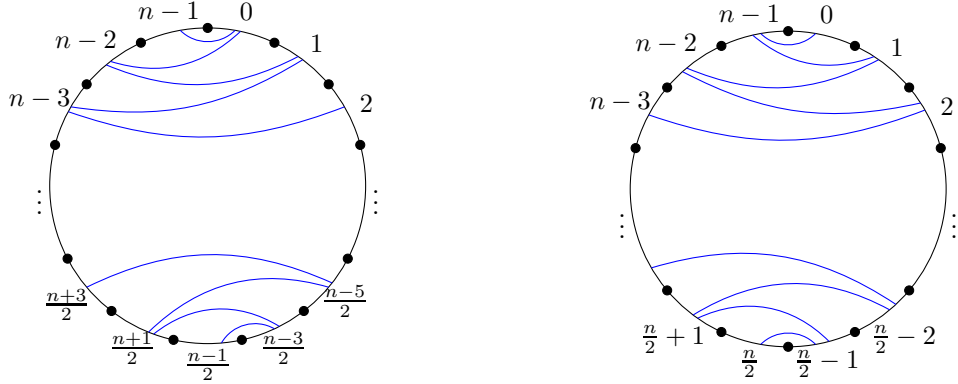

Consider now the disk equipped with a set of $n$ marked points $\Lambda_n$ on its boundary as before. Take the collection of disjoint, properly embedded, arcs as given in Figure \ref{swingingarcs}. Since these satisfy the hypotheses of Theorem \ref{Auroux}, the collection of products of two distinct arcs generates $\W^2_n$. Denote by $\mathcal{B}_n$ the endomorphism algebra of such collection.

\begin{theorem}\label{qefn}
    For $n\geq 3$, there is a quasi-equivalence of categories $\perf(\tg_n)\simeq \W^2_n$, induced by the isomorphism of $\kk$-algebras
    \begin{equation*}
        \tg_n \xlongrightarrow{\cong} \mathcal{B}_n.
    \end{equation*}
\end{theorem}
        
\begin{proof}
    We start by giving a one-to-one correspondence between the vertices $\{L_{I,J}\}$ of the quiver in Figure \ref{ACquiver} (whose path algebra is $\tg_n$) and pairs of Lagrangian arcs given in Figure \ref{swingingarcs}, whose endomorphism algebra is $\mathcal{B}_n$. Denote by $I=|a-b|$ the arc with endpoints labelled $a$ and $b$ on the boundary of the disk; for fixed $n$, there are $n-1$ such arcs and ${n-1}\choose 2$ unordered pairs of them. If $ab \times cd$ is such a pair, and $I=|a-b|$, $J=|c-d|$, $I<J$, we claim that the bijection that associates
    \begin{equation*}
        L_{I,J} \leftrightarrow I\times J
    \end{equation*}
    for $1\leq I<J \leq n-1$ gives rise to an isomorphism of $\kk$-algebras.

    Fix a product of arcs $I\times J$ as above. It is easy to see from the description in Figure \ref{swingingarcs} that the following holds:
    \begin{itemize}
        \item If $I,J\equiv 1\mod 2$, there is a Reeb chord giving rise to a morphism from $I \times J$ to each of the four products $I \times (J\pm 1)$, $(I\pm 1) \times J$, whenever these exist;
        \item  If $I,J\equiv 0\mod 2$, there is a Reeb chord from each of the four products $I \times (J\pm 1)$, $(I\pm 1) \times J$ to $I \times J$, whenever these exist;
        \item If $I,J\equiv 1\mod 2$, there is a Reeb chord giving rise to a morphism from $I \times J$ to each of the four products $(I\pm 1) \times (J\pm 1)$, whenever these exist.
    \end{itemize}
    Such morphisms are in one-to-one correspondence with the paths of length 1 and 2 in the quiver, as described in Proposition \ref{Dynkinquiverprop}. Denote  by $x$ (resp. $y$) the morphism space generated by the Reeb chord going from $I$ to $I+1$ (resp. $J$ to $J+1$), and by $id$ the identity morphism; the (unique) morphism from $I \times J$ to $(I + 1) \times (J + 1)$ then arises as the compositions:
    \begin{equation*}
        \begin{tikzcd}
            I \times J+1 \arrow[r,"x \otimes id"]  & I+1 \times J+1 \\
            I \times J  \arrow[u,"id \otimes y"] \arrow[r, "x \otimes id"]& I+1 \times J \arrow[u,"id \otimes y"] 
        \end{tikzcd}
    \end{equation*}
   The same argument can be applied to the squares: 
    \begin{equation*}
        \begin{tikzcd}
            I \times J \arrow{d}  \arrow{r} & I+1 \times J \arrow{d}   \\
            I   \arrow{r} \times J-1 & I+1 \times J-1 
        \end{tikzcd} 
        \begin{tikzcd}
            I-1 \times J+1  & I \times J+1 \arrow{l}\\
            I-1 \times J \arrow{u} & I \times J \arrow{u} \arrow{l}
        \end{tikzcd}
        \begin{tikzcd}
            I-1 \times J \arrow{d}  & I \times J \arrow{d}  \arrow{l}\\
            I-1 \times J-1 & I \times J-1  \arrow{l}
        \end{tikzcd}
    \end{equation*}
    which reflect the commutativity relations in the quiver. This establishes an isomorphism of (ungraded) algebras between $\tg_n$ and $\mathcal{B}_n$. Using an argument completely analogous to the one made in Proposition \ref{prop:degree0}, one can view $\mathcal{B}_n$ as an $A_{\infty}$-algebra concentrated in degree 0, so that $\tg_n \cong \mathcal{B}_n$ as graded algebras. As the generators of $\mathcal{B}_n$ and $\tg_n$ generate $\W^2_n$ and $\perf(\tg_n)$ respectively, the equivalence of triangulated categories follows.
\end{proof}

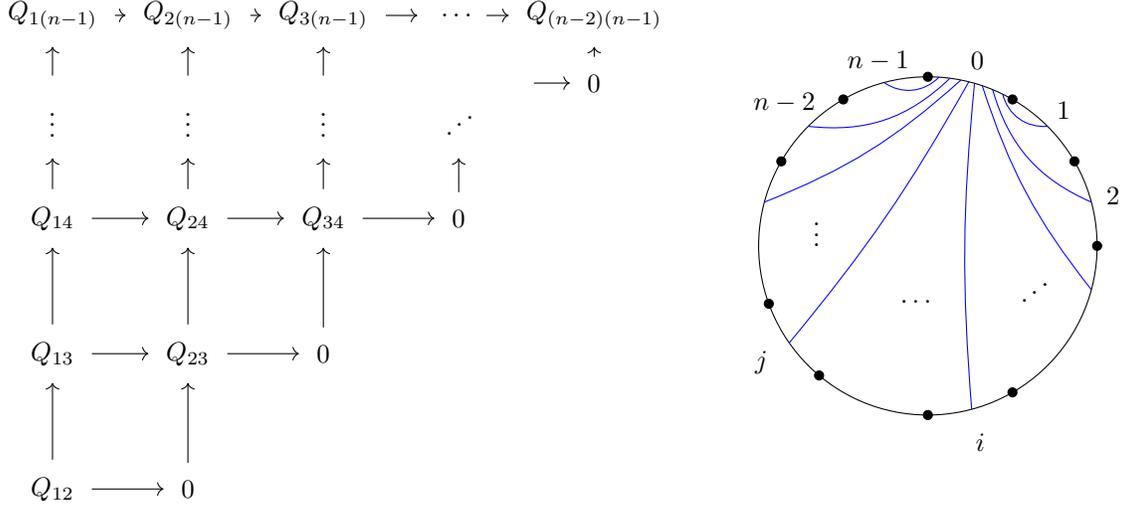
\begin{figure}
    \centering
    \begin{minipage}{.6 \textwidth}
        \centering
    \begin{tikzpicture}
        [scale=1.8,align=center, every path/.style={shorten >= 1mm, shorten <= 1mm}, font=\small, label distance=1pt,
        ]

        \node (00) at (0,0) {$\Q_{12}$};
        \node (01) at (0,1) {$\Q_{13}$};
        \node (11) at (1,1) {$\Q_{23}$};
        \node (02) at (0,2) {$\Q_{14}$};
        \node (12) at (1,2) {$\Q_{24}$};
        \node (22) at (2,2) {$\Q_{34}$};

        \node (aa) at (0,3.5) {$\Q_{1(n-1)}$};
        \node (bb) at (1,3.5) {$\Q_{2(n-1)}$};
        \node (cc) at (2,3.5) {$\Q_{3(n-1)}$};
       
        \node (ee) at (4,3.5) {$\Q_{(n-2)(n-1)}$};

        \node (10) at (1,0) {$0$};
        \node (21) at (2,1) {$0$};
        \node (32) at (3,2) {$0$};

        \node (43) at (4,3) {$0$};

        \path
        (00) edge[color=black,->] (01)
        (01) edge[color=black,->] (11)
        (01) edge[color=black,->] (02)

        (02) edge[color=black,->] (12)
        (11) edge[color=black,->] (12)
        (12) edge[color=black,->] (22)

        (02) edge[color=black,->] (0,2.5)
        (12) edge[color=black,->] (1,2.5)
        (22) edge[color=black,->] (2,2.5)

        (0,3) edge[color=black,->] (aa)
        (1,3) edge[color=black,->] (bb)
        (2,3) edge[color=black,->] (cc)

        (00) edge[color=black,->] (10)
        (10) edge[color=black,->] (11)
        (11) edge[color=black,->] (21)
        (21) edge[color=black,->] (22)
        (22) edge[color=black,->] (32)
        (32) edge[color=black,->] (3,2.5)

        (3.5,3) edge[color=black,->] (43)
        (43) edge[color=black,->] (ee)
        ;    
        
        \node at(0,2.75){$\vdots$};
        \node at(1,2.75){$\vdots$};
        \node at(2,2.75){$\vdots$};

        \node at(3,3.5){$\dots$};
        \node at(3,2.75){$\iddots$};

        \path
        (aa) edge[color=black,->] (bb)
        (bb) edge[color=black,->] (cc)
        (cc) edge[color=black,->] (2.75,3.5)
        (3.15,3.5) edge[color=black,->] (ee)
        ;
    \end{tikzpicture}
\end{minipage}
\qquad
\begin{minipage}{.34 \textwidth}
    \centering
    \begin{tikzpicture}[scale=.75,align=center, v/.style={draw,shape=circle, fill=black, minimum size=1.2mm, inner sep=0pt, outer sep=0pt}, font=\small, label distance=1pt,
        ]
    \draw (0,0) circle (3cm);

    \node[v] at(90:3cm){};
    \node[v] at(60:3cm){};
    \node[v] at(30:3cm){};
    \node[v] at(00:3cm){};

    \node[v] at(120:3cm){};
    \node[v] at(150:3cm){};

    \node[v] at(-90:3cm){};
    \node[v] at(300:3cm){};
    \node[v] at(200:3cm){};
    \node[v] at(230:3cm){};

    \draw[color=blue] (63.75:3cm)[inner sep=0pt] to[bend right=40] (45:3cm);
    \draw[color=blue] (67.5:3cm)[inner sep=0pt] to[bend right=25] (15:3cm);
    \draw[color=blue] (71.25:3cm)[inner sep=0pt] to[bend right=10] (-15:3cm);
    \draw[color=blue] (74:3cm)[inner sep=0pt] to[bend right=5] (285:3cm);
    \draw[color=blue] (76:3cm)[inner sep=0pt] to[bend left=5] (215:3cm);
    \draw[color=blue] (78.75:3cm)[inner sep=0pt] to[bend left=10] (165:3cm);
    \draw[color=blue] (82.5:3cm)[inner sep=0pt] to[bend left=25] (135:3cm);
    \draw[color=blue] (86.25:3cm)[inner sep=0pt] to[bend left=40] (105:3cm);

    \node at (75:3.4cm){0};
    \node at (45:3.4cm){1};
    \node at (15:3.4cm){2};

    \node at (285:3.6cm){$i$};
    \node at (215:3.6cm){$j$};

    \node at (135:3.6cm){$n-2$};
    \node at (105:3.4cm){$n-1$};

    \node at(260:1cm){$\dots$};
    \node at(170:2cm){$\vdots$};
    \node at(-20:2cm){$\iddots$};

    \end{tikzpicture}
\end{minipage}
\caption{(left) Quiver whose path algebra is $\Gamma_n$, with all possible commutativity relations. (right) Iyama generators of $\W^2_n$, whose endomorphism algebra is $\mathcal{A}_{n}$. The isomorphism of graded algebras $\Gamma_n \cong \mathcal{A}_n$ arises from the bijection $Q_{ij}\leftrightarrow 0i\times 0j$.}
\label{Iyamagens2}
\end{figure}

\cite[Theorem~1]{DJL} gives a quasi-equivalence of triangulated $A_{\infty}$-categories $\W^{2}_n \simeq \perf(\Gamma_n)$. Underlying this equivalence, there is a quasi-isomorphism (\cite[Theorem~2.2.3]{DJL}) of differential graded $\kk$-algebras
\begin{equation*}
    \Gamma_n \xlongrightarrow{\cong} \mathcal{A}_{n}
\end{equation*}
between the Auslander algebra $\Gamma_n$ and the endomorphism algebra $\mathcal{A}_{n}$ of a distinguished set of generators of $\W^{2}_n$, given by the collection of Lagrangians
\begin{equation*}
    \{0i \times 0j, 1\leq i < j \leq n-1\}
\end{equation*}
(Figure \ref{Iyamagens2}); we call these \emph{Iyama generators} of $\W^{2}_n$. This collection of pairs of arcs satisfies the hypotheses of Theorem \ref{Auroux} and thus generates the partially wrapped Fukaya category. The above equivalence gives a correspondence between the Iyama generators of $\perf(\Gamma_n)$ given in Section \ref{section:mainresults} (Figure \ref{algebrasgamman}) and the Iyama generators of $\W^{2}_n$.

\begin{remark}\label{Serresymp}
    Following \cite{DJL}, we note that there is a natural symplectomorphism $(\D,\Lambda_n) \to (\D,\Lambda_n)$ cyclically permuting points in $\Lambda_n$ (when these are fixed to be the $n^{th}$ roots of unity), given by rotation by $\frac{2\pi}{n}$. This lifts to a graded symplectomorphism $\Sym^2(\D) \to \Sym^2(\D)$ preserving $\Lambda^{(2)}_n$, which in turn can be used to extract a natural autoequivalence of $\W^{2}_n$. \cite[Proposition~2.5.1]{DJL} states that this autoequivalence agrees with the Serre functor on $\perf(\Gamma_n)$.
\end{remark}

\subsection{The combing algorithm}\label{section:algorithm}
Consequence of Theorem \ref{qefn} and \cite[Theorem~1]{DJL} is an abstract equivalence of $A_{\infty}$-categories:
\begin{equation*}
    \perf(\tg_n) \xrightarrow{\simeq} \perf(\Gamma_n).
\end{equation*}
The aim of this section is to make this equivalence explicit: we provide a constructive algorithm relating the two natural collections of generators of the respective categories to each other. In view of Corollary \ref{cor:maineq}, we will take three parallel approaches in presenting this algorithm, detailed in Sections \ref{algorithmdisk}, \ref{section:quiveralgorithm} and \ref{section:vcyclesalgorithm} respectively. Essentially, we are providing the same algorithmic description that, to a set of generators of $\W^2_n$, $\perf(\Gamma_n)$ and $\F(f_n)$ respectively, associates another collection of generators of $\W^2_n$, $\perf(\Gamma_n)$ and $\F(f_n)$.

\subsubsection{Algorithm on $\W^2_n$}\label{algorithmdisk}
Let us first describe the algorithm that, to the A'Campo generators of $\W^2_n$, associates the Iyama generators. The algorithm consists in a series of \emph{mutations} on the generating collection, and relies on the existence of the triangles (\ref{trianglesAuroux}) provided by Auroux. Suppose $pq \times ij$ and $pq \times ik$ are objects in a given collection $\G$ of generators, with $u$ a generator of the morphism space between them. A \emph{mutation} on $\G$ replaces the pair $(pq \times ij, pq \times ik)$ with the mapping cone of $u$ and either of the two original objects (the choice will be specified every time).

\begin{proposition}[Combing algorithm]\label{combalgorithm}
    Let $\G$ be the collection of generators of $\W^2_n$ given in Theorem \ref{qefn}. Then there exists a series of mutations on $\G$ that replaces it with the Iyama generators of $\W^2_n$.
\end{proposition}

\begin{proof}
    The proof of the claim is constructive. It differs slightly between $n$ even and odd; in both cases, we will give the algorithm inductively on $n$ (with an induction step of size 1).  Note that, for the base case $n=3$, the A'Campo generator $01 \times 02$ is exactly the distinguished Iyama generator.

Suppose first $n+1$ even. After applying the symplectomorphism by rotation as in Remark \ref{Serresymp} (which, on the disk model, is simply a relabelling of the arcs in $\D$), the A'Campo generators of $\W^2_n$ are as in Figure \ref{combalgorithmfigeven}; we call this Step 0.

\begin{itemize}
    \item STEP $A$. Forgetting the arc $(n-1)n$, the pairs of remaining arcs are exactly as they appear in Figure \ref{swingingarcs} (left). Denote by $\G$ the collection of generators, and by $\Gg$ the sub-collection:
    \begin{equation*}
       \Gg:= \left\{pq \times rs \mid pq \neq (n-1)n \text{ and } rs \neq (n-1)n\right\}.
    \end{equation*}
    Assume inductively that there exists a series of mutations on  $\Gg$ that replaces $\Gg$ with the following sub-collection of the Iyama generators:
    \begin{equation*}
        \left\{0i \times 0j \mid i \neq n \text{ and } j \neq n\right\}.
    \end{equation*}
    Additionally, for $n>3$ and increasing $1\leq h \leq n-3$, perform a mutation on $\G \setminus \Gg$ that replaces (in order) the generator $\frac{h+1}{2}\frac{2n-h-3}{2}\times (n-1)n$ (resp. $\frac{h}{2}\frac{2n-h-4}{2}\times (n-1)n$) with the object $0 \frac{h+1}{2} \times (n-1)n$ (resp. and $0 \frac{2n-h-4}{2} \times (n-1)n$) for $h$ odd (resp. even). Call this mutation Step $A_h$. We can perform such mutations thanks to the existence of the following exact triangles (provided by Auroux), for $h$ odd and even respectively:
    \begin{equation*}
        \begin{tikzcd}
            &  0\frac{2n-h-3}{2}\times (n-1)n \arrow[r] & \frac{h+1}{2}\frac{2n-h-3}{2}\times (n-1)n \arrow[d] \\
            &             & 0 \frac{h+1}{2} \times (n-1)n \arrow[ul, "{[1]}"]
        \end{tikzcd}
    \end{equation*}
    \begin{equation*}
       \begin{tikzcd}
           &  \frac{h}{2}\frac{2n-h-4}{2}\times (n-1)n \arrow[r] & 0\frac{h}{2}\times (n-1)n \arrow[d] \\
            &             & 0 \frac{2n-h-4}{2} \times (n-1)n\arrow[ul, "{[1]}"]
        \end{tikzcd}
    \end{equation*}
    where $[1]$ denotes a shift in grading by one, with respect to the brane structure of the Lagrangians discussed in Section \ref{section:W2n}. The cone of the horizontal morphism (whose source and target belong to $\G \setminus \Gg$ after Step $A_{h-1}$) is, for both $h$ odd and even, exactly the new object of the distinguished collection.
    \item STEP $B$. We now perform the simultaneous mutations on $\G \setminus \Gg$ that replace each of the generators $0p\times (n-1)n$ with $0p \times 0n$, for $p \in \{1,\dots, n-1\}$; this is again an admissible mutation due to the existence of the exact triangles:
    \begin{equation*}
        0p \times (n-1)n \to 0p \times 0(n-1) \to 0p \times 0n \xrightarrow{[1]} 0p \times (n-1)n
        \end{equation*}
    for all $1 \leq p \leq n-2$, given by (\ref{trianglesAuroux}),   as well as the quasi-isomorphism given by (\ref{isoAuroux}):
    \begin{equation*}
        0(n-1) \times (n-1)n \simeq 0(n-1)\times 0n.
    \end{equation*}
    \end{itemize}

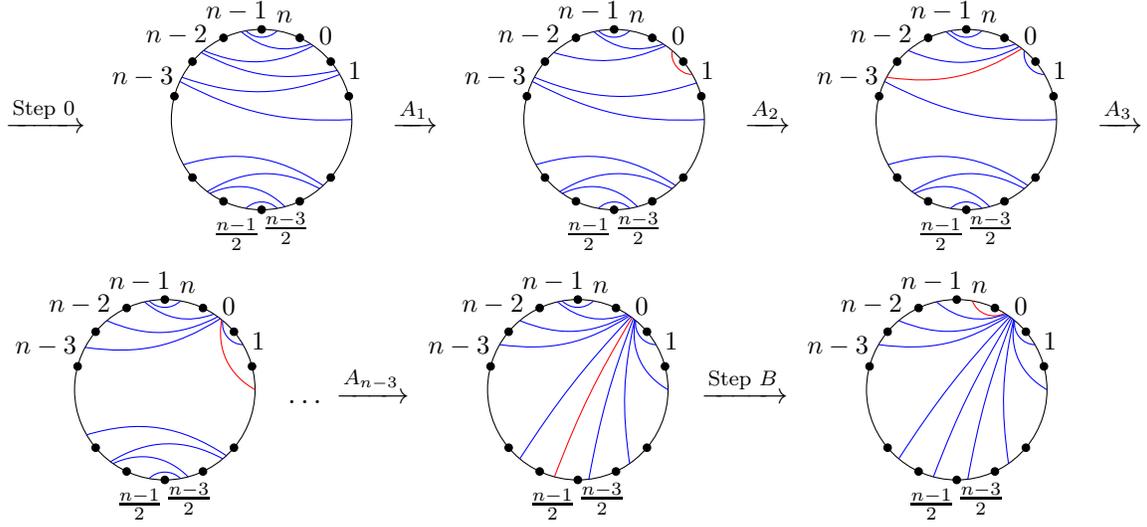
\begin{figure}
      $\xrightarrow{\text{Step 0}}$
      \begin{minipage}{.24 \textwidth}
          \begin{tikzpicture}[scale=.4,align=center, v/.style={draw,shape=circle, fill=black, minimum size=1mm, inner sep=0pt, outer sep=0pt}, font=\small, label distance=1pt,
              ]
            \draw (0,0) circle (3cm);
            \node[v] at(165:3cm){};
            \node[v] at(140:3cm){};
            \node[v] at(115:3cm){};
            \node[v] at(90:3cm){};
            \node[v] at(65:3cm){};
            \node[v] at(40:3cm){};
            \node[v] at(15:3cm){};
        
            \node[v] at(220:3cm){};
            \node[v] at(245:3cm){};
            \node[v] at(270:3cm){};
            \node[v] at(295:3cm){};
            \node[v] at(320:3cm){};

        \node at (77.5:3.5cm){$n$};
        \node at (52.5:3.5cm){0};
        \node at (27.5:3.5cm){1};
    
        \node at (282.5:3.7cm){$\frac{n-3}{2}$};
        \node at (257.5:3.9cm){$\frac{n-1}{2}$};
    
        \node at (160:4.2cm){$n-3$};
        \node at (135:4cm){$n-2$};
        \node at (103:3.7cm){$n-1$};

        \draw[color=blue] (80:3cm)[inner sep=0pt] to[bend left=45] (100:3cm);
        \draw[color=blue] (103:3cm)[inner sep=0pt] to[bend right=30] (57:3cm);
        \draw[color=blue] (130:3cm)[inner sep=0pt] to[bend right=25] (55:3cm);
        \draw[color=blue] (33:3cm)[inner sep=0pt] to[bend left=20] (132:3cm);
        \draw[color=blue] (152:3cm)[inner sep=0pt] to[bend right=20] (30:3cm);
        \draw[color=blue] (155:3cm)[inner sep=0pt] to[bend right=15] (0:3cm);
    
        \draw[color=blue] (210:3cm)[inner sep=0pt] to[bend left=25] (313:3cm);
        \draw[color=blue] (233:3cm)[inner sep=0pt] to[bend left=30] (310:3cm);
        \draw[color=blue] (235:3cm)[inner sep=0pt] to[bend left=35] (285:3cm);
        \draw[color=blue] (260:3cm)[inner sep=0pt] to[bend left=45] (280:3cm);
    \end{tikzpicture}
    \end{minipage}
    $\xrightarrow{A_1}$
    \begin{minipage}{.24 \textwidth}
    \begin{tikzpicture}[scale=.4,align=center, v/.style={draw,shape=circle, fill=black, minimum size=1mm, inner sep=0pt, outer sep=0pt}, font=\small, label distance=1pt,
        ]
        \draw (0,0) circle (3cm);
        \node[v] at(165:3cm){};
        \node[v] at(140:3cm){};
        \node[v] at(115:3cm){};
        \node[v] at(90:3cm){};
        \node[v] at(65:3cm){};
        \node[v] at(40:3cm){};
        \node[v] at(15:3cm){};
            
        \node[v] at(220:3cm){};
        \node[v] at(245:3cm){};
        \node[v] at(270:3cm){};
        \node[v] at(295:3cm){};
        \node[v] at(320:3cm){};
    
        \node at (77.5:3.5cm){$n$};
        \node at (52.5:3.5cm){0};
        \node at (27.5:3.5cm){1};
    
        \node at (282.5:3.7cm){$\frac{n-3}{2}$};
        \node at (257.5:3.9cm){$\frac{n-1}{2}$};
            
        \node at (160:4.2cm){$n-3$};
        \node at (135:4cm){$n-2$};
        \node at (103:3.7cm){$n-1$};

        \draw[color=blue] (80:3cm)[inner sep=0pt] to[bend left=45] (100:3cm);
        \draw[color=blue] (103:3cm)[inner sep=0pt] to[bend right=30] (57:3cm);
        \draw[color=red] (50:3cm)[inner sep=0pt] to[bend right=50] (30:3cm);
        \draw[color=blue] (55:3cm)[inner sep=0pt] to[bend left=20] (132:3cm);
        \draw[color=blue] (152:3cm)[inner sep=0pt] to[bend right=20] (24:3cm);
        \draw[color=blue] (155:3cm)[inner sep=0pt] to[bend right=15] (0:3cm);
    
        \draw[color=blue] (210:3cm)[inner sep=0pt] to[bend left=25] (313:3cm);
        \draw[color=blue] (233:3cm)[inner sep=0pt] to[bend left=30] (310:3cm);
        \draw[color=blue] (235:3cm)[inner sep=0pt] to[bend left=35] (285:3cm);
        \draw[color=blue] (260:3cm)[inner sep=0pt] to[bend left=45] (280:3cm);
                            
        \end{tikzpicture}
        \end{minipage}
    $\xrightarrow{A_2}$
    \begin{minipage}{.24 \textwidth}
        \begin{tikzpicture}[scale=.4,align=center, v/.style={draw,shape=circle, fill=black, minimum size=1mm, inner sep=0pt, outer sep=0pt}, font=\small, label distance=1pt,
        ]
        \draw (0,0) circle (3cm);
        \node[v] at(165:3cm){};
        \node[v] at(140:3cm){};
        \node[v] at(115:3cm){};
        \node[v] at(90:3cm){};
        \node[v] at(65:3cm){};
        \node[v] at(40:3cm){};
        \node[v] at(15:3cm){};

        \node[v] at(220:3cm){};
        \node[v] at(245:3cm){};
        \node[v] at(270:3cm){};
        \node[v] at(295:3cm){};
        \node[v] at(320:3cm){};
    
        \node at (77.5:3.5cm){$n$};
        \node at (52.5:3.5cm){0};
        \node at (27.5:3.5cm){1};
    
        \node at (282.5:3.7cm){$\frac{n-3}{2}$};
        \node at (257.5:3.9cm){$\frac{n-1}{2}$};
    
        \node at (160:4.2cm){$n-3$};
        \node at (135:4cm){$n-2$};
        \node at (103:3.7cm){$n-1$};
            
        \draw[color=blue] (80:3cm)[inner sep=0pt] to[bend left=45] (100:3cm);
        \draw[color=blue] (103:3cm)[inner sep=0pt] to[bend right=30] (57:3cm);
        \draw[color=blue] (50:3cm)[inner sep=0pt] to[bend right=50] (30:3cm);
        \draw[color=blue] (130:3cm)[inner sep=0pt] to[bend right=25] (54:3cm);
        \draw[color=red] (152:3cm)[inner sep=0pt] to[bend right=20] (52:3cm);
        \draw[color=blue] (155:3cm)[inner sep=0pt] to[bend right=15] (0:3cm);
    
        \draw[color=blue] (210:3cm)[inner sep=0pt] to[bend left=25] (313:3cm);
        \draw[color=blue] (233:3cm)[inner sep=0pt] to[bend left=30] (310:3cm);
        \draw[color=blue] (235:3cm)[inner sep=0pt] to[bend left=35] (285:3cm);
        \draw[color=blue] (260:3cm)[inner sep=0pt] to[bend left=45] (280:3cm);
        \end{tikzpicture}
        \end{minipage}
    $\xrightarrow{A_3}$
    \begin{minipage}{.24 \textwidth}
        \begin{tikzpicture}[scale=.4,align=center, v/.style={draw,shape=circle, fill=black, minimum size=1mm, inner sep=0pt, outer sep=0pt}, font=\small, label distance=1pt,
            ]
            \draw (0,0) circle (3cm);
            \node[v] at(165:3cm){};
            \node[v] at(140:3cm){};
            \node[v] at(115:3cm){};
            \node[v] at(90:3cm){};
            \node[v] at(65:3cm){};
            \node[v] at(40:3cm){};
            \node[v] at(15:3cm){};
        
            \node[v] at(220:3cm){};
            \node[v] at(245:3cm){};
            \node[v] at(270:3cm){};
            \node[v] at(295:3cm){};
            \node[v] at(320:3cm){};
    
            \node at (77.5:3.5cm){$n$};
            \node at (52.5:3.5cm){0};
            \node at (27.5:3.5cm){1};
        
            \node at (282.5:3.7cm){$\frac{n-3}{2}$};
            \node at (257.5:3.9cm){$\frac{n-1}{2}$};
        
            \node at (160:4.2cm){$n-3$};
            \node at (135:4cm){$n-2$};
            \node at (103:3.7cm){$n-1$};
    
        \draw[color=blue] (80:3cm)[inner sep=0pt] to[bend left=45] (100:3cm);
        \draw[color=blue] (103:3cm)[inner sep=0pt] to[bend right=30] (57:3cm);
        \draw[color=blue] (50:3cm)[inner sep=0pt] to[bend right=50] (30:3cm);
        \draw[color=blue] (130:3cm)[inner sep=0pt] to[bend right=25] (54:3cm);
        \draw[color=blue] (152:3cm)[inner sep=0pt] to[bend right=20] (52:3cm);
        \draw[color=red] (51:3cm)[inner sep=0pt] to[bend right=35] (0:3cm);

        \draw[color=blue] (210:3cm)[inner sep=0pt] to[bend left=25] (313:3cm);
        \draw[color=blue] (233:3cm)[inner sep=0pt] to[bend left=30] (310:3cm);
        \draw[color=blue] (235:3cm)[inner sep=0pt] to[bend left=35] (285:3cm);
        \draw[color=blue] (260:3cm)[inner sep=0pt] to[bend left=45] (280:3cm);
    
        \end{tikzpicture}
        \end{minipage}
        $\dots$ $\xrightarrow{A_{n-3}}$
        \begin{minipage}{.24 \textwidth}
            \begin{tikzpicture}[scale=.4,align=center, v/.style={draw,shape=circle, fill=black, minimum size=1mm, inner sep=0pt, outer sep=0pt},font=\small, label distance=1pt,
                ]
            \draw (0,0) circle (3cm);

            \node[v] at(165:3cm){};
            \node[v] at(140:3cm){};
            \node[v] at(115:3cm){};
            \node[v] at(90:3cm){};
            \node[v] at(65:3cm){};
            \node[v] at(40:3cm){};
            \node[v] at(15:3cm){};

            \node[v] at(220:3cm){};
            \node[v] at(245:3cm){};
            \node[v] at(270:3cm){};
            \node[v] at(295:3cm){};
            \node[v] at(320:3cm){};

        \node at (77.5:3.5cm){$n$};
        \node at (52.5:3.5cm){0};
        \node at (27.5:3.5cm){1};
    
        \node at (282.5:3.7cm){$\frac{n-3}{2}$};
        \node at (257.5:3.9cm){$\frac{n-1}{2}$};
    
        \node at (160:4.2cm){$n-3$};
        \node at (135:4cm){$n-2$};
        \node at (103:3.7cm){$n-1$};

        \draw[color=blue] (80:3cm)[inner sep=0pt] to[bend left=45] (100:3cm);
        \draw[color=blue] (103:3cm)[inner sep=0pt] to[bend right=30] (57:3cm);
        \draw[color=blue] (130:3cm)[inner sep=0pt] to[bend right=25] (56:3cm);
        \draw[color=blue] (55:3cm)[inner sep=0pt] to[bend left=20] (150:3cm);
        \draw[color=blue] (50:3cm)[inner sep=0pt] to[bend right=50] (30:3cm);
        \draw[color=blue] (50.5:3cm)[inner sep=0pt] to[bend right=35] (0:3cm);

        \draw[color=blue] (230:3cm)[inner sep=0pt] to[bend left=5] (54:3cm);
        \draw[color=red] (255:3cm)[inner sep=0pt] to[bend left=5] (53:3cm);
        \draw[color=blue] (277:3cm)[inner sep=0pt] to[bend left=5] (52:3cm);
        \draw[color=blue] (305:3cm)[inner sep=0pt] to[bend left=10] (51:3cm);
        \end{tikzpicture}
        \end{minipage}
        $\xrightarrow{\text{Step $B$}}$
    \begin{minipage}{.2 \textwidth}
        \begin{tikzpicture}[scale=.4,align=center, v/.style={draw,shape=circle, fill=black, minimum size=1mm, inner sep=0pt, outer sep=0pt}, font=\small, label distance=1pt,
            ]
        \draw (0,0) circle (3cm);
          
        \node[v] at(165:3cm){};
        \node[v] at(140:3cm){};
        \node[v] at(115:3cm){};
        \node[v] at(90:3cm){};
        \node[v] at(65:3cm){};
        \node[v] at(40:3cm){};
        \node[v] at(15:3cm){};

        \node[v] at(220:3cm){};
        \node[v] at(245:3cm){};
        \node[v] at(270:3cm){};
        \node[v] at(295:3cm){};
        \node[v] at(320:3cm){};

        \node at (77.5:3.5cm){$n$};
        \node at (52.5:3.5cm){0};
        \node at (27.5:3.5cm){1};
    
        \node at (282.5:3.7cm){$\frac{n-3}{2}$};
        \node at (257.5:3.9cm){$\frac{n-1}{2}$};
    
        \node at (160:4.2cm){$n-3$};
        \node at (135:4cm){$n-2$};
        \node at (103:3.7cm){$n-1$};
            
        \draw[color=red] (80:3cm)[inner sep=0pt] to[bend right=50] (58:3cm);
        \draw[color=blue] (103:3cm)[inner sep=0pt] to[bend right=30] (57:3cm);
        \draw[color=blue] (130:3cm)[inner sep=0pt] to[bend right=25] (56:3cm);
        \draw[color=blue] (55:3cm)[inner sep=0pt] to[bend left=20] (150:3cm);
        \draw[color=blue] (50:3cm)[inner sep=0pt] to[bend right=50] (30:3cm);
        \draw[color=blue] (50.5:3cm)[inner sep=0pt] to[bend right=35] (0:3cm);

        \draw[color=blue] (230:3cm)[inner sep=0pt] to[bend left=5] (54:3cm);
        \draw[color=blue] (255:3cm)[inner sep=0pt] to[bend left=5] (53:3cm);
        \draw[color=blue] (277:3cm)[inner sep=0pt] to[bend left=5] (52:3cm);
        \draw[color=blue] (305:3cm)[inner sep=0pt] to[bend left=10] (51:3cm);
        \end{tikzpicture}
        \end{minipage}
      \caption{(Inductive) combing algorithm for $n+1$ even, where each $A_h$ refers to Step $A_h$. We highlighted the relevant mutation at each step in red.}
  \label{combalgorithmfigeven}
  \end{figure}

  This concludes the algorithm for $n+1$ even: by induction, we have constructed the Iyama generators of $\W^2_{n+1}$ of the form $0p \times 0q$, for $1\leq p < q \leq n-1$, while following Steps $A$ and $B$ we recover $0p\times 0n$, for $p \in \{ 1,\dots n-1 \}$. See Figure \ref{combalgorithmfigeven} for a pictorial description of the mutations at each step.

Suppose now $n+1$ odd and fix the choice of grading structures on the generators that determine an endomorphism algebra $\mathcal{B}_n$ concentrated in degree 0. Define $\Gg$ to be the following sub-collection of the collection of generators $\G$:
\begin{equation*}
    \Gg:= \left\{pq \times rs \mid pq \neq 0n \text{ and } rs \neq 0n\right\}
 \end{equation*}
 and assume inductively that there exists a series of mutations on $\Gg$ that replaces $\Gg$ with the following sub-collection of the Iyama generators:
 \begin{equation*}
    \left\{0i \times 0j \mid i \neq n \text{ and } j \neq n\right\}.
\end{equation*}

For $h\in \{1,\dots,n-2\}$, in order, we define Step $h$, consisting of the mutation on $\G \setminus \Gg$ that replaces the generator $\frac{h+1}{2}\frac{2n-h-1}{2}\times 0n$ (resp. $\frac{h}{2}\frac{2n-h-2}{2}\times 0n$) with $0 \frac{h+1}{2} \times 0n$ (resp. $0 \frac{2n-h-2}{2} \times 0n$), for $h$ odd and even respectively. As in the even case, these new objects are the cones of the horizontal morphisms in the exact triangles provided by Auroux:
  \begin{equation*}
      \begin{tikzcd}
          &  0\frac{2n-h-1}{2}\times 0n  \arrow[r] & \frac{h+1}{2}\frac{2n-h-1}{2}\times 0n \arrow[d] \\
          &             &  0 \frac{h+1}{2} \times 0n \arrow[ul, "{[1]}"]
      \end{tikzcd}
  \quad
      \begin{tikzcd}
          &  \frac{h}{2}\frac{2n-h-2}{2}\times 0n  \arrow[r] & 0\frac{h}{2}\times 0n  \arrow[d] \\
          &             &  0 \frac{2n-h-2}{2} \times 0n \arrow[ul, "{[1]}"]
      \end{tikzcd}
  \end{equation*}
This concludes the inductive construction of the Iyama generators $0p\times 0q$, for $1\leq p < q \leq n$. See Figure \ref{combalgorithmfigodd} for a pictorial description of the mutations at each step.
\begin{figure}
    \begin{minipage}{.29 \textwidth}
        \centering
        \begin{tikzpicture}[scale=.4,align=center, v/.style={draw,shape=circle, fill=black, minimum size=1mm, inner sep=0pt, outer sep=0pt}, font=\small, label distance=1pt,
            ]
            \draw (0,0) circle (3cm);
        \node[v] at(140:3cm){};
        \node[v] at(115:3cm){};
        \node[v] at(90:3cm){};
        \node[v] at(65:3cm){};
        \node[v] at(40:3cm){};
        \node[v] at(15:3cm){};

        \node[v] at(15:3cm){};

        \node[v] at(232.5:3cm){};
        \node[v] at(257.5:3cm){};
        \node[v] at(282.5:3cm){};
        \node[v] at(307.5:3cm){};

        \node[v] at(180:3cm){};
  
        \node at (77.5:3.4cm){$0$};
        \node at (52.5:3.4cm){1};

        \node at (0:4cm){$\frac{h+1}{2}$};
        \node at (190:4.5cm){$\frac{2n-h-3}{2}$};
        \node at (170:4.5cm){$\frac{2n-h-1}{2}$};

        \node at (295:3.7cm){$\frac{n-2}{2}$};
        \node at (270:3.7cm){$\frac{n}{2}$};
        \node at (245:3.7cm){$\frac{n+2}{2}$};

        \node at (133:3.9cm){$n-1$};
        \node at (102.5:3.4cm){$n$};
        
        \draw[color=blue] (80:3cm)[inner sep=0pt] to[bend left=45] (100:3cm);
        \draw[color=blue] (78:3cm)[inner sep=0pt] to[bend left=30] (128:3cm);
        \draw[color=blue] (72:3cm)[inner sep=0pt] to[bend right=40] (50:3cm);

        \draw[color=blue] (168:3cm)[inner sep=0pt] to[bend right=15] (75:3cm);
        \draw[color=blue] (172:3cm)[inner sep=0pt] to[bend right=0] (5:3cm);
        \draw[color=blue] (190:3cm)[inner sep=0pt] to[bend left=5] (0:3cm);

        \draw[color=blue] (248:3cm)[inner sep=0pt] to[bend left=30] (320:3cm);
        \draw[color=blue] (250:3cm)[inner sep=0pt] to[bend left=35] (297:3cm);
        \draw[color=blue] (275:3cm)[inner sep=0pt] to[bend left=45] (295:3cm);
        
        \end{tikzpicture}
        \end{minipage}
    $\xrightarrow[h]{\text{Step}}$
        \begin{minipage}{.29 \textwidth}
            \centering
            \begin{tikzpicture}[scale=.4,align=center, v/.style={draw,shape=circle, fill=black, minimum size=1mm, inner sep=0pt, outer sep=0pt}, font=\small, label distance=1pt,
                ]
            \draw (0,0) circle (3cm);
            \node[v] at(140:3cm){};
            \node[v] at(115:3cm){};
            \node[v] at(90:3cm){};
            \node[v] at(65:3cm){};
            \node[v] at(40:3cm){};
            \node[v] at(15:3cm){};
    
            \node[v] at(15:3cm){};
    
            \node[v] at(232.5:3cm){};
            \node[v] at(257.5:3cm){};
            \node[v] at(282.5:3cm){};
            \node[v] at(307.5:3cm){};
    
            \node[v] at(180:3cm){};
      
            \node at (77.5:3.4cm){$0$};
            \node at (52.5:3.4cm){1};

            \node at (0:4cm){$\frac{h+1}{2}$};
            \node at (190:4.5cm){$\frac{2n-h-3}{2}$};
            \node at (170:4.5cm){$\frac{2n-h-1}{2}$};

            \node at (295:3.7cm){$\frac{n-2}{2}$};
            \node at (270:3.7cm){$\frac{n}{2}$};
            \node at (245:3.7cm){$\frac{n+2}{2}$};

            \node at (133:3.9cm){$n-1$};
            \node at (102.5:3.4cm){$n$};

            \draw[color=blue] (80:3cm)[inner sep=0pt] to[bend left=45] (100:3cm);
            \draw[color=blue] (78:3cm)[inner sep=0pt] to[bend left=30] (128:3cm);
            \draw[color=blue] (72:3cm)[inner sep=0pt] to[bend right=40] (50:3cm);

            \draw[color=blue] (168:3cm)[inner sep=0pt] to[bend right=15] (75:3cm);
            \draw[color=red] (74:3cm)[inner sep=0pt] to[bend right=20] (5:3cm);
            \draw[color=blue] (190:3cm)[inner sep=0pt] to[bend left=5] (0:3cm);

            \draw[color=blue] (248:3cm)[inner sep=0pt] to[bend left=30] (320:3cm);
            \draw[color=blue] (250:3cm)[inner sep=0pt] to[bend left=35] (297:3cm);
            \draw[color=blue] (275:3cm)[inner sep=0pt] to[bend left=45] (295:3cm);

            \end{tikzpicture}
            \end{minipage}
        $\xrightarrow[h+1]{\text{Step}}$
        \begin{minipage}{.28 \textwidth}
            \centering
            \begin{tikzpicture}[scale=.4,align=center, v/.style={draw,shape=circle, fill=black, minimum size=1mm, inner sep=0pt, outer sep=0pt}, font=\small, label distance=1pt,
                ]
            \draw (0,0) circle (3cm);
            \node[v] at(140:3cm){};
            \node[v] at(115:3cm){};
            \node[v] at(90:3cm){};
            \node[v] at(65:3cm){};
            \node[v] at(40:3cm){};
            \node[v] at(15:3cm){};
    
            \node[v] at(15:3cm){};

            \node[v] at(232.5:3cm){};
            \node[v] at(257.5:3cm){};
            \node[v] at(282.5:3cm){};
            \node[v] at(307.5:3cm){};
    
            \node[v] at(180:3cm){};
            
            \node at (77.5:3.4cm){$0$};
            \node at (52.5:3.4cm){1};

            \node at (0:4cm){$\frac{h+1}{2}$};
            \node at (190:4.5cm){$\frac{2n-h-3}{2}$};
            \node at (170:4.5cm){$\frac{2n-h-1}{2}$};

            \node at (295:3.7cm){$\frac{n-2}{2}$};
            \node at (270:3.7cm){$\frac{n}{2}$};
            \node at (245:3.7cm){$\frac{n+2}{2}$};

            \node at (133:3.9cm){$n-1$};
            \node at (102.5:3.4cm){$n$};
            
            \draw[color=blue] (80:3cm)[inner sep=0pt] to[bend left=45] (100:3cm);
            \draw[color=blue] (78:3cm)[inner sep=0pt] to[bend left=30] (128:3cm);
            \draw[color=blue] (72:3cm)[inner sep=0pt] to[bend right=40] (50:3cm);

            \draw[color=blue] (168:3cm)[inner sep=0pt] to[bend right=15] (76:3cm);
            \draw[color=blue] (74:3cm)[inner sep=0pt] to[bend right=20] (5:3cm);
            \draw[color=red] (190:3cm)[inner sep=0pt] to[bend right=20] (75:3cm);

            \draw[color=blue] (248:3cm)[inner sep=0pt] to[bend left=30] (320:3cm);
            \draw[color=blue] (250:3cm)[inner sep=0pt] to[bend left=35] (297:3cm);
            \draw[color=blue] (275:3cm)[inner sep=0pt] to[bend left=45] (295:3cm);

            \end{tikzpicture}
            \end{minipage}
    \caption{Pictorial description of Steps $h$ and $h+1$ of the combing algorithm on $\W^2_{n+1}$, for $n+1$ and $h$ both odd. In red, the relevant mutation at each step.}
\label{combalgorithmfigodd}
\end{figure}
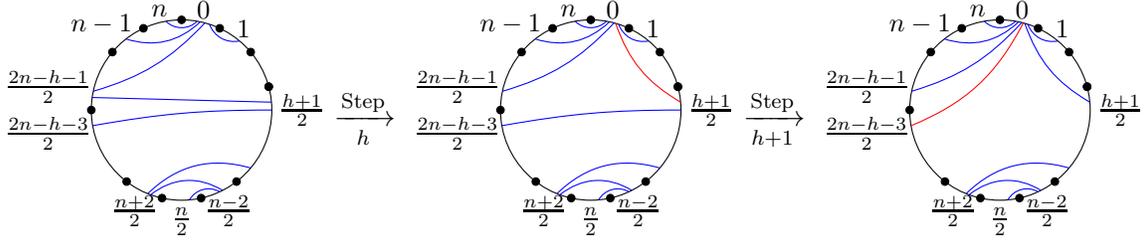
\end{proof}

\begin{remark}
    As shown in Theorem \ref{qefn}, the endomorphism algebra $\mathcal{B}_n$ is concentrated in degree 0. Fixing the unique (up to shift) choice of grading structures on the generators that realises such grading of $\mathcal{B}_n$, each step of the algorithm (for $n$ both even and odd) introduces a new generator of degree shifted by one. However, \cite[Proposition~2.2.26]{DJL} shows that we can perform a global shift so that the endomorphism algebra of the final collection of generators is again concentrated in degree 0.
\end{remark}

\subsubsection{A tilting complex for $\Gamma_n$}\label{section:quiveralgorithm}

In this section, we provide a purely algebraic proof of the quasi-equivalence of triangulated categories
\begin{equation*}
\perf(\tg_n) \simeq \perf(\Gamma_n)
\end{equation*}
which we already abstractly know to hold due to a combination of Theorem \ref{qefn} and \cite[Theorem~1]{DJL}. Our results fall under what is known as \emph{tilting theory}, which has been a broad field of study in the representation theory of finite-dimensional algebras. We use results by Rickard \cite{Rickard} on the Morita theory for derived categories, which rely on the previous foundations laid out by Happel \cite{Happel} and later generalised in \cite{Miyashita} and \cite{CPS}. Throughout this section, we will use the formalism of \cite{Keller} and \cite[I]{Hartshorne1966}. We make the following preliminary remarks. Consider the collection of indecomposable projective $\Gamma_n$-modules:
\begin{equation*}
    \Pp_{hk}:= \Gamma_n e_{(n+1-k)(n+1-h)}
\end{equation*}
where $e_{(n+1-k)(n+1-h)}$ ($h<k$) denotes the idempotent element of $\Gamma_n$ corresponding to the path of length 0 starting at the vertex $\Q_{(n+1-k)(n+1-h)}$ of the quiver with relations depicted in Figure \ref{Iyamagens2}. Viewed as complexes concentrated in degree 0, the indecomposable projectives form a distinguished collection of generators of $\perf(\Gamma_n)$. It is an immediate check (see \cite[Section~II.2]{ARS}) that, as an $A_{\infty}$-algebra concentrated in degree 0, the endomorphism algebra of this collection of $\Gamma_n$-modules is isomorphic to $\Gamma_n^{\text{op}}$. In particular, given two irreducible projective $\Gamma_n$-modules $P_{hk}$ and $P_{rs}$, there is at most one morphism from one to the other, and this exists exactly when there is a non-zero path from the vertex $\Q_{(n+1-s)(n+1-r)}$ to $\Q_{(n+1-k)(n+1-h)}$ in the quiver with relations in Figure \ref{Iyamagens2}, which is exactly when the indexes satisfy the following relations:
\begin{equation}\label{eq:existencemorphismcondition}
    h\leq r < k \leq s.
\end{equation}

Objects of $\perf(\Gamma_n)$ are bounded complexes of projective $\Gamma_n$-modules. Given an arbitrary complex $(K,d_{K})$ and an integer $t$, we denote by $(K[t],d_{K[t]})$ the complex with components $(K[t])^p:=(K)^{p+t}$ and differential $d_{K[t]}:=(-1)^{t}d_K$. We recall that, given two arbitrary complexes $(K,d_{K})$ and $(L,d_{L})$ a \emph{morphism of complexes} $u=(u^p): K \to L$ is a collection  of maps $u^p: K^p \to L^p$ compatible with the differentials:
\begin{equation*}
    d_L \circ u^p = u^{p+1} \circ d_K.
\end{equation*}
We recall whenever such morphism exists, the \emph{mapping cone} of $u$ is defined to be the complex
\begin{equation*}
    \Cone(u):=\Cone(K \xrightarrow{u} L):= L\oplus K[1]
\end{equation*}
with components $L^p \oplus K^{p+1}$ equipped with differential:
\begin{equation*}
    \begin{pmatrix}
        d_L & u\\
        0 & -d_K
    \end{pmatrix}.
\end{equation*}
For two complexes $(K,d_{K}),(L,d_{L})$ of $\Gamma_n$-modules, the \emph{morphism complex} $\Hom^{\bullet}_{\Gamma_n}(K,L)$ is defined as:
\begin{equation*}
    \Hom^{\bullet}_{\Gamma_n}(K,L)^q := \prod_{p\in \Z} \Hom_{\Gamma_{n}}(K^p,L^{p+q})
\end{equation*}
with differential given by:
\begin{equation*}
    d\varphi = d_{L} \circ \varphi + (-1)^{q} \varphi \circ d_{K}
\end{equation*}
for $\varphi \in \Hom^{\bullet}_{\Gamma_n}(K,L)^q$. The $\Ext$-group $\Ext^t(K,L)$ is defined to be the $t^{th}$ cohomology of the corresponding morphism complex. There is a canonical isomorphism:
\begin{equation*}
    \text{Ext}^t(K,L)\cong \text{Hom}_{\perf(\Gamma_n)}(K,L[t])
\end{equation*}
for any $t \in \Z$ (\cite[Chapter~I, Section~6]{Hartshorne1966}). We will now define a special collection of objects in $\perf(\Gamma_n)$.

\begin{definition}[Complexes $\K_{ij\ell m}$]\label{def:tiltingcomplex}
    Fix $n > 2$. For $n$ odd (resp. even) and for indexes $i,j,\ell,m$, satisfying the following:
    \begin{equation*}
        \begin{cases}
            0\leq \ell \leq i \leq \frac{n-3}{2} &(\text{resp. } 0\leq \ell \leq i \leq \frac{n-2}{2}) \\
            j\in\{n-2-i,n-1-i\} &(\text{resp. } j\in\{n-1-i,n-i\}) \\
            m\in\{n-2-\ell,n-1-\ell\} &(\text{resp. } m\in\{n-1-\ell,n-\ell\}) \\
        \end{cases}
    \end{equation*}
    with $j\leq m\leq n-1$, not both $\ell=i$ and $j=m$, we define $\K_{ij\ell m}$ to be the following complex of projective $\Gamma_n$-modules, concentrated in degrees $-2,-1,0,1$:
    \begin{equation}\label{eq:tiltingcomplexdef}
        \K_{ij\ell m}:= P_{\ell i} \to P_{\ell j} \to P_{im} \to P_{jm}
    \end{equation}
   with the additional convention that $P_{hh}=P_{0h}=0$ for any index $h$. Each differential is given by the unique morphism between indecomposable projectives, and it vanishes exactly when either source or target is zero.
\end{definition}

For a fixed $n>2$, define $T_n$ to be the complex $T_n:=\bigoplus \K_{ij\ell m}$, where the sum ranges over all the objects defined in Definition \ref{def:tiltingcomplex}. The main result of this section is the following.

\begin{proposition}\label{prop:equivperf}
    The equivalence of triangulated categories
    \begin{equation*}
        \perf(\Gamma_n) \simeq \perf(\tg_n)
    \end{equation*}
    is explicitly realised by the functor $\perf(\Gamma_n) \to \perf(\tg_n)$ sending $T_n$ to $\tg_n$.
\end{proposition}

In order to prove this equivalence, we require some preliminary lemmas.

\begin{lemma}\label{lemma:conescomplex}
    The complexes appearing in Definition \ref{def:tiltingcomplex} arise as iterated mapping cones of projective $\Gamma_n$-modules. In particular, for appropriate $i,j,\ell,m$ and viewing each projective $\Gamma_n$-modules as a complex concentrated in degree 0, we have the following isomorphisms in $\perf(\Gamma_n)$:
\begin{equation}\label{eq:conescomplex}
   \K_{ij\ell m} \cong \Cone\left( \Cone\left( P_{\ell i} \xrightarrow{u} P_{\ell j}\right) \xrightarrow{z} \Cone\left( P_{im} \xrightarrow{v} P_{jm}\right)[1] \right)
\end{equation}
where each $u,v,z$ is the unique (possibly zero, if either the source or the target is zero) morphism between projective $\Gamma_n$-modules, and $[1]$ denotes a shift by $1$ of the complex. 
\end{lemma}

\begin{proof}
    This follows immediately from the definition of the mapping cone of a morphism of complexes. The shift in degree is necessary for $z$ to be a morphism of complexes.
\end{proof}

Denote by $\text{add}(T_n)$ the closure of $T_n$ in $\perf(\Gamma_n)$ under taking direct summands of finite directed sums, and by $\langle \text{add}(T_n)\rangle$ the smallest triangulated subcategory of $\perf(\Gamma_n)$ containing $\text{add}(T_n)$. We recall that a \emph{tilting complex} $T$ for $\Gamma_n$ is defined (\cite[Section~6]{Rickard}) to be an object of $\perf(\Gamma_n)$ satisfying the following conditions:
\begin{enumerate}[i.]
    \item $\text{Hom}_{\perf(\Gamma_n)}(T,T[t])=0$ for $t\neq 0$;
    \item $\text{add}(T)$ generates $\perf(\Gamma_n)$ as a triangulated category.
\end{enumerate}

\begin{lemma}\label{lemma:tiltingcomplex}
    The complex $T_n$ is a tilting complex of $\perf(\Gamma_n)$.
\end{lemma}

Lemma \ref{lemma:tiltingcomplex} amounts Lemma \ref{lemma:tiltingcomplex1} and \ref{lemma:tiltingcomplex2}.

\begin{lemma}\label{lemma:tiltingcomplex1}
    $T_n$ is an exceptional object in $\perf(\Gamma_n)$.
\end{lemma}

\begin{proof}
Fix $n$ and $K:=\K_{ij \ell m}$, $L:=\K_{hkrs}$ two arbitrary complexes defined in Definition \ref{def:tiltingcomplex}, and recall that the indexes satisfy the following relations:
\begin{equation}\label{eq:existencomplexcondition}
    \begin{cases}
    0\leq \ell \leq i, \quad  j\leq m, \quad 0 \leq r \leq h, \quad k\leq s, \\
    n-2-i \leq j \leq n-1-i \text{ for $n$ odd } \quad \text{(resp. $n-1-i \leq j \leq n-i$ for $n$ even)}, \\
    n-2-\ell  \leq m \leq n-1-\ell  \text{ for $n$ odd } \quad \text{(resp. $n-1-\ell  \leq m \leq n-\ell $ for $n$ even)}, \\
    n-2-h \leq k \leq n-1-h \text{ for $n$ odd } \quad \text{(resp. $n-1-h \leq k \leq n-h$ for $n$ even)}, \\
    n-2-r \leq s \leq n-1-r \text{ for $n$ odd } \quad \text{(resp. $n-1-r \leq s \leq n-r$ for $n$ even)}.
\end{cases}
\end{equation}

We prove that $\{\K_{ij \ell m}\}$ is an exceptional collection, i.e.\ that $Ext^t(K,L)$ vanishes for any $t\neq0$. By conditions (\ref{eq:existencemorphismcondition}) and (\ref{eq:existencomplexcondition}), one can verify that the morphism complex $\text{Hom}^{\bullet}:=\text{Hom}_{\Gamma_n}^{\bullet}(K,L)$ is concentrated in degrees $-1,0,1,2$. We compute its cohomology in degrees $-1,1$ and $2$ and verify that it  vanishes.

By (\ref{eq:existencemorphismcondition}) and (\ref{eq:existencomplexcondition}), the only possible non-zero generator of $\text{Hom}^{-1}$ is the unique morphism $\pi: P_{im}\to P_{rk}$, whenever this exists. If it does, one can verify using (\ref{eq:existencemorphismcondition}) that there also exist (unique and non-zero) morphisms $\beta: P_{\ell j}\to P_{rk}$ and $\gamma: P_{im}\to P_{hs}$, generators of $\text{Hom}^{0}$, such that $d\pi=\beta + \gamma$. In this case, the cohomology of the morphism complex in degree $-1$ vanishes. This is also trivially true if $\pi$ vanishes.

In order to compute $\text{Ext}^1(K,L)$, we distinguish the following cases:
\begin{enumerate}[i.]
    \item There are no morphisms $\xi: P_{\ell i}\to P_{hs}$ and $\eta: P_{\ell j}\to P_{ks}$;
    \item There is no morphism $\xi: P_{\ell i}\to P_{hs}$ and there is a unique morphism $\eta: P_{\ell j}\to P_{ks}$;
    \item There is no morphism $\eta: P_{\ell j}\to P_{ks}$ and there is a unique morphism $\xi: P_{\ell i}\to P_{hs}$.
\end{enumerate}

Computations, following immediately from conditions (\ref{eq:existencemorphismcondition}) and (\ref{eq:existencomplexcondition}), show that $\xi: P_{\ell i}\to P_{hs}$ and $\eta: P_{\ell j}\to P_{ks}$ cannot co-exist. If i.\ holds, computations following immediately from (\ref{eq:existencemorphismcondition}) and (\ref{eq:existencomplexcondition}) show that the following statements are true:
\begin{itemize}
    \item If there exists a morphism $\lambda: P_{\ell i}\to P_{rk}$, then there exists a generator $\alpha: P_{\ell i}\to P_{rh}$ of $\text{Hom}^{0}$;
    \item If there exists a morphism $\nu: P_{im}\to P_{ks}$, then there exists a generator $\delta: P_{jm}\to P_{ks}$ of $\text{Hom}^{0}$;
    \item If there exists a morphism $\mu: P_{\ell j}\to P_{hs}$, then at least one of the following is true:
    \begin{itemize}
        \item There exist morphisms $\alpha: P_{\ell i}\to P_{rh}$ and $\beta: P_{\ell j}\to P_{rk}$, generators of $\text{Hom}^{0}$;
        \item There exist morphisms $\gamma: P_{im}\to P_{hs}$ and $\delta: P_{jm}\to P_{ks}$, generators of $\text{Hom}^{0}$;
        \item There is a morphism $\beta: P_{\ell j}\to P_{rk}$ and there is no morphism $P_{\ell i}\to P_{rk}$;
        \item There is a morphism $\gamma: P_{im}\to P_{hs}$ and there is no morphism $P_{im}\to P_{ks}$.
    \end{itemize}
\end{itemize}

Similarly, if ii.\ holds, by (\ref{eq:existencemorphismcondition}) and (\ref{eq:existencomplexcondition}) the following statements are true:

\begin{itemize}
    \item If there exists a morphism $\lambda: P_{\ell i}\to P_{rk}$, then there exists a generator $\alpha: P_{\ell i}\to P_{rh}$ of $\text{Hom}^{0}$;
    \item If there exist morphisms $\mu: P_{\ell j}\to P_{hs}$ and $\nu: P_{im}\to P_{ks}$, then there exists a generator $\gamma: P_{im}\to P_{hs}$ of $\text{Hom}^{0}$;
\end{itemize}

Furthermore, if iii.\ holds, by (\ref{eq:existencemorphismcondition}) and (\ref{eq:existencomplexcondition}) the following statements are true:

\begin{itemize}
    \item If there exists a morphism $\nu: P_{im}\to P_{ks}$, then there exists a generator $\delta: P_{jm}\to P_{ks}$ of $\text{Hom}^{0}$;
    \item If there exist morphisms $\lambda: P_{\ell i}\to P_{rk}$ and $\mu: P_{\ell j}\to P_{hs}$, then there exists a generator $\beta: P_{\ell j}\to P_{rk}$ of $\text{Hom}^{0}$.
\end{itemize}

The above cases exhaust all possible computational cases for $\text{Ext}^1(K,L)$, and directly imply that the latter vanishes. Finally, if either morphism $\xi: P_{\ell i}\to P_{hs}$ or $\eta: P_{\ell j}\to P_{ks}$ exists, computations following immediately from (\ref{eq:existencemorphismcondition}) and (\ref{eq:existencomplexcondition}) show that there exists a generator $\mu: P_{\ell j}\to P_{hs}$ of $\text{Hom}^{1}$. It directly follows that $\text{Ext}^2(K,L)$ also vanishes (it trivially does if neither $\xi$ nor $\eta$ exist) and concludes the proof.\end{proof}

\begin{lemma}\label{lemma:tiltingcomplex2}
    $\textnormal{add}(T_n)$ generates $\perf(\Gamma_n)$ as a triangulated category.
\end{lemma}

\begin{proof}
It suffices to show that the distinguished collection of projective $\Gamma_n$-modules is contained in $\langle\textnormal{add}(T_n)\rangle$. We construct an iterative proof of the claim, in the following way: fixing $n$ odd (resp. even) we prove, inductively on the index $\ell$ and for $m\in\{n-2-\ell,n-1-\ell\}$ (resp. $m\in\{n-1-\ell,n-\ell\}$), that the projectives $P_{\ell i},P_{\ell j},P_{im},P_{jm}$, for all $i\geq \ell$, $j\in\{n-2-i,n-1-i\}$ (resp. $j\in\{n-1-i,n-i\}$), $j\leq m$, belong to $\langle\textnormal{add}(T_n)\rangle$. Each iteration (for fixed $\ell$) is itself proven inductively on the index $i\geq \ell$ and $j\in\{n-2-i,n-1-i\}$ (resp. $j\in\{n-1-i,n-i\}$). Each iteration terminates after a finite (depending on $n$ and fixed $\ell$) number of steps at $i=\frac{n-3}{2}, j=\frac{n-1}{2}$ (resp. at $i=\frac{n-2}{2}, j=\frac{n}{2}$). The induction on $\ell$ also terminates after a finite (depending on $n$) number of steps, when one reaches $\ell=\frac{n-3}{2}, m=\frac{n+1}{2}$ (resp. $\ell=\frac{n-2}{2}, m=\frac{n+2}{2}$) and concludes the proof.
    
We first fix $\ell=0$. Assuming first $n$ odd and fixing $m=n-1$, $P_{(n-2)(n-1)}$ is, up to a shift in degree, isomorphic to the complex $\K_{0(n-2)0(n-1)}$ defined in (\ref{eq:tiltingcomplexdef}), so it is a summand in $T_n$. We assume, inductively on $i$, that the projectives $P_{i(n-1)}$ and $P_{j(n-1)}$, for $j=n-2-i$, are contained in $\langle\textnormal{add}(T_n)\rangle$. From Lemma \ref{lemma:conescomplex} we know that $\K_{(i+1)j0(n-1)}$ is isomorphic to the mapping cone of the unique morphism $P_{(i+1)(n-1)}\to P_{j(n-1)}$, therefore $P_{(i+1)(n-1)}$ is contained in $\langle\textnormal{add}(T_n)\rangle$. Moreover, $\K_{(i+1)(j-1)0(n-1)}$ is isomorphic to the mapping cone of the unique morphism $P_{(i+1)(n-1)}\to P_{(j-1)(n-1)}$, so $P_{(j-1)(n-1)}$ also belongs to $\langle\textnormal{add}(T_n)\rangle$. This proves that $P_{(i+1)(n-1)}$ and $P_{(j-1)(n-1)}$, for $j-1=n-2-(i+1)$ belong to $\langle\textnormal{add}(T_n)\rangle$. Similarly, now fixing $m=n-2$ for $n$ odd, $P_{1(n-2)}$ is isomorphic to the complex $\K_{1(n-2)0(n-1)}$, and we can inductively on $i$ prove that all $P_{i(n-2)}$ and $P_{j(n-2)}$ ($j=n-2-i$), are contained in $\langle\textnormal{add}(T_n)\rangle$. This concludes the base case for the induction on $\ell$ for $n$ odd. The even case is completely analogous, with the only difference being that the base case only consists of $\ell=0$ and $m=n-1$.

Assume now, inductively on $\ell$, that the projectives $P_{\ell i},P_{\ell j},P_{im}, P_{jm}$ belong to $\langle\textnormal{add}(T_n)\rangle$ for $m\in\{n-2-\ell,n-1-\ell\}$ (resp. $m\in\{n-1-\ell,n-\ell\}$) and for all $i\geq \ell$ and $j\in\{n-2-i,n-1-i\}$ (resp. $j\in\{n-1-i,n-i\}$), $j\leq m$. Fixing $m=n-2-\ell$ (resp. $m=n-1-\ell$), and similarly to the base case, we can prove inductively on $i$ that $P_{(\ell+1) i}$ and $P_{(\ell+1) j}$ belong to $\langle\textnormal{add}(T_n)\rangle$ for all $i \geq \ell+1$, $j\in\{n-2-i,n-1-i\}$ (resp. $j\in\{n-1-i,n-i\}$); to show this, we use the isomorphisms (\ref{eq:conescomplex}) for $\K_{ij(\ell+1)m}$. Consequently, using the isomorphisms (\ref{eq:conescomplex}) for $\K_{ij(\ell+1)(m-1)}$, one can prove inductively on $i$ that $P_{i(m-1)}$ and $P_{j(m-1)}$ belong to $\langle\textnormal{add}(T_n)\rangle$ for all $i \geq \ell$, $j\in\{n-2-i,n-1-i\}$ (resp. $j\in\{n-1-i,n-i\}$), which concludes the proof.
\end{proof}

This concludes the proof of Lemma \ref{lemma:tiltingcomplex}. The following Lemma is the last property needed to prove Proposition \ref{prop:equivperf}.

\begin{lemma}\label{lemma:isomT}
    The endomorphism algebra of $T_n$ over $\perf(\Gamma_n)$ is isomorphic to $\tg_n$.
\end{lemma}

\begin{proof}
    Fix $n$ and fix an arbitrary complex $\K_{ij\ell m}$ defined in (\ref{eq:tiltingcomplexdef}). It is clear that there exists a morphism from $\K_{ij\ell m}$ to another complex defined in (\ref{eq:tiltingcomplexdef}) whenever the latter is one of the following:
    \begin{equation*}
        \K_{(i+1)j\ell m} \qquad \qquad \K_{i(j+1)\ell m} \qquad \qquad \K_{ij(\ell+1) m} \qquad \qquad \K_{ij\ell (m+1)}
    \end{equation*}
    Furthermore, the following facts hold:
    \begin{enumerate}[i.]
        \item The composition of the morphisms:
            \begin{alignat*}{3}
                u=&(u^{-2},id,u^{0},id):& \K_{ij\ell m}\to &\K_{(i+1)j\ell m} \\
                v=&(id, id,v^{0},v^{1}):& \K_{(i+1)j\ell m} \to &\K_{(i+1)j\ell (m+1)}
            \end{alignat*}
            (whenever all above complexes exist, and where each $u^p, v^{p}$ is the unique morphism between projective $\Gamma_n$-modules and $id$ is the identity) gives rise to a morphism of complexes:
            \begin{equation*}
                z=(u^{-2},id,v^{0}u^{0},v^{1}):\K_{ij\ell m} \to \K_{(i+1)j\ell (m+1)}.
            \end{equation*}
            Moreover, the composition of the unique morphisms:
            \begin{align*}
                a=(id,id,a^{0},a^{1}): &\K_{ij\ell m}\to \K_{ij\ell (m+1)} \\
                b=(b^{-2}, id,b^{0},id): &\K_{ij\ell (m+1)} \to \K_{(i+1)j\ell (m+1)}
            \end{align*}
            (whenever all above complexes exist, and where each $a^p, b^{p}$ is the unique morphism between projective $\Gamma_n$-modules) gives rise to a morphism:
            \begin{equation*}
                c=(b^{-2},id,b^{0}a^{0},a^{1}):\K_{ij\ell m} \to \K_{(i+1)j\ell (m+1)}.
            \end{equation*}
            As $b^{-2}=u^{-2}$, $a^{1}=v^{1}$ and $b^{0}a^{0}=v^{0}u^{0}$ (following from the uniqueness of the morphism $P_{im} \to P_{(i+1)(m+1)}$), the two compositions give rise to the same morphism of complexes. Uniqueness of the morphism $z=c$ follows from uniqueness of the morphisms between projective $\Gamma_n$-modules;
        \item Similarly to i., the morphism of complexes $\K_{ij\ell m} \to \K_{(i+1)j(\ell+1)m}$ (when the latter complex exists) is unique and arises as both the composition $\K_{ij\ell m}\to \K_{(i+1)j\ell m} \to \K_{(i+1)j(\ell+1)m}$ and $\K_{ij\ell m} \to \K_{ij(\ell+1)m} \to \K_{(i+1)j(\ell+1)m}$;
        \item Similarly, the morphism of complexes $\K_{ij\ell m} \to \K_{i(j+1)\ell (m+1)}$ (when this exists) is unique and arises as both the composition $\K_{ij\ell m}\to \K_{i(j+1)\ell m} \to \K_{i(j+1)\ell (m+1)}$ and $\K_{ij\ell m} \to \K_{ij\ell (m+1)} \to \K_{i(j+1)\ell (m+1)}$;
        \item Similarly, the morphism of complexes $\K_{ij\ell m} \to \K_{i(j+1)(\ell+1)m}$ (when this exists) is unique and arises as both the composition $\K_{ij\ell m}\to \K_{i(j+1)\ell m} \to \K_{i(j+1)(\ell+1)m}$ and $\K_{ij\ell m} \to \K_{ij(\ell+1)m} \to \K_{i(j+1)(\ell+1)m}$.
    \end{enumerate}
    Finally, we claim that there is no morphism from $\K_{ij \ell m}$ to any other complex $\K_{hkrs}$ ($0\leq r \leq h < k \leq s \leq n-1$) defined in (\ref{eq:tiltingcomplexdef}) and not listed above. If this holds, it is clear that $\text{End}(T_n)\cong \tg_n$ as (ungraded) algebras.
    
    To prove the claim, suppose that there is a morphism $u=(u^{-2},u^{-1},u^{0},u^{1}): \K_{ij \ell m} \to \K_{hkrs}$, with $r>\ell +1$ (resp. $h>i+1$, $s>m+1$, $k>j+1$). By (\ref{eq:existencomplexcondition}), in particular this implies $s<m$ (resp. $k<j$, $r<\ell$, $h<i$) and therefore $u^{0}=u^{1}=0$ (resp. $u^{-1}=u^{1}=0$, $u^{-2}=u^{-1}=0$, $u^{-2}=u^{0}=0$). By (\ref{eq:existencomplexcondition}) and (\ref{eq:existencemorphismcondition}) one can check that it directly follows that $u^{-2}=u^{-1}=0$ (resp. $u^{-2}=u^{0}=0$, $u^{0}=u^{1}=0$, $u^{-1}=u^{1}=0$), and $u$ cannot exist. This exhausts all possible morphisms.\end{proof}

\begin{proof}[Proof of Proposition \ref{prop:equivperf}]
    This follows from results of Rickard \cite[Theorems~2.12, 6.4]{Rickard}, and from Lemma \ref{lemma:tiltingcomplex} and \ref{lemma:isomT}.
\end{proof}

\subsubsection{The vanishing cycles on the regular fibre}\label{section:vcyclesalgorithm}
In Section \ref{ACconf} we gave a description of $\Sigma_n$, the Milnor fibre of $f_n$, together with a collection of vanishing cycles (associated to a distinguished basis of vanishing paths), whose corresponding Lefschetz thimbles generate the Fukaya-Seidel category $\F(f_n)$. In this section we perform a series of mutations on the given basis of vanishing paths; each mutation consists in a Hurwitz move on the distinguished collection, and it gives rise to a new basis. The effect of each mutation on the vanishing cycles is that of a symplectic Dehn twist, as prescribed in Section \ref{section:background}.

Let $\{V_{I,J}\}$ be the vanishing cycles given in Section \ref{ACconf} and $\{\DAC_{I,J}\}$ their associated Lefschetz thimbles. Under the equivalence in Theorem \ref{qefn}, each $\DAC_{I,J}$ corresponds to the generator $ij \times \ell m$ of $\W^2_n$, where $I=|i-j|$ and $J=|\ell -m|$. Following this notation, we redefine each cycle as $V_{ij,\ell m}:=V_{I,J}$ and we denote by $\gamma_{ij,\ell m}$ its corresponding path. We now construct a new collection of vanishing paths $\gamma_{I,J}$ and cycles $\{\Iyama_{I,J}\}$ (which we call \emph{Iyama vanishing cycles}), whose associated Lefschetz thimbles $\{\DIy_{I,J}\}$ are those corresponding, under the equivalence $\F(f_n) \simeq \W^2_n$ of Corollary \ref{cor:maineq}, to the generators $\{0I \times 0J\}$ of $\W^2_n$. We do so in two steps: in Proposition \ref{orderIyamagenerators} we describe the mutations on the paths, while in Proposition \ref{prop:homeofibre} we track their geometric effect on the collection of cycles.

\begin{proposition}\label{orderIyamagenerators}
   The vanishing paths $\{\gamma_{I,J}\}$ are obtained from $\{\gamma_{ij,\ell m}\}$ via a series of mutations arising from Hurwitz moves. Moreover, their total order is the lexicographic one:
    \begin{equation*}
        \gamma_{I,J} <\gamma_{K,L} \quad \text{if and only if} \quad (J < L) \text{ or } (J=L, I < K).
    \end{equation*}
\end{proposition}

\begin{remark}\label{rk:notation}
    The series of mutations involved in the above statement correspond to those prescribed in Section \ref{algorithmdisk}, under the equivalence $\F(f_n)\cong \W^2_n$ of Corollary \ref{cor:maineq}. However, a Hurwitz-type move on a pair of paths requires them to be consecutively ordered. For this reason, while defining mutations on a given collection, we have to keep track of the total order of the objects, which is not done in the algorithm given in Proposition \ref{combalgorithm}.
\end{remark}

\begin{proof}
    Assume first $n+1$ even. Denote by $\G$ the distinguished collection of vanishing paths prescribed in Section \ref{ACconf}. Further denote ${\gamma_{ij,\ell m}}$ in $\G$ by $\gamma_{ij,\ell m}\n$, $\gamma_{ij,\ell m}\s$ and $\gamma_{ij,\ell m}\p$ respectively depending on whether they are of negative, saddles or positive type. We choose a further perturbation of the Morsification of $f_n$ that separates all critical values (Remark \ref{rk:furtherperturbation}) so that the order (after applying Step 0 of the algorithm in the proof of Proposition \ref{combalgorithm}, which can be interpreted as a relabelling) is the following:
    \begin{itemize}
        \item The paths $\gamma_{ij,(n-1)n}\n$ come after the other paths of same type, and are themselves ordered for increasing $j$, with the first one being $\gamma_{\frac{n-3}{2}\frac{n-1}{2},(n-1)n}\n$ and the last one being $\gamma_{0(n-2),(n-1)n}\n$;
        \item The saddles $\gamma_{ij,0(n-1)}\s$ all come after the other saddles, and are ordered for increasing $j$;
        \item  The saddles $\gamma_{ij,(n-1)n}\s$ all come after the remaining saddles, and are ordered for increasing $i$;
        \item The positive paths $\gamma_{ij,0(n-1)}\p$ come after the other paths of same type, and are ordered for increasing $i$;
        \item The remaining paths are ordered arbitrarily.
    \end{itemize}

    The positive paths $\gamma_{ij,\ell m}\p$ ($\ell m \neq 0(n-1)$) are disjoint from the saddles $\gamma_{ij,(n-1)n}\s$ and $\gamma_{ij,0(n-1)}\s$, and the paths $\gamma_{ij,(n-1)n}\n$ are disjoint from the remaining saddles and positive paths (Figure \ref{disjointcycles}, left); we perform appropriate Hurwitz moves on the paths, so that $\gamma_{ij,\ell m}$ (for $\ell m \neq 0(n-1)$ and $\ell m \neq (n-1)n$) all come before $\gamma_{ij,(n-1)n}$ and $\gamma_{ij,0(n-1)}$ in the total ordering of this distinguished collection of paths. Denote by $\Gg$ the sub-collection of vanishing paths such that:
    \begin{equation*}
       \G \setminus \Gg=\left\{{\gamma_{ij,(n-1)n}}\right\} \cup \left\{{\gamma_{ij,0(n-1)}}\right\}.
    \end{equation*}
    Assume inductively that we have constructed a series of mutations on $\Gg$ (acting by Hurwitz moves) that replaces $\Gg$ with the following collection of paths, whose order is the lexicographic one:
    \begin{equation*}
        \left\{\gamma_{I,J} \mid I,J \neq n-1 \text{ and } I,J \neq n\right\}.
    \end{equation*}

    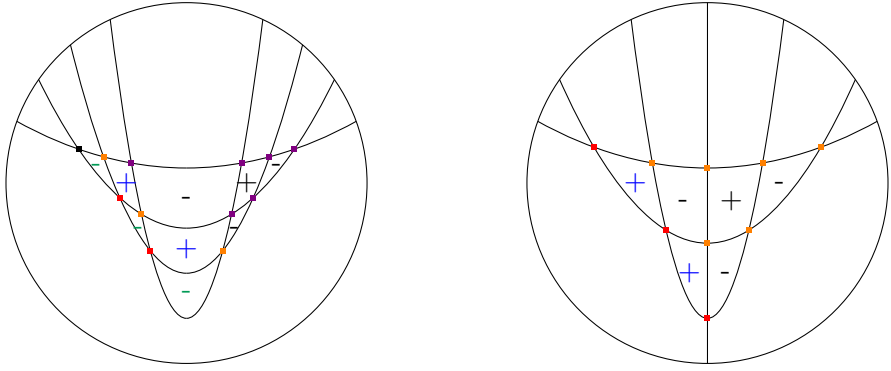
\begin{figure}
        \centering
\begin{minipage}{.45 \textwidth}
\centering
\begin{tikzpicture}[scale=.8,v/.style={draw,shape=circle, fill=black, minimum size=1.2mm, inner sep=0pt, outer sep=0pt}]

\draw (0,0) circle (3cm);

\draw[name path = p8] (0,.25) parabola (20:3cm);
\draw[name path = p1]  (0,.25) parabola (160:3cm);

\draw[name path = p7] (0,-.75) parabola (35:3cm);
\draw[name path = p2] (0,-.75) parabola (145:3cm);

\draw[name path = p6] (0,-1.5) parabola (50:3cm);
\draw[name path = p3] (0,-1.5) parabola (130:3cm);

\draw[name path = p5] (0,-2.25) parabola (65:3cm);
\draw[name path = p4] (0,-2.25) parabola (115:3cm);

\path [name intersections = {of = p1 and p2, by = A}];
\node [fill = black, inner sep = 1pt] at (A) {};

\path [name intersections = {of = p2 and p3, by = B}];
\node [fill = red, inner sep = 1pt] at (B) {};

\path [name intersections = {of = p3 and p4, by = C}];
\node [fill = red, inner sep = 1pt] at (C) {};

\path [name intersections = {of = p1 and p3, by = D}];
\node [fill = orange, inner sep = 1pt] at (D) {};

\path [name intersections = {of = p2 and p4, by = E}];
\node [fill = orange, inner sep = 1pt] at (E) {};

\path [name intersections = {of = p5 and p6, by = F}];
\node [fill = orange, inner sep = 1pt] at (F) {};

\path [name intersections = {of = p1 and p4, by = G}];
\node [fill = violet, inner sep = 1pt] at (G) {};

\path [name intersections = {of = p5 and p7, by = H}];
\node [fill = violet, inner sep = 1pt] at (H) {};

\path [name intersections = {of = p5 and p8, by = I}];
\node [fill = violet, inner sep = 1pt] at (I) {};

\path [name intersections = {of = p6 and p7, by = J}];
\node [fill = violet, inner sep = 1pt] at (J) {};

\path [name intersections = {of = p6 and p8, by = K}];
\node [fill = violet, inner sep = 1pt] at (K) {};

\path [name intersections = {of = p7 and p8, by = L}];
\node [fill = violet, inner sep = 1pt] at (L) {};


\node[color=ForestGreen] at (0,-1.8){-};

\node[color=blue] at (0,-1.1){+};
\node at (.8,-.75){-};
\node[color=ForestGreen] at (-.8,-.75){-};

\node at (0,-.25){-};
\node at (1,0){+};
\node[color=blue] at (-1,0){+};
\node at (1.5,.3){-};
\node[color=ForestGreen] at (-1.5,.3){-};

\end{tikzpicture}
\end{minipage}
\centering
\begin{minipage}{.45 \textwidth}
\centering
\begin{tikzpicture}[scale=.8]
    \draw (0,0) circle (3cm);

\draw[name path = p7] (0,.25) parabola (20:3cm);
\draw[name path = p1] (0,.25) parabola (160:3cm);

\draw[name path = p6] (0,-1) parabola (35:3cm);
\draw[name path = p2] (0,-1) parabola (145:3cm);

\draw[name path = p5] (0,-2.25) parabola (65:3cm);
\draw[name path = p3] (0,-2.25) parabola (115:3cm);

\draw[name path = p4] (90:3cm) to (-90:3cm);

\path [name intersections = {of = p1 and p2, by = A}];
\node [fill = red, inner sep = 1pt] at (A) {};

\path [name intersections = {of = p2 and p3, by = B}];
\node [fill = red, inner sep = 1pt] at (B) {};

\path [name intersections = {of = p3 and p4, by = C}];
\node [fill = red, inner sep = 1pt] at (C) {};

 \path [name intersections = {of = p1 and p3, by = D}];
 \node [fill = orange, inner sep = 1pt] at (D) {};

 \path [name intersections = {of = p2 and p4, by = E}];
 \node [fill = orange, inner sep = 1pt] at (E) {};

 \path [name intersections = {of = p5 and p6, by = F}];
 \node [fill = orange, inner sep = 1pt] at (F) {};

 \path [name intersections = {of = p1 and p4, by = G}];
 \node [fill = orange, inner sep = 1pt] at (G) {};

 \path [name intersections = {of = p5 and p7, by = H}];
 \node [fill = orange, inner sep = 1pt] at (H) {};

 \path [name intersections = {of = p6 and p7, by = I}];
 \node [fill = orange, inner sep = 1pt] at (I) {};


\node at (.3,-1.5){-};
\node[color=blue] at (-.3,-1.5){+};

\node at (.4,-.3){+};
\node at (1.2,0){-};
\node at (-.4,-.3){-};
\node[color=blue] at (-1.2,0){+};

\end{tikzpicture}
\end{minipage}
\caption{The $r$-divides associated to the Milnor fibre and cycles $\{V_{I,J}\}$ for (left) $n+1$ even and (right) $n+1$ odd. (left) In red and orange respectively, the saddles $V_{ij,(n-1)n}\s$ and $V_{ij,0(n-1)}\s$; in black, the saddle $V_{0(n-1),n(n-1)}\s$, in green, the negative cycles $V_{ij,(n-1)n}\n$ and in blue, the positive cycles $V_{ij,0(n-1)}\p$. (right) In red and blue, the cycles $V_{ij,0n}$.}
\label{disjointcycles}
    \end{figure}

The cycle associated to $\gamma_{0(n-1),(n-1)n}\s$ is disjoint from any of the positive cycles: perform appropriate Hurwitz moves so that it comes as the last ordered path. Retracing Step $A$ of the algorithm given in Proposition \ref{combalgorithm}, we now construct a series of mutations on $\G \setminus \Gg$, and we call each iterated move ``Step $A_h$'', for $1\leq h \leq n-3$. For $h$ odd, we define the following Hurwitz move:
\begin{equation*}
    \left(\gamma_{0\frac{2n-h-3}{2}, pq},\gamma_{\frac{h+1}{2}\frac{2n-h-3}{2}, pq}\right) \mapsto  \left(\tau_{\gamma_{0\frac{2n-h-3}{2}, pq}}\gamma_{\frac{h+1}{2}\frac{2n-h-3}{2}, pq},\gamma_{0\frac{2n-h-3}{2}, pq}\right)
\end{equation*}
for  $pq\in \{0(n-1),(n-1)n\}$. Denote the paths resulting from this move as:
\begin{equation*}
    \gamma_{0\frac{h+1}{2},pq}\z:=\tau_{\gamma_{0\frac{2n-h-3}{2}, pq}\y} \gamma_{\frac{h+1}{2}\frac{2n-h-3}{2}, pq}\z
\end{equation*}
for appropriate $\textstyle{\circ},\scriptstyle{\triangle} \in \{-, \boldsymbol{\cdot}, +\}$. Similarly, for $h$ even, appropriate Hurwitz moves allow us to define:
\begin{equation*}
    \gamma_{0\frac{2n-h-4}{2},pq}\z:=\tau^{-1}_{\gamma_{0\frac{h}{2},pq}\y}\gamma_{\frac{h}{2}\frac{2n-h-4}{2},pq}\z.
\end{equation*}
Additionally, after Step $A_h$ (for $h<n-4$), and before Step $A_{h+1}$, we perform the following mutations:
\begin{itemize}
    \item If $h\equiv 1 \mod 2$, and for all $i=\frac{h+3}{2},\dots, \frac{n-3}{2}$, a mutation of the paths $\gamma_{0\frac{2n-h-3}{2},pq}\y$ around the paths $\gamma_{i(n-1-i),pq}\z$;
    \item If $h\equiv 0 \mod 2$, and for all $i=\frac{h+2}{2},\dots, \frac{n-3}{2}$, a mutation of the paths $\gamma_{0\frac{h}{2}, pq}\y$ around the paths $\gamma_{i(n-2-i), pq}\z$;
\end{itemize}
for appropriate $\textstyle{\circ},\scriptstyle{\triangle} \in \{-, \boldsymbol{\cdot}, +\}$. One can check that these mutations consists in Hurwitz moves on paths corresponding to pairwise disjoint cycles, and guarantee that Step $A_h$ is a Hurwitz move on consecutively ordered paths.

\begin{remark}
    To help in keeping track of the iterated mutations, we have carefully chosen our notation so that each Lefschetz thimble associated to the path $\gamma_{pq,rs}$ \textbf{at each step} is the one corresponding to $pq \times rs$ under the equivalence $\F(f_n)\simeq \W^2_n$.
\end{remark}

We can now drop the sign superscript; after Step $A_{n-3}$, the paths are ordered as:
\begin{equation*}
    \dots < \left\{\gamma_{0i,(n-1)n}\right\} < \left\{\gamma_{0j,(n-1)n}\right\} < \left\{\gamma_{0i,0(n-1)}\right\} < \left\{\gamma_{0j,(n-1)n}\right\} < \gamma_{0(n-1),(n-1)n},
\end{equation*}
ordered for increasing $1\leq i \leq \frac{n-3}{2}$ and $\frac{n-1}{2}\leq j \leq n-2$. Before defining Step $B$ of the algorithm, we perform appropriate Hurwitz moves between paths corresponding to disjoint vanishing cycles so that the order is:
\begin{equation*}
    \gamma_{0i,(n-1)n} < \gamma_{0i,0(n-1)} < \gamma_{0(i+1),(n-1)n}
\end{equation*}
for all $i$. The mutations constituting Step $B$ of the algorithm consist in the simultaneous Hurwitz moves, for all $i$:
\begin{equation*}
   \left(\gamma_{0i,(n-1)n} ,\gamma_{0i,0(n-1)}\right) \mapsto  \left(\gamma_{0i,0(n-1)},\tau^{-1}_{\gamma_{0i,0(n-1)}}\gamma_{0i,(n-1)n}\right).
\end{equation*}
Define $\gamma_{0i,0n}:= \tau^{-1}_{\gamma_{0i,0(n-1)}}\gamma_{0i,(n-1)n}$ the resulting path for each $i$. Finally, we perform appropriate Hurwitz moves on disjoint cycles, so that the final order of the paths is:
\begin{equation*}
    \gamma_{0i,0(n-1)} < \gamma_{0i,0n} < \gamma_{0(i+1),0(n-1)}
\end{equation*}
for all $i$. Define $\gamma_{I,n-1}:=\gamma_{0I,0(n-1)}$ and $\gamma_{I,n}:=\gamma_{0I,0n}$ the final vanishing paths. This concludes the even case: we have inductively constructed a series of mutations on $\G$, explicitly describing the vanishing paths associated to the cycles we will call the Iyama ones.

We deal with the odd case in an analogous way; denote by $\G$ the distinguished collections of vanishing paths $\gamma_{ij,\ell m}$ prescribed in Section \ref{section:ACmorsification}, and by $\Gg$ the sub-collection:
\begin{equation*}
    \Gg:= \left\{\gamma_{ij,\ell m} \mid \ell m \neq 0n \right\}.
\end{equation*}
With the exception of those corresponding to paths belonging to $\G \setminus \Gg$, all the positive cycles are disjoint from the saddles corresponding to $\gamma_{ij,0n}$ (Figure \ref{disjointcycles}, right): perform appropriate Hurwitz moves, so that the objects of $\G \setminus \Gg$ are the last ordered  paths of $\G$. Using an inductive argument, assume that we have constructed a series of mutations on $\Gg$ (acting by Hurwitz moves) such that the final paths are:
\begin{equation*}
    \left\{\gamma_{I,J}:=\gamma_{0I,0J} \mid I,J \neq n\right\}.
\end{equation*}

In analogy to the even case, we choose a further perturbation of the Lefschetz fibration that separates all critical values, in such a way that the order on the paths $\gamma_{ij,0n}$ is for increasing $j$ for saddles, and for increasing $i$ for positive ones. Retracing the algorithm in Proposition \ref{combalgorithm}, the mutation constituting Step $h$ consists in the Hurwitz move:
\begin{align*}
    \left(\gamma_{0\frac{2n-h-1}{2},0n},\gamma_{\frac{h+1}{2}\frac{2n-h-1}{1},0n}\right) &\mapsto   \left(\tau_{\gamma_{0\frac{2n-h-1}{2},0n}}\gamma_{\frac{h+1}{2}\frac{2n-h-1}{1},0n},\gamma_{0\frac{2n-h-1}{2},0n}\right) \\
    \left(\gamma_{\frac{h}{2}\frac{2n-h-2}{2},0n},\gamma_{0\frac{h}{2},0n}\right) &\mapsto   \left(\tau^{-1}_{\gamma_{0\frac{h}{2},0n}}\gamma_{\frac{h}{2}\frac{2n-h-2}{2},0n},\gamma_{0\frac{h}{2},0n}\right)
\end{align*}
for $h$ odd and even respectively. In analogy to the even case, this mutation is followed by appropriate Hurwitz moves between disjoint vanishing cycles. At each step, for $h$ odd and even respectively, define:
\begin{align*}
    \gamma_{0\frac{h+1}{2},0n}&:=\tau_{\gamma_{0\frac{2n-h-1}{2}, 0n}}\gamma_{\frac{h+1}{2}\frac{2n-h-1}{1}, 0n} \\
    \gamma_{0\frac{2n-h-2}{2}, 0n}&:=\tau^{-1}_{\gamma_{0\frac{h}{2}, 0n}}\gamma_{\frac{h}{2}\frac{2n-h-2}{2}, 0n}.
\end{align*}
The constructed paths $\gamma_{0i,0n}$ are the desired $\gamma_{I,n}$, as given by  $\gamma_{I,n}:=\gamma_{0I,0n}$. Moreover, the final order of $\{\gamma_{0i,0n}\}$ is the following:
\begin{equation*}
    \dots < \gamma_{01,0n} < \gamma_{02,0n} < \dots < \gamma_{\frac{n-2}{2},0n} < \gamma_{\frac{n}{2},0n} < \dots < \gamma_{0(n-1),0n}.
\end{equation*}
\end{proof}

In analogy with Remark \ref{rk:notation}, denote by $V_{pq,rs}$ the vanishing cycle associated to path $\gamma_{pq,rs}$ for any $\gamma_{pq,rs}$ defined in the proof of Proposition \ref{orderIyamagenerators}. Further define $\{\Iyama_{I,J}\}$ to be the vanishing cycles associated to the constructed distinguished collection $\{\gamma_{I,J}\}$.

\begin{figure}
    \centering
    \begin{minipage}{.55 \textwidth}
        \centering
        \begin{tikzpicture}
            [remember picture,scale=.7,stop/.style={shape=circle, fill=black, inner sep=.7pt}]
            \tikzstyle{reverseclip}=[insert path={(current page.north east) --
(current page.south east) -- (current page.south west) -- (current page.north west) -- (current page.north east)}
            ]
    
    \draw (-3,-12.2) to (-3,-12.5) to(1.75,-12.5) to (1.75,-12.2);
    \node at (-.55,-12.8) {$\Xnn$};

    \draw (6,-12.2) to (6,-12.5) to(1.85,-12.5) to (1.85,-12.2);
    \node at (3.925,-12.8) {$\Ynn$};

    \draw (0,1) arc (90:45:.5 and 1) coordinate (A1);
    \draw (0,-1) arc (-90:-45:.5 and 1) coordinate (A2);
    \draw (0,1) arc (-270:-185:.5 and 1) coordinate(B1) ;
    \draw (0,-1) arc (-90:-108:.5 and 1) coordinate(B2);
    \draw[dotted] (B1) to[out=-90, in=142, looseness=.85] (B2);
     \draw (0,1) to (2.7,1);
     \draw (0,-1) to (2.7,-1);
    \draw[dotted] (2.7,1) arc (-270:-90:.5 and 1);
    \draw[dotted] (2.7,1) to[out=0,in=100] (3,.6);
    \draw (3,.4) to (3,.1);
    \draw (3,-.1) to (3,-.4);
    \draw (2.7,-1) to[out=0,in=-100] (3,-.6);

    \draw (0,-4) arc (90:45:.5 and 1) coordinate (A3);
    \draw (0,-6) arc (-90:-45:.5 and 1) coordinate (A4);
    \draw (0,-4) arc (-270:-255:.5 and 1) coordinate(B3) ;
    \draw (0,-6) arc (-90:-105:.5 and 1) coordinate(B4) ;
    \draw[dotted] (0,-4) arc (-270:-90:.5 and 1) ;
    \draw (0,-4) to (2.7,-4);
    \draw (0,-6) to (2.7,-6);
    \draw[dotted] (2.7,-4) arc (-270:-90:.5 and 1);

    \draw (2.7,-4) to[out=0,in=100] (3,-4.4);
    \draw (3,-4.6) to (3,-4.9);
    \draw (3,-5.1) to (3,-5.4);
    \draw[dotted] (3,-5.6) to[out=-100,in=0] (2.7,-6);
            
    \draw (0,-7) arc (90:45:.5 and 1) coordinate (A5);
    \draw (0,-9) arc (-90:-45:.5 and 1) coordinate (A6);
    \draw (0,-7) arc (-270:-255:.5 and 1) coordinate(B5) ;
    \draw (0,-9) arc (-90:-105:.5 and 1) coordinate(B6) ;
    \draw[dotted] (0,-7) arc (-270:-90:.5 and 1) ;
    \draw (0,-7) to (2.7,-7);
    \draw (0,-9) to (2.7,-9);
    \draw[dotted] (2.7,-7) arc (-270:-90:.5 and 1);

    \draw (2.7,-7) to[out=0,in=100] (3,-7.4);
    \draw (3,-7.6) to (3,-7.9);
    \draw[dotted]  (3,-8.1) to (3,-8.4);
    \draw (3,-8.6) to[out=-100,in=0] (2.7,-9);

    \draw (0,-10) arc (90:45:.5 and 1) coordinate (A7);
    \draw (0,-12) arc (-90:-45:.5 and 1) coordinate (A8);
    \draw (0,-10) arc (-270:-255:.5 and 1) coordinate(B7) ;
    \draw (0,-12) arc (-90:-175:.5 and 1) coordinate(B8);
    \draw[dotted] (B7) to[in=90, out=-142, looseness=.85] (B8);
    \draw (0,-10) to (2.7,-10);
    \draw (0,-12) to (2.7,-12);
    \draw[dotted] (2.7,-10) arc (-270:-90:.5 and 1);

    \draw (2.7,-10) to[out=0,in=90,looseness=.8] (3,-10.9);
    \draw[dotted]  (3,-11.1) to (3,-11.3);
    \draw (3,-11.5) to (3,-11.6);
    \draw (3,-11.8) to[out=-100,in=0] (2.7,-12);
            
    \draw[loosely dotted] (1,-2) to (1,-3);
    \draw[loosely dotted] (2,-2) to (2,-3);
    
    \draw (A1) to[bend right=60]coordinate[pos=0.85] (glueing1) (A8);
    
    \draw (A2) to[out=200, in=160, looseness=1.5] (-.3,-2);
    \draw (-.3,-3) to[out=200, in=160, looseness=1] (A3);
    \draw (A4) to[out=200, in=160, looseness=2] (A5);
    \draw (A6) to[out=200, in=160, looseness=2] coordinate[pos=0.55] (glueing2) (A7);
            
    \draw[dashed] (glueing1) to (glueing2);
            
    \begin{scope}[decoration={markings,mark=at position 0.8 with {\arrow{<}}}] 
    \draw[color=gray,postaction={decorate}] (1.8,-1) arc (-90:90:.5 and 1) ;
    \draw[color=gray,dotted] (1.8,1) arc (-270:-90:.5 and 1) ;

    \draw[color=gray,postaction={decorate}] (1.8,-6) arc (-90:90:.5 and 1) ;
    \draw[color=gray,dotted] (1.8,-4) arc (-270:-90:.5 and 1) ;

    \draw[color=gray,postaction={decorate}] (1.8,-9) arc (-90:90:.5 and 1) ;
    \draw[color=gray,dotted] (1.8,-7) arc (-270:-90:.5 and 1) ;

    \draw[color=gray,postaction={decorate}] (1.8,-12) arc (-90:90:.5 and 1) ;
    \draw[color=gray,dotted] (1.8,-10) arc (-270:-90:.5 and 1) ;
    \end{scope}
            
        
    \begin{scope}[decoration={markings,mark=at position 0.1 with {\arrow{>}}}] 
    \draw[color=gray,postaction={decorate}] (3,-.5) to[out=170,in=50] (-1,-1) to[out=230,in=90] (-2,-5) to[out=270,in=120] (-1,-10) to[out=300,in=190] (3,-11);
    \draw[color=gray,postaction={decorate}] (3,0) to[out=180,in=60] (-.8,-.8) to[out=240,in=100] (-1.5,-5) to[out=280,in=140] (0,-8) to[out=320,in=200](3,-8);
    \draw[color=gray,postaction={decorate}] (3,.5) to[out=190,in=60] (-.3,-.3) to[out=240,in=-240] (-.5,-4.5) to[out=-60,in=180] (3,-5.5);
    \end{scope}
            
    \begin{scope}[decoration={markings,mark=at position 0.2 with {\arrow{>}}}] 
    \draw[color=gray,postaction={decorate}] (3,-5) to[out=180,in=40] (-.5,-5.5) to[out=220,in=120] (-.8,-10) to[out=300,in=180] (3,-11.4);
    \draw[color=gray,postaction={decorate}] (3,-4.5) to[out=210,in=60] (-.5,-6) to[out=240,in=120] (-.5,-7) to[out=300,in=170] (3,-8.5);
    \end{scope}
            
    \begin{scope}[decoration={markings,mark=at position 0.2 with {\arrow{>}}}] 
    \draw[color=gray,postaction={decorate}] (3,-7.5) to[out=180,in=50] (-.5,-8) to[out=230,in=130] (-.5,-10) to[out=310,in=150] (.5,-10.8) to[out=330,in=170] (3,-11.7);
    \end{scope}
            
    \draw[color=black] (3,-.4) to[out=350, in=370](3,-11.1);
    \draw[color=gray] (3,-.5) to[out=350, in=370](3,-11);
    \draw[color=black] (3,-.6) to[out=350, in=370](3,-10.9);

    \begin{pgfinterruptboundingbox}
    \path[clip]  (4.55,-1.5) -- (5.3,-2.5)  -- (5.2,-2) -- (4.53,-1.15) -- cycle[reverseclip];
    \end{pgfinterruptboundingbox}

    \begin{pgfinterruptboundingbox}
    \path[clip]   (4.7,-9.9) -- (4.4,-10.17) --  (4.18,-10.6)  --  (4.55,-10.3)  -- cycle[reverseclip];
    \end{pgfinterruptboundingbox}
    
    \begin{pgfinterruptboundingbox}
    \path[clip]    (4.27,-10.3) -- (4.05,-10.45) --  (4,-10.75)  --  (4.22,-10.57) -- cycle[reverseclip];
    \end{pgfinterruptboundingbox}
        
    \draw[color=black] (3,.1) to[out=0, in=380](3,-8.1);
    \draw[color=gray] (3,0) to[out=0, in=380](3,-8);
    \draw[color=black] (3,-.1) to[out=0, in=380](3,-7.9);

    \begin{pgfinterruptboundingbox}
    \path[clip]  (4.27,-.38) -- (4.45,-.9)  -- (4.7,-1.2) -- (4.6,-.7)-- cycle[reverseclip];
    \end{pgfinterruptboundingbox}
        
    \draw[color=black] (3,.4) to[out=10, in=0](3,-5.4);
    \draw[color=gray] (3,.5) to[out=10, in=0](3,-5.5);
    \draw[color=black] (3,.6) to[out=10, in=0] (3,-5.6);

    \begin{pgfinterruptboundingbox}
    \path[clip] (3.92,-4.92) -- (3.79,-5.05) --  (3.87,-5.27)  --  (4.01,-5.15) -- cycle[reverseclip];
    \end{pgfinterruptboundingbox}
 
    \begin{pgfinterruptboundingbox}
    \path[clip]  (4.15,-7.08) -- (3.95,-7.2) --  (3.85,-7.6)  --  (4.1,-7.4) -- cycle[reverseclip];
    \end{pgfinterruptboundingbox}

    \begin{pgfinterruptboundingbox}
    \path[clip]  (3.55,-7.55) -- (3.35,-7.7) --  (3.53,-7.84)  --  (3.71,-7.72) -- cycle[reverseclip];
    \end{pgfinterruptboundingbox}
    
    \draw[color=black] (3,-4.4)to[out=30, in=350](3,-8.6);
    \draw[color=gray] (3,-4.5) to[out=30, in=350](3,-8.5);
    \draw[color=black] (3,-4.6) to[out=30, in=350] (3,-8.4);
          
    \begin{pgfinterruptboundingbox}
     \path[clip]  (3.78,-7.77) -- (3.67,-8) --  (3.79,-8.18)  --  (3.9,-7.95) -- cycle[reverseclip];
    \end{pgfinterruptboundingbox}

    \begin{pgfinterruptboundingbox}
    \path[clip]  (3.75,-5.08) -- (3.57,-5.25) --  (3.75,-5.4)   --  (3.87,-5.2) -- (3.95,-5.6)-- (4.1,-5.7) -- (4.05,-5.33) -- cycle[reverseclip];
    \end{pgfinterruptboundingbox}

    \begin{pgfinterruptboundingbox}
    \path[clip]   (4.62,-6.22) -- (4.5,-6.45) --  (4.57,-6.65)  --  (4.7,-6.5)  -- cycle[reverseclip];
    \end{pgfinterruptboundingbox}
       
    \draw[color=black] (3,-4.9) to[out=0, in=0](3,-11.5);
    \draw[color=gray] (3,-5) to[out=0, in=0](3,-11.4);
    \draw[color=black] (3,-5.1) to[out=0, in=0] (3,-11.3);

    \begin{pgfinterruptboundingbox}
    \path[clip]   (4.2,-10.62) -- (3.9,-10.9) --  (3.7,-11.3)  --  (4.05,-11.1)   -- cycle[reverseclip];
    \end{pgfinterruptboundingbox}

    \draw[color=black]  (3,-7.4) to[out=0, in=350](3,-11.8);
    \draw[color=gray]  (3,-7.5)  to[out=0, in=350](3,-11.7);
    \draw[color=black] (3,-7.6)  to[out=0, in=350] (3,-11.6);
            
    \node at (2,1.3) {$\LL_{1,2}$};
    \node at (5,1) {$\LL_{1,n-2}$};
    \path (4.3,.7) edge[->] (4,.4);
    \node at (7,-1.3) {$\LL_{1,n-1}$};
    \path (6,-1.3) edge[->](5,-1.3);
    \node at (3.5,-1.5) {$\LL_{1,n}$};
    \path (3.5,-1.2) edge[->](3.5,-.8);

    \node at (2,-3.7) {$\LL_{n-3,n-2}$};
    \node at (2,-6.7) {$\LL_{n-2,n-1}$};
    \node at (2,-9.7) {$\LL_{n-1,n}$};

    \node at (6.5,-10.9) {$\LL_{n-3,n}$};
    \path (5.2,-10.9) edge[->](4.5,-10.9);
    \node at (6,-11.7) {$\LL_{n-2,n}$};
    \path (4.7,-11.7) edge[->](4,-11.7);

    \end{tikzpicture}
            \vskip .5cm
        (a)
        \end{minipage}
        \qquad
        \begin{minipage}{.35\textwidth}
        \centering
        \begin{tikzpicture}
        \draw (0,0) to (3,0);
        \draw (0,.5) to (3,-.5);
        \draw (0,-.5) to (3,0.5);
        \draw (0,1) to (3,-1);
        \draw (0,-1) to (3,1);
        
        \node at (3.5,1) {$\LL_{1,J}$};
        \node at (3.5,.5) {$\LL_{2,J}$};
        \node at (3.5,0) {$\LL_{3,J}$};
        \node at (3.5,-.5) {$\LL_{4,J}$};
        \node at (3.5,-1) {$\LL_{5,J}$};
        \end{tikzpicture}
        \vskip .5cm
        (b) 
        \vskip 1.5 cm
        \begin{tikzpicture}
        [align=center, v/.style={draw,shape=circle, fill=black, minimum size=1.2mm, inner sep=0pt, outer sep=0pt}, font=\small, label distance=1pt,
        ]
        \node[v, ] (1) at (60:1.5cm) {};
        \node[v,] (2) at (0:1.5cm) {};
        \node[v, ] (3) at (-60:1.5cm) {};
        \node[v,] (4) at (-120:1.5cm) {};
        \node[v,] (5) at (180:1.5cm) {};
        \node[v, ] (6) at (120:1.5cm) {};
        \path (1) edge (2);
        \path (1) edge (3);
        \path (2) edge (3);
        \path (1) edge (4);
        \path (2) edge (4);
        \path (3) edge (4);
        \path (1) edge (5);
        \path (2) edge (5);
        \path (3) edge (5);
        \path (4) edge (5);
        \path (1) edge (6);
        \path (2) edge (6);
        \path (3) edge (6);
        \path (4) edge (6);
        \path (5) edge (6);
        
        \node at (60:2cm) {$\LL_{1,2}$};
        \node at (0:2cm) {$\LL_{2,3}$};
        \node at (-60:2cm) {$\LL_{3,4}$};
        \node at (240:2cm) {$\LL_{4,5}$};
        \node at (180:2cm) {$\LL_{5,6}$};
        \node at (120:2cm) {$\LL_{6,7}$};
        
        \end{tikzpicture}
        \vskip .5cm
        (c)
        \end{minipage}
    \caption{(a) The Milnor fibre $\Sigma_{n+1}$ of $f_{n+1}$ with the Iyama vanishing cycles. It is made up of $n-1$ cylinders, whose zero-sections $\Iyama_{I,I+1}$ bound two connected components $\Ynn$ and $\Xnn$. $\Xnn$ is a 1-punctured surface of genus 0, $\Ynn$ is the thickening of the complete graph associated to the set of $n-1$ vertices. (b) Schematic representation of the restriction of the vanishing cycles to a cylinder (unperturbed, so that multiple cycles intersect in one point). (c) (Unperturbed) Schematic representation of the complete graph pattern (all diagonals in an $n$-gon) of the cycles restricted to $\Xn$ ($n=8$), where $\LAMBDA$ is depicted a collection of points.}
    \label{fig:Iyamasurface}
    \end{figure}
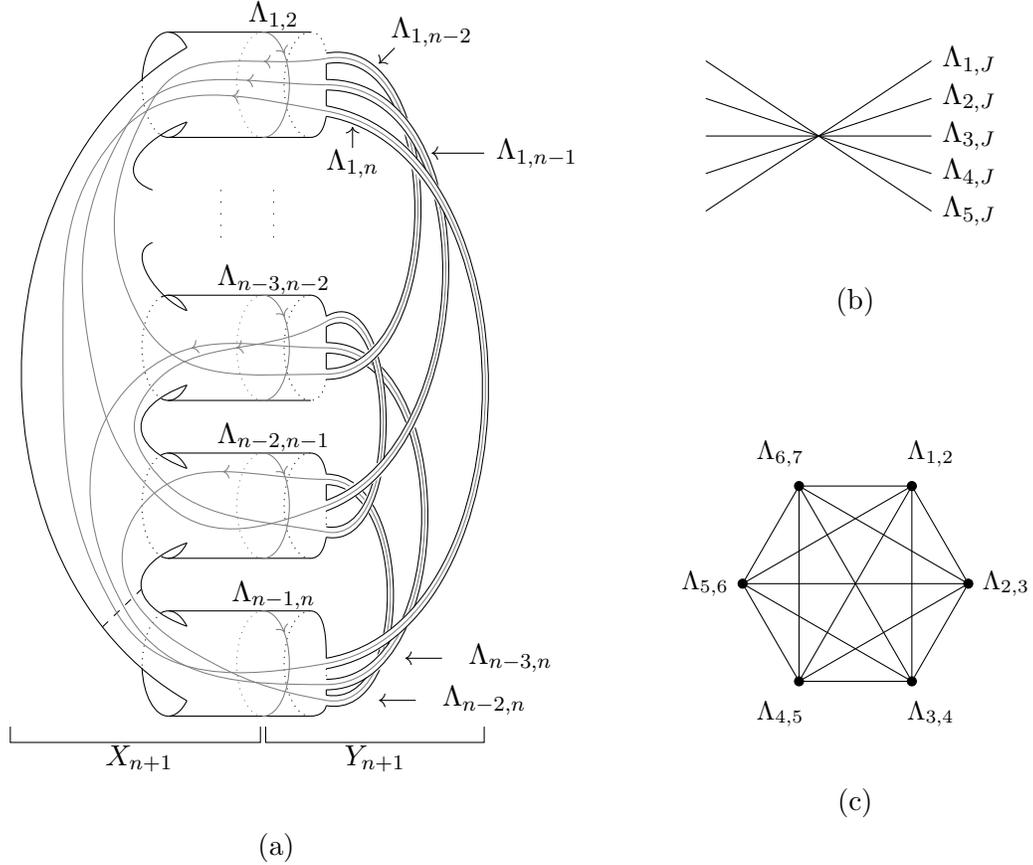

\begin{proposition}\label{prop:homeofibre}
    The Milnor fibre $\SIy_n$ of $f_n$ equipped with the Iyama vanishing cycles $\Iyama_{I,J}$ ($1\leq I < J \leq n-1$) is homeomorphic to the surface in Figure \ref{fig:Iyamasurface}.
\end{proposition}

\begin{proof}
    The case $n=3$ is trivial: $\SIy_3$ is topologically a cylinder, with a single vanishing cycle going around its waist, which is exactly what is prescribed in Figure \ref{fig:Iyamasurface} for such $n$. We prove the claim by constructing $\SIy_{n+1}$ from $\SIy_{n}$ via a glueing process. Denote by $\LAMBDA'$ and $\LAMBDA$ respectively the collection of cycles $\LAMBDA':=\left\{\Iyama_{I,I+1},I=1,\dots,n-2 \right\}$ and $\LAMBDA:=\left\{\Iyama_{I,I+1},I=1,\dots,n-1 \right\}$. We assume that $\SIy_{n}$ results from the glueing
    \begin{equation*}
        \SIy_{n}=\Xn \#_{\LAMBDA'} \Yn
    \end{equation*} 
    of two connected surfaces $\Xn$ and $\Yn$ along $\LAMBDA'$. We also assume that $\Xn$ is homeomorphic to a $1$-punctured surface of genus 0 and $n-2$ boundary components (the collection $\LAMBDA'$) (see Figure \ref{fig:Iyamasurface}). We first prove that there is a Riemann surface $\Xnn$ of genus 0 constructed from $\Xn$ and embedding naturally into $\SIy_{n+1} \setminus \LAMBDA$.

    Suppose $n+1$ is even. The (by assumption) connected component $\Xn \subset \Sigma_{n+1}$ is depicted in Figure \ref{fig:inductivehomeomorphismeven} (a) (shaded in grey). From this one can see that $\Xnn$ is obtained by glueing (dashed line in Figures \ref{fig:Iyamasurface} and \ref{fig:inductivehomeomorphismeven}) $\Xn$ to a surface of genus 0 bounded by $\Iyama_{n-1,n}$ (shaded in yellow in Figure \ref{fig:inductivehomeomorphismeven}).

\begin{figure}
    \centering
    \begin{minipage}{1\textwidth}
        \centering
        \begin{tikzpicture}[scale=.63, 
            cylinderop/.pic={
            \draw (-.45,1.125) to (-.45,-1.125); 
            \draw (.45,1.125) to (.45,-1.125); 
            \draw (-1.125,.45) to (-.63,.45); 
            \draw (-.27,.45) to (.27,.45); 
            \draw (1.125,.45) to (.63,.45); 
            \draw (-1.125,-.45) to (-.63,-.45); 
            \draw (-.27,-.45) to (.27,-.45); 
            \draw (1.125,-.45) to (.63,-.45);
            }, cylinder/.pic={
            \draw (-.45,1.125) to (-.45,-1.125); 
            \draw (.45,1.125) to (.45,-1.125); 
            \draw (-.45,1.125) to (.45,1.125);
            \draw (-1.125,.45) to (-.63,.45); 
            \draw (-.27,.45) to (.27,.45); 
            \draw (1.125,.45) to (.63,.45); 
            \draw (-1.125,-.45) to (-.63,-.45); 
            \draw (-.27,-.45) to (.27,-.45); 
            \draw (1.125,-.45) to (.63,-.45);
            }
            ]
        
        \draw (-5,4) pic {cylinder} (0,4) pic {cylinder} (5,4) pic {cylinder} (10,4) pic {cylinder} (0,-2) pic {cylinderop} (5,-2) pic {cylinderop} (10,-2) pic {cylinderop};
        
        \draw(-6.77,4.7) to (-6.77,3.3);

        \draw[name path=line11] (-3.22,4.715) to (-1.78,3.285);
        \draw[name path=line12] (-3.22,3.285) to (-1.78,4.715);
        \path [name intersections={of=line11 and line12,by=int1}];

        \draw[name path=line21] (1.78,4.715) to (3.22,3.285);
        \draw[name path=line22] (1.78,3.285) to (3.22,4.715);
        \path [name intersections={of=line21 and line22,by=int2}];

        \draw[name path=line31] (6.78,4.715) to (8.22,3.285);
        \draw[name path=line32] (6.78,3.285) to (8.22,4.715);
        \path [name intersections={of=line31 and line32,by=int3}];
        
        \draw[name path=line101] (1.78,-1.285) to (3.22,-2.715);
        \draw[name path=line102] (1.78,-2.715) to (3.22,-1.285);
        \path [name intersections={of=line101 and line102,by=int10}];

        \draw[name path=line111] (6.78,-1.285) to (8.22,-2.715);
        \draw[name path=line112] (6.78,-2.715) to (8.22,-1.285);
        \path [name intersections={of=line111 and line112,by=int11}];
        
        \draw[name path=line61] (-0.715,2.215) to (.715,-.215);
        \draw[name path=line62] (-0.715,-.215) to (.715,2.215);
        \path [name intersections={of=line61 and line62,by=int6}];

        \draw[name path=line71](4.285,2.215) to (5.715,-.215);
        \draw[name path=line72] (4.285,-.215) to (5.715,2.215);
        \path [name intersections={of=line71 and line72,by=int7}];

        \draw[name path=line81] (9.285,2.215) to (10.715,-.215);
        \draw[name path=line82] (9.285,-.215) to (10.715,2.215);
        \path [name intersections={of=line81 and line82,by=int8}];

        \draw[name path=line51] (-5.715,2.215) to (-1.78,-1.285);
        \draw[name path=line52] (-4.285,2.215) to (-1.78,-2.715);
        \path [name intersections={of=line51 and line52,by=int5}];

            \draw[dashed] (-1.5,4.715) to (-1.5,3.285);

        \draw[color=red](-5.85,3.3) to[bend left=45] (-5.85,4.7);
        \draw[color=red] (-5.6,4.8) to[bend left=45] (-4.4,4.8);
        \draw[color=red](-4.2,4.7) to[bend left=45] (-4.2,3.3);
        \draw[color=red] (-4.4,3.2) to[bend left=45] (-5.6,3.2);
               
        \draw[color=blue] (-4.05,3.3) to[out=10,in=200] (int1) to[out=380,in=140] (-.9,3.3);
        \draw[color=blue] (-.7,3.2) to[out=280,in=110] (int6) to[out=290,in=120] (.65,-1.25);
        \draw[color=blue] (.8,-1.35) to[bend left=30] (.8,-2.7);
        \draw[color=blue]  (.6,-2.8) to[bend left=30] (-.6,-2.8);
        \draw[color=blue] (-.8,-2.7) to[out=160, in=-55, looseness=1.2] (int5) to[out=130, in=250, looseness=.8] (-4.3,3.15);
        
        \draw[color=magenta] (-4.1,3.35) to[out=50,in=160] (int1) to[out=340,in=170] (-1,3.3);
        \draw[color=magenta] (-.57,3.2) to[out=320,in=90,looseness=.6] (int6) to[out=270,in=60] (-.55,-1.25);
        \draw[color=magenta] (-.8,-1.4) to[out=190, in=-50, looseness=1.1] (int5) to[out=130, in=230, looseness=1.2] (-4.4,3.15);
        
        \draw[color=orange](-.85,3.3) to[bend left=45] (-.85,4.7);
        \draw[color=orange] (-.6,4.8) to[bend left=45] (.6,4.8);
        \draw[color=orange](.8,4.7) to[bend left=45] (.8,3.3);
        \draw[color=orange] (.6,3.2) to[bend left=45] (-.5,3.2);

        \draw[color=ForestGreen](4.2,3.3) to[bend left=45] (4.2,4.7);
        \draw[color=ForestGreen](4.5,4.8) to[bend left=45] (5.6,4.8);
        \draw[color=ForestGreen](5.8,4.7) to[bend right=5] (int3);
        \draw[color=ForestGreen] (int3) to[bend left=5] (9.2,3.3);
        \draw[color=ForestGreen] (9.4,3.2) to[bend left=5] (int8);
        \draw[color=ForestGreen] (int8) to[bend left=30] (9.4,-1.2);
        \draw[color=ForestGreen] (9.2,-1.3) to[bend left=5] (int11);
        \draw[color=ForestGreen] (int11) to[bend left=30] (5.8,-1.3);
        \draw[color=ForestGreen] (5.6,-1.2) to[bend left=5] (int7);
        \draw[color=ForestGreen] (int7) to[bend right=5] (4.4,3.2);

        \draw (10.8,-2.7) to [bend left=5] (11.75,-2);
        \draw (9.4,-2.8)to [bend right=45] (10.6,-2.8);
        \draw  (5.9,-1.3) to[bend left=5] (int11) to[bend right=5](9.2,-2.7);
        \draw  (5.5,-1.2) to[bend left=30] (int7) to[bend left=10] (5.6,3.2);
        \draw (5.8,3.3) to[bend left=5] (int3) to[bend right=5] (9.2,4.7);
        \draw (9.4,4.8) to [bend left=45] (10.6,4.8);
        \draw (10.8,4.7) to [bend right=5] (11.75,4.2);

        \draw (9.4,-2.9) to [bend left=5] (10,-3.8);
        \draw (9.2,-2.6) to [bend left=45] (9.2,-1.4);
        \draw  (9.4,-1.1) to[bend right=5] (int8) to[bend left=5]  (10.6,3.2);
        \draw  (10.8,3.3) to [bend left=5] (11.75,3.8);

        \draw (10.2,-3.8) to [bend left=5] (10.6,-2.9);
        \draw (10.8,-2.5) to [bend left=30] (11.75,-1.8);

        \draw[color=cyan] (-.65,-2.8) to (0,-3.8);
        \draw[color=cyan] (-.8,-2.6) to[out=150, in=-53, looseness=1.2] (int5) to[out=133, in=245, looseness=.8] (-4.35,3.15);
        \draw[color=cyan] (-4.05,3.35) to[out=15,in=190] (int1) to[out=370,in=150] (-1.05,3.35);
        \draw[color=cyan] (-.65,3.2) to[out=280,in=100] (int6) to[out=280,in=140] (.6,-1.3);
        \draw[color=cyan] (.9,-1.4) to[out=-40,in=180] (int10) to[out=0,in=150] (4.2,-2.7);
        \draw[color=cyan] (4.4,-2.75) to[out=-30,in=110] (5.2,-3.8);
        
        \draw[color=green] (-.67,-2.78) to[out=300] (-.1,-3.8);
        \draw[color=green] (-.75,-2.55) to[out=140, in=-52, looseness=1.2] (int5) to[out=135, in=245, looseness=.8] (-4.45,3);
        \draw[color=green] (-4,3.4) to[out=30,in=180] (int1) to[out=0,in=165] (-1.15,3.37);
        \draw[color=green] (-.6,3.15) to[out=290,in=95] (int6) to[out=275,in=140] (.5,-1.25);
        \draw[color=green] (.9,-1.3) to[out=-30,in=170] (int10) to[out=-10,in=150] (4.1,-2.7);
        \draw[color=green] (4.5,-2.75) to[out=-30,in=110] (5.3,-3.8);
        
        \draw[fill=yellow,draw=none, fill opacity=.3] (-6.8,3.3) to (-5.85,3.3) to[bend left=45] (-5.85,4.7) to (-6.8,4.7);
        \draw[fill=yellow,draw=none, fill opacity=.3] (-4.4,4.7) to (-5.6,4.7) to (-5.6,4.8) to[bend left=45] (-4.4,4.8);
        \draw[fill=yellow,draw=none, fill opacity=.3] (-4.4,3.3) to (-5.6,3.3) to (-5.6,3.2) to[bend right=45] (-4.4,3.2);
        \draw[fill=yellow,draw=none, fill opacity=.3] (-4.2,4.7) to[bend left=45] (-4.2,3.3) to (-3.22,3.3) to (int1) to (-3.22,4.715);
        \draw[fill=yellow,draw=none, fill opacity=.3] (-1.78,4.715) to (int1)  to (-1.78,3.285) to (-1.5,3.285) to (-1.5,4.715) ;

        \draw[fill=gray,draw=none, fill opacity=.3] (-1.5,3.285) to (-1.5,4.715)to (-.85,4.7) to[bend right=45] (-.85,3.3);
        \draw[fill=gray,draw=none, fill opacity=.3] (-.6,4.7) to  (-.6,4.8) to[bend left=45] (.6,4.8) to (.6,4.7);
        \draw[fill=gray,draw=none, fill opacity=.3] (-.5,3.3) to (-.5,3.2) to[bend right=45] (.6,3.2) to  (.55,3.3);
        \draw[fill=gray,draw=none, fill opacity=.3](.8,4.715) to[bend left=45] (.8,3.285) to (1.6,3.285) to (1.6,4.715);
        \draw[fill=gray,draw=none, fill opacity=.3] (1.6,4.715) to (1.6,3.285) to (1.78,3.285) to (int2) to(1.78,4.715);
        \draw[fill=gray,draw=none, fill opacity=.3] (int2) to (3.22,3.3) to (4.2,3.3) to[bend left=45] (4.2,4.7) to (3.22,4.715);
        \draw[fill=gray,draw=none, fill opacity=.3](4.5,4.8) to[bend left=45] (5.6,4.8) to (5.6,4.7) to (4.5,4.7);
        \draw[fill=gray,draw=none, fill opacity=.3] (int7) to[bend right=5] (4.4,3.2) to  (4.4,3.3) to (5.6,3.3) to (5.6,3.2) to[bend right=10] (int7);
        \draw[fill=gray,draw=none, fill opacity=.3](5.8,4.7)  to[bend right=5] (int3) to (6.78,4.715);
        \draw[fill=gray,draw=none, fill opacity=.3] (5.8,3.3) to[bend left=5] (int3) to  (6.78,3.285);
        \draw[fill=gray,draw=none, fill opacity=.3]  (int3) to[bend right=5] (9.2,4.7) to  (8.22,4.7);
        \draw[fill=gray,draw=none, fill opacity=.3] (int3) to[bend left=5] (9.2,3.3)to  (8.22,3.285);
        \draw[fill=gray,draw=none, fill opacity=.3] (9.4,4.8) to [bend left=45] (10.6,4.8) to (10.6,4.7) to (9.4,4.7);
        \draw[fill=gray,draw=none, fill opacity=.3] (10.8,4.7) to [bend right=5] (11.75,4.2) to (11.75,4.7);
        \draw[fill=gray,draw=none, fill opacity=.3]  (5.5,-1.2) to[bend left=30] (int7) to[bend right=5] (5.6,-1.2);
        \draw[fill=gray,draw=none, fill opacity=.3]  (5.8,-1.3) to[bend right=30] (int11) to[bend right=5] (5.9,-1.3);
        \draw[fill=gray,draw=none, fill opacity=.3] (9.2,-1.4) to[bend right=45] (9.2,-2.6) to  (9.2,-2.7) to[bend left=5] (int11) to[bend right=5] (9.2,-1.3);

        \draw[fill=gray,draw=none, fill opacity=.3] (10.8,-2.5) to [bend left=30] (11.75,-1.8) to (11.75,-2) to[bend right=5] (10.8,-2.7);

        \draw[fill=gray,draw=none, fill opacity=.3]  (10,-3.8)to [bend right=5](9.4,-2.9) to (9.4,-2.8)to [bend right=45] (10.6,-2.8) to (10.6,-2.9) to[bend right=5] (10.2,-3.8);

        \draw[fill=gray,draw=none, fill opacity=.3] (9.4,3.3) to (9.4,3.2) to[bend left=5] (int8) to[bend left=5]  (10.6,3.2) to (10.6,3.3);

        \draw[fill=gray,draw=none, fill opacity=.3] (9.4,-1.1)  to[bend right=5] (int8) to[bend left=30]  (9.4,-1.2);
        \draw[fill=gray,draw=none, fill opacity=.3]  (10.8,3.3) to [bend left=5] (11.75,3.8) to  (11.75,3.3);

        \node at (12.7,1) {$\dots$};
        \node at (7.5,-4) {$\vdots$};

        \node[color=red] at (-6.8,5.3) {$\LL_{n-1,n}$};
        \node[color=orange] at (-2,5.3) {$\LL_{n-2,n-1}$};
        \node[color=ForestGreen] at (2.9,5.3) {$\LL_{n-3,n-2}$};
        \node at (8.5,5.3) {$\LL_{1,2}$};

        \node[color=magenta] at (-6.5,0) {$V_{0(n-2),(n-1)n}$};
        \path (-5.25,.5) edge[->, color=magenta] (-4.8,1.75);

        \node at (-5.5,-2.5) {$\left\{V_{0i,(n-1)n}\right\}$};
        \path (-4.5,-2) edge[->] (-2.75,-1.5);
        
        \end{tikzpicture} \vskip .0cm
            (a)
        \vskip .2cm
    \centering
        \begin{tikzpicture}
            [scale=.8,align=center, v/.style={draw,shape=circle, fill=black, minimum size=1.2mm, inner sep=0pt, outer sep=0pt},
            vr/.style={draw,shape=circle, fill=red, minimum size=1.2mm, inner sep=0pt, outer sep=0pt},
            vo/.style={draw,shape=circle, fill=orange, minimum size=1.2mm, inner sep=0pt, outer sep=0pt},
            vg/.style={draw,shape=circle, fill=ForestGreen, minimum size=1.2mm, inner sep=0pt, outer sep=0pt},
            every path/.style={},
            font=\small, label distance=1pt,
            every loop/.style={distance=1cm, label=right:}
            ]
            \draw (-3,-2) -- (-3,2) -- (3,2) -- (3,-2) -- (-3,-2);
            \draw[dashed] (-3,-1) to[out=-20,in=20,looseness=4] (-3,1);
            
            \node[v, ] (1) at (60:1.5cm) {};
            \node[v,] (2) at (0:1.5cm) {};
            \node[vg, ] (3) at (-60:1.5cm) {};
            \node[vo,] (4) at (-120:1.5cm) {};
            \node[vr,] (5) at (180:1.5cm) {};
            \node[v, ] (6) at (120:1.5cm) {};

            \path (1) edge (2);
            \path (1) edge (3);
            \path (2) edge (3);
            \path[dotted] (1) edge (4);
            \path[densely dotted] (2) edge (4);
            \path (3) edge (4);
            \path (1) edge (6);
            \path (2) edge (6);
            \path (3) edge (6);
            \path[loosely dotted] (4) edge (6);

            \path[color=cyan] (4) edge[bend right=20] (5);
            \path[color=gray] (4) edge[bend right=10] (5);
            \path[color=green] (4) edge[bend right=0] (5);
            \path[color=magenta] (4) edge[bend left=10] (5);
            \path[color=blue] (4) edge[bend right=30] (5);
    
            \draw[fill=yellow,draw=none, fill opacity=.3] (-3,-1) to[out=-20,in=20,looseness=4] (-3,1);
            \draw[fill=gray,draw=none, fill opacity=.1] (-3,-1) to[out=-20,in=20,looseness=4] (-3,1) to (-3,2) to (3,2) to (3,-2) to (-3,-2);

            \node[color=red] at (180:2.3cm) {$\LL_{n-1,n}$};
            \node[color=orange] at (240:2cm) {$\LL_{n-2,n-1}$};
            \node[color=ForestGreen] at (300:2cm) {$\LL_{n-3,n-2}$};

            \node at (190:5cm) {$V_{0i,(n-1)n}$};
            \path (193:4cm) edge[->] (205:1.7cm);

            \node at (0,-3) {(b)};
                    
        \end{tikzpicture}\qquad
    \begin{tikzpicture}
        [scale=.8,align=center, v/.style={draw,shape=circle, fill=black, minimum size=1.2mm, inner sep=0pt, outer sep=0pt},
        vr/.style={draw,shape=circle, fill=red, minimum size=1.2mm, inner sep=0pt, outer sep=0pt},
        vo/.style={draw,shape=circle, fill=orange, minimum size=1.2mm, inner sep=0pt, outer sep=0pt},
        vg/.style={draw,shape=circle, fill=ForestGreen, minimum size=1.2mm, inner sep=0pt, outer sep=0pt},
        every path/.style={},
        font=\small, label distance=1pt,
        every loop/.style={distance=1cm, label=right:}
        ]
        \draw (-3,-2) -- (-3,2) -- (3,2) -- (3,-2) -- (-3,-2);
        \draw[dashed] (-3,-1) to[out=-20,in=20,looseness=4] (-3,1);
            
        \node[v, ] (1) at (60:1.5cm) {};
        \node[v,] (2) at (0:1.5cm) {};
        \node[vg, ] (3) at (-60:1.5cm) {};
        \node[vo,] (4) at (-120:1.5cm) {};
        \node[vr,] (5) at (180:1.5cm) {};
        \node[v, ] (6) at (120:1.5cm) {};
        \path (1) edge (2);
        \path (1) edge (3);
        \path (2) edge (3);
        \path[dotted] (1) edge (4);
        \path[densely dotted] (2) edge (4);
        \path (3) edge (4);
        \path[color=cyan] (1) edge (5);
        \path (2) edge (5);
        \path[color=green] (3) edge (5);
        \path[color=magenta] (4) edge (5);
        \path (1) edge (6);
        \path (2) edge (6);
        \path (3) edge (6);
        \path[loosely dotted] (4) edge (6);
        \path[color=blue] (5) edge (6);
            
        \draw[fill=yellow,draw=none, fill opacity=.3] (-3,-1) to[out=-20,in=20,looseness=4] (-3,1);
        \draw[fill=gray,draw=none, fill opacity=.1] (-3,-1) to[out=-20,in=20,looseness=4] (-3,1) to (-3,2) to (3,2) to (3,-2) to (-3,-2);

        \node at (0,-3) {(c)};
        \end{tikzpicture}
    \end{minipage}
    \caption{(a) (part of the) Milnor fibre of $f_{n+1}$ and vanishing cycles (after Step $A$) for $n+1$ even. (b) The vanishing cycles restricted to $\Xnn$ after Step $A$. Dotted, the vanishing cycles $V_{0i, 0(n-1)}$. (c) The Iyama vanishing cycles on $\Xnn$.}
    \label{fig:inductivehomeomorphismeven}
\end{figure}
    
    We now describe the restrictions of the vanishing cycles $\Iyama_{I,J}$ to $\Xnn$. For $J \neq n$, these are (by assumption on $\Xn$ and up to Hamiltonian isotopy) entirely contained in $\Xn$ and are schematically depicted in Figure \ref{fig:Iyamasurface} (c). We perform Step $A$ of the algorithm described in Proposition \ref{orderIyamagenerators} on the vanishing cycles $\{V_{pq, (n-1)n}\}$: we leave $V_{0(n-1), (n-1)n}$ and $V_{0(n-2), (n-1)n}$ unchanged, and iteratively perform appropriate Dehn twists of the remaining ones. The restrictions to $\Xnn$ of the vanishing cycles $\{V_{0k,(n-1)n}, 1\leq k \leq n-2\}$ at the end of Step $A$ are depicted in Figure \ref{fig:inductivehomeomorphismeven} (a), where we can see each cycle $V_{0k, (n-1)n}$ intersecting $\Iyama_{n-1,n}$, entering $\Xnn$ and exiting it after crossing $\Iyama_{n-2,n-1}$. We drew the same configuration of such cycles restricted to $\Xnn$ in Figure \ref{fig:inductivehomeomorphismeven} (b). Step $B$ of the algorithm consists in a (left) Dehn twist of $V_{0p, (n-1)n}$ around $V_{0p, 0(n-1)}=\Iyama_{p,n-1}$, for all $1\leq p \leq n-2$. Again by assumption, the restrictions of $\{\Iyama_{p,n-1}\}$ to $\Xn \subset \Xnn$ are arcs joining the cycles $\Iyama_{p,p+1}$ and $\Iyama_{n-2,n-1}$ (dotted in Figure \ref{fig:inductivehomeomorphismeven} (b)). The final configuration of the arcs (restrictions of the vanishing cycles to $\Xnn$) is given in Figure \ref{fig:inductivehomeomorphismeven} (c).

    The case $n+1$ odd is completely analogous: by iteratively assuming that the collection of cycles $\LAMBDA'$ bounds a punctured surface $\Xn$ of genus 0, we can naturally embed this into $\Xnn \subset \Sigma_{n+1}$ and prove that $\Xnn$ is a punctured surface of genus 0, on which the restrictions of the vanishing cycles form a complete graph pattern.

    Let us now turn our attention to the complement $\Yn$ of $\Xn$ in $\SIy_n$, now for any $n$. We momentarily isotope the vanishing cycles restricted to $\Xn$ so that all the arcs intersecting each cycle in $\LAMBDA$, do so in one point (as depicted in Figure \ref{fig:Iyamasurface} (b)). We can then construct the \emph{ribbon graph} $\RR_{\Yn}$ associated to the restriction of the vanishing cycles to $\Yn$; its thickening is a surface of Euler characteristic determined by the Euler characteristic of the ribbon graph, which naturally embeds into $\Yn$. The Euler characteristic of a ribbon graph $\RR$ is:
    \begin{equation*}
        \chi(\RR)=\RR_0-\RR_1,
    \end{equation*}
    where $\RR_0$ is the number of vertices, and $\RR_1$ the number of edges of the graph; in our case, the number of vertices is the number of intersection points, of which we have $n-2$ (one for each cycle $\Iyama_{I,I+1}$), while the number of edges is equal to the total number of vanishing cycle, which is ${n-1} \choose 2$. For $\Iyama_{I,I+1}$, the corresponding edge in the ribbon graph is a loop based at the intersection point, while for any other vanishing cycle it is an arc connecting two distinct intersection points, corresponding to the restriction of the vanishing cycle to $\Yn$. The Euler characteristic of $\tilde{\RR}_{\Yn}$, the thickening of $\RR_{\Yn}$, is:
    \begin{equation*}
        \chi(\tilde{\RR}_{\Yn})=\frac{(3-n)(n-2)}{2}.
    \end{equation*}

    The surface $\tilde{\Sigma}_{n}:=\Xn \#_{\LAMBDA'}\tilde{\RR}_{\Yn}$ obtained from glueing $\Xn$ to $\tilde{\RR}_{\Yn}$ along $\LAMBDA'$  embeds into $\SIy_n$ and has Euler characteristic equal to it (see Section \ref{section:regularfibre}), hence the two surfaces are homeomorphic; more precisely, the former is a deformation retract of the latter. Note that the ribbon graph $\RR_{\Yn}$ is connected, in particular making $\Yn$ connected.
\end{proof}

This concludes the description of the Milnor fibre $\SIy_n$ of $f_n$, equipped with a collection of vanishing cycles that bound the Lefschetz thimbles corresponding to the Iyama generators of $\F(f_n)$.
 
    \subsection{The Fukaya category $\F(f_n)$}
    \subsubsection{Objects and morphism spaces}
    In Section \ref{section:vcyclesalgorithm} we gave an iterative construction of the regular fibre $\SIy_n$ of the Lefschetz fibration $f_n$, as well as the vanishing cycles $\LL_{I,J}$, $1 \leq I <J \leq n-1$, (bounding the thimbles $\DIy_{I,J}$ that are generators of the Fukaya-Seidel category) associated to the critical points. In this notation, we have the following:
    \begin{itemize}
        \item The \emph{waist Lagrangians}, i.e.\ vanishing cycles encircling the ``waist'' cylinders in Figure \ref{fig:Iyamasurface}, are $\LL_{I,I+1}$,  $I \in \{1,\dots, n-2\}$;
        \item Each vanishing cycle $\LL_{I,J}$ with $I < J-1$ corresponds to a vanishing cycle in Figure \ref{fig:Iyamasurface} entering the cylinders whose respective waist Lagrangians are $\LL_{I,I+1}$ and $\LL_{J-1,J}$. 
    \end{itemize}
    The morphism spaces between these generators are given as follows:
    \begin{itemize}
        \item For each object $L=\DIy_{I,J}$, $CF^*(L,L)$ is generated by the identity morphism;
        \item For any ordered pair $\LL_{I,J} < \LL_{I',J'}$, $CF^*(\DIy_{I,J},\DIy_{I',J'}) \cong CF^*(\LL_{I,J},\LL_{I',J'})$ is non-zero exactly whenever $I\leq I' < J \leq J'$, in which case it is one-dimensional and generated by the single intersection point between the two vanishing cycles.
    \end{itemize}
   
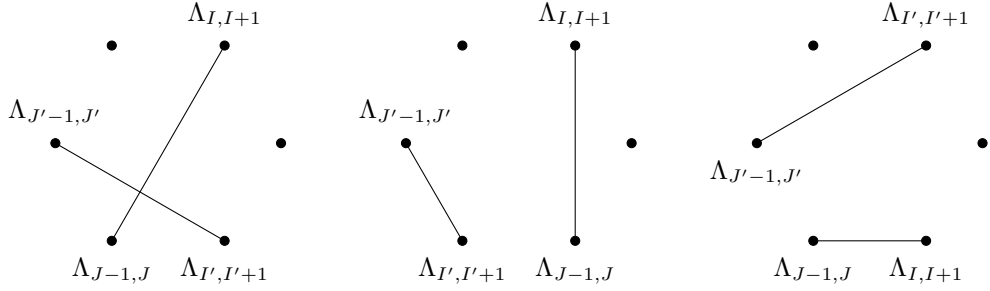
\begin{figure}
        \centering
    \begin{minipage}{.3 \textwidth}
    \centering
        \begin{tikzpicture}
        [align=center, v/.style={draw,shape=circle, fill=black, minimum size=1.2mm, inner sep=0pt, outer sep=0pt},
        every path/.style={},
        font=\small, label distance=1pt,
        every loop/.style={distance=1cm, label=right:}
        ]
        \node[v, label=90:$\LL_{I,I+1}$] (1) at (60:1.5cm) {};
        \node[v,] (2) at (0:1.5cm) {};
        \node[v, label=270:$\LL_{I',I'+1}$] (3) at (-60:1.5cm) {};
        \node[v,label=270:$\LL_{J-1,J}$] (4) at (-120:1.5cm) {};
        \node[v,label=90:$\LL_{J'-1,J'}$] (5) at (180:1.5cm) {};
        \node[v, ] (6) at (120:1.5cm) {};
        \path (1) edge (4);
        \path (3) edge (5);
        
        \end{tikzpicture}
    \end{minipage}
    \begin{minipage}{.3 \textwidth}
    \centering
    \begin{tikzpicture}
        [align=center, v/.style={draw,shape=circle, fill=black, minimum size=1.2mm, inner sep=0pt, outer sep=0pt},
        every path/.style={},
        font=\small, label distance=1pt,
        every loop/.style={distance=1cm, label=right:}
        ]
        \node[v, label=90:$\LL_{I,I+1}$] (1) at (60:1.5cm) {};
        \node[v,] (2) at (0:1.5cm) {};
        \node[v, label=270:$\LL_{J-1,J}$] (3) at (-60:1.5cm) {};
        \node[v,label=270:$\LL_{I',I'+1}$] (4) at (-120:1.5cm) {};
        \node[v,label=90:$\LL_{J'-1,J'}$] (5) at (180:1.5cm) {};
        \node[v, ] (6) at (120:1.5cm) {};
        \path (1) edge (3);
        \path (4) edge (5);
        
        \end{tikzpicture}
    \end{minipage}
    \begin{minipage}{.3 \textwidth}
    \centering
    \begin{tikzpicture}
        [align=center, v/.style={draw,shape=circle, fill=black, minimum size=1.2mm, inner sep=0pt, outer sep=0pt},
        every path/.style={},
        font=\small, label distance=1pt,
        every loop/.style={distance=1cm, label=right:}
        ]
        \node[v, label=90:$\LL_{I',I'+1}$] (1) at (60:1.5cm) {};
        \node[v,] (2) at (0:1.5cm) {};
        \node[v, label=270:$\LL_{I,I+1}$] (3) at (-60:1.5cm) {};
        \node[v,label=270:$\LL_{J-1,J}$] (4) at (-120:1.5cm) {};
        \node[v,label=27 0:$\LL_{J'-1,J'}$] (5) at (180:1.5cm) {};
        \node[v, ] (6) at (120:1.5cm) {};
        \path (1) edge (5);
        \path (3) edge (4);
        
        \end{tikzpicture}
    \end{minipage}
    \caption{(left) Intersection point generating the Floer complex $CF^*(\LL_{I,J},\LL_{I',J'})$ for $I<I'<J-1$, as seen on the restrictions of the cycles to $\Xn$. (centre) The vanishing cycles $\LL_{I,J}$ and $\LL_{I',J'}$ not intersecting on $\Xn$, for $J<I'-1$; as no intersection happens on $\Sigma_n \setminus \Xn$, the Floer complex is trivial. (right) No intersection giving rise to trivial Floer complex, for $I>I'$.}
    \label{fig:intersectioncycles}
    \end{figure}

    For $I=I'$, $J=J'$ or $I'=J-1$, the unique intersection point between the two vanishing cycles is the one in Figure \ref{fig:Iyamasurface} (b), which is seen in the cylinder whose waist Lagrangian is $\LL_{I,I+1}$, $\LL_{J-1,J}$ or $\LL_{I',I'+1}$ respectively. Figure \ref{fig:intersectioncycles} exhausts all other cases, where the intersection (or lack thereof) is seen on the connected component $\Xn\subset \SIy_n$.

\subsubsection{Grading, spin structure and $A_{\infty}$-products}\label{section:grading}
In this final section we discuss the brane structures we equip our final generating collection of Lefschetz thimbles with (existence and uniqueness of which was discussed in Section \ref{ACconf}), and we compute the non-vanishing $A_{\infty}$-products of their endomorphism algebra.

\begin{lemma}\label{mod2grading}
    There exists a $\Z_2$-grading of Floer complexes such that each of them is concentrated in degree 0.
\end{lemma}

\begin{proof}
A choice of orientation of the vanishing cycles endows transverse intersection points between them with a $\Z_2$-grading; the rule for this is given in \cite[Section~2.3]{LT}. We can choose orientations of $\Lambda_{I,J}$ as given in Figure \ref{fig:Iyamasurface}; with respect to this choice, the Floer complexes $CF^*(L_i,L_j)$ associated to ordered pairs of vanishing cycles $L_i<L_j$ lie in even degree.\end{proof}

Suppose now $L_1=\LL_{I_1,J_1} < L_2=\LL_{I_2,J_2}< L_3=\LL_{I_3,J_3}$ are three ordered vanishing cycles, with composition given by
\begin{equation}\label{composition}
    \mu^2: CF^*(L_2,L_3) \otimes CF^*(L_1,L_2) \to CF^*(L_1,L_3).
\end{equation}
As described in Section \ref{ACconf}, the coefficients of the compositions are given by the signed count of immersed triangles bounded by the (counter-clockwise) ordered union of the vanishing cycles. 

\begin{proposition}\label{prop:composition}
    The composition map (\ref{composition}) is non-zero exactly when the ordered triple $\LL_{I_1,J_1} < \LL_{I_2,J_2}< \LL_{I_3,J_3}$ satisfies the following relations:
\begin{equation}\label{inequalitiescomposition}
    I_1 \leq I_2 \leq I_3 \leq J_1 -1 \leq J_2 - 1 \leq J_3 - 1.
\end{equation}
When this holds, (\ref{composition}) is given by $y_{23}\otimes y_{12} \mapsto y_{13}$, where each $y_{ij}$ is a fixed generator of the Floer complex $CF^*(L_i,L_j)$
\end{proposition}

\begin{proof}
Given the ordered triple of vanishing cycles, the conditions (\ref{inequalitiescomposition}) are necessary in order for all the Floer complexes in (\ref{composition}) to be non-zero, and in particular for $\mu^2$ to be non-zero. On the other hand, (\ref{inequalitiescomposition}) are sufficient for each $CF^*(L_i,L_j)$ to be non-zero. When this holds, there is a single obvious holomorphic triangle contributing to the product. We can distinguish three different cases of this happening, based on whether the triangle is entirely contained in $\Xn \subset \SSn$, whether it is entirely contained in a cylinder of the surface $\SSn$, or whether it is partially contained in both. The first case is only verified when all the inequalities (\ref{inequalitiescomposition}) are strict: in this case, the triangle appears as in Figure \ref{fig:triangles} (left). The second case is verified whenever all three cycles enter the same cylinder and intersect as in Figure \ref{fig:Iyamasurface} (b): when this happens, the triangle appears as illustrated in Figure \ref{fig:triangles} (middle). This is exactly verified whenever, in addition to (\ref{inequalitiescomposition}), one of the following holds:
\begin{equation*}
    I_1=I_2=I_3, \quad J_1=J_2=J_3, \quad I_2=I_3=J_1-1, \quad I_3=J_1-1=J_2-1.
\end{equation*}

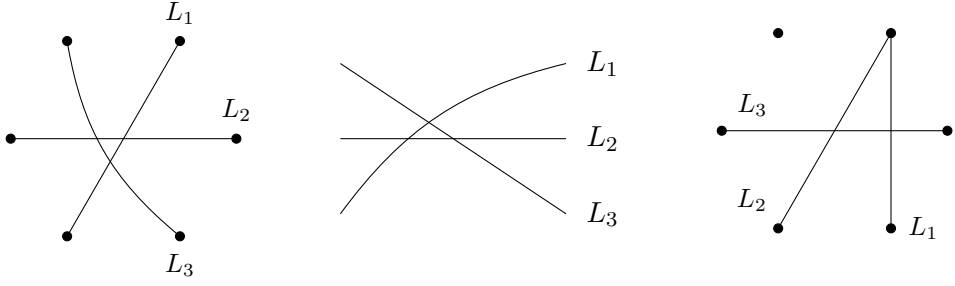
\begin{figure}
    \centering
\begin{minipage}{.3 \textwidth}
\centering
    \begin{tikzpicture}
    [align=center, v/.style={draw,shape=circle, fill=black, minimum size=1.2mm, inner sep=0pt, outer sep=0pt},
    every path/.style={},
    font=\small, label distance=1pt,
    every loop/.style={distance=1cm, label=right:}
    ]
    \node[v, label=90:$L_1$] (1) at (60:1.5cm) {};
    \node[v,label=90:$L_2$] (2) at (0:1.5cm) {};
    \node[v, label=270:$L_3$] (3) at (-60:1.5cm) {};
    \node[v,] (4) at (-120:1.5cm) {};
    \node[v] (5) at (180:1.5cm) {};
    \node[v] (6) at (120:1.5cm) {};
    \path (1) edge (4);
    \path (2) edge (5);
    \path[bend left=20] (3) edge (6);
    
    \end{tikzpicture}
\end{minipage}
\begin{minipage}{.3 \textwidth}
    \centering
\begin{tikzpicture}
    \node[label=right:$L_1$] at (3,1){};
    \node[label=right:$L_2$] at (3,0){};
    \node[label=right:$L_3$] at (3,-1){};
    \draw (0,0) to (3,0);
    \draw (0,1) to (3,-1);
    \draw[bend left=20] (0,-1) to (3,1);
    \end{tikzpicture}
\end{minipage}
\begin{minipage}{.3 \textwidth}
    \centering
        \begin{tikzpicture}
        [align=center, v/.style={draw,shape=circle, fill=black, minimum size=1.2mm, inner sep=0pt, outer sep=0pt},
        every path/.style={},
        font=\small, label distance=1pt,
        every loop/.style={distance=1cm, label=right:}
        ]
        \node[v] (1) at (60:1.5cm) {};
        \node[v] (2) at (0:1.5cm) {};
        \node[v, label=0:$L_1$] (3) at (-60:1.5cm) {};
        \node[v,label=100:$L_2$] (4) at (-120:1.5cm) {};
        \node[v,label=40:$L_3$] (5) at (180:1.5cm) {};
        \node[v] (6) at (120:1.5cm) {};
        \path (1) edge (3);
        \path (1) edge (4);
        \path (2) edge (5);
        
        \end{tikzpicture}
    \end{minipage}
\caption{The obvious triangles contributing to $\mu^2$ in (left) $\Xn \subset \SSn$, (centre) the waist region and in (right) $\Xn \subset \SSn$, where here the intersection point coinciding with the vertex of the $n$-gon represents the unique intersection point in the waist region.}
\label{fig:triangles}
\end{figure}

Finally, in all other cases satisfying (\ref{inequalitiescomposition}) the triangle appears as in Figure \ref{fig:triangles} (right), with one or more vertices at a vertex of the $n$-gon.

The triangles described above all contribute to the product (\ref{composition}). By the open mapping theorem, there is no other triangle contributing to it, so that (\ref{composition}) is given by
\begin{equation*}
    y_{23}\otimes y_{12} \mapsto \pm y_{13},
\end{equation*}
where the sign depends on the orientation of the moduli spaces of such holomorphic triangles. In order to pick the sign with which each triangle contributes, and following \cite[Section~7]{SeidelArt4}, we pick for each vanishing cycle $L_i$, of an additional point $\star_i \in L_i$ that is strictly distinct from any intersection point. This marked point endows $L_i$ with a non-trivial spin structure, which is trivialised away from it. We can make this choice so that all the points $\star_i$ are away from triangles. More precisely, if $u: D \to \SSn$ is a pseudo-holomorphic map from the 3-punctured disk, mapping each boundary component to an arc of the Lagrangians $L_i$, such that the counter-clockwise ordering is preserved, we can choose each $\star_i$ so that $u^{-1}(\star_i)$ is empty for any $i$. This choice is possible, as we know all triangles to be entirely contained in $\Sigma_n \setminus \Yn$: we can place each $\star_i \in L_i$ on the restriction of $L_i$ to $\Yn$. Following \cite[Section~7]{SeidelArt4}, we know that if $L_1$, $L_2$ and $L_3$ bounding the triangle are oriented following its natural orientation, the contribution of the latter to (\ref{composition}) is positive. Moreover, changing the orientation of $L_2$ (resp. $L_3$) changes the contribution of such triangle by a factor of $(-1)^{|y_{12}|}$ (resp. $(-1)^{|y_{23}|}$), where $|y|\in \Z_2$ denotes the mod 2 degree of the generator $y$ of the corresponding Floer complex. By Lemma \ref{mod2grading}, the grading can be chosen so that all Floer complexes lie in even degree, making the orientation of moduli space of triangles independent of the choice of orientation of vanishing cycles. Finally, having taken marked points $\star_i$ away from the boundary of triangles, these do not change the signs of their contribution to the composition, making all of these positive.
\end{proof}

\begin{corollary}\label{cor:commutativity}
    Whenever the triples $\LL_{I_1,J_1} < \LL_{I_2,J_2} < \LL_{I_3,J_3}$ and $\LL_{I_1,J_1} < \LL_{I'_2,J'_2} < \LL_{I_3,J_3}$ satisfy (\ref{inequalitiescomposition}), the compositions (\ref{commsquare}) commute.
\begin{equation}\label{commsquare}
    \begin{tikzcd}
        \LL_{I_2,J_2} \arrow{r}  & \LL_{I_3,J_3} \\
        \LL_{I_1,J_1}  \arrow{u} \arrow{r}& \LL_{I'_2,J'_2} \arrow{u}
        \end{tikzcd}    
\end{equation}
\end{corollary}

\begin{proof}
    The two compositions are given by the maps:
    \begin{equation*}
        y_{23}\otimes y_{12} \mapsto y_{13}, \quad y_{2'3}\otimes y_{12'} \mapsto y_{13}.
    \end{equation*}
\end{proof}

\begin{proposition}\label{Zgrading}
    There exists a $\Z$-grading of $\F(f_n)$ such that all Floer complexes are concentrated in degree 0.
\end{proposition}

\begin{proof}
Following the same argument made in Proposition \ref{prop:degree0}, we fix a grading of $\LL_{1,2}$ and we iteratively shift the gradings of the other vanishing cycles, so that the corresponding morphism spaces are concentrated in degree 0. By commutativity relations given by Corollary \ref{cor:commutativity}, the claim follows.\end{proof}

\begin{corollary}\label{higherprods}
    All higher $A_{\infty}$-products, except for the composition, vanish.\qed
\end{corollary}

\begin{remark}
    The above is a chain-level description of the morphism, but as the differential vanishes everywhere it is also a description on cohomology level.
\end{remark}

\bibliographystyle{alpha}
\bibliography{refs}
\end{document}